\documentclass[12pt]{article}

\usepackage{filecontents}
\usepackage{ifpdf}
\usepackage{latexsym}
\usepackage{amssymb}
\usepackage{pifont}
\usepackage{amsopn}
\usepackage{graphicx}
\usepackage{tikz}

\usepackage[empty]{fullpage}

\ifpdf
\usepackage[pdftex,pdfstartview=FitH,pdfborderstyle={/S/B/W 1}, colorlinks=true,linkcolor=blue,urlcolor=blue,citecolor=blue]{hyperref}

\makeatletter
\newcommand\fake@math{}
\def\fake@math#1\){[math]}
\makeatother

\else
\usepackage[hypertex,colorlinks,]{hyperref}
\fi

\usepackage{bookmark}

\newtheorem{theorem}{Theorem}[section]
\newtheorem{proposition}[theorem]{Proposition}
\newtheorem{lemma}[theorem]{Lemma}
\newtheorem{corollary}[theorem]{Corollary}
\newtheorem{definition}[theorem]{Definition}
\newtheorem{notation}[theorem]{Notation}
\newtheorem{claim}[theorem]{Claim}
\newtheorem{convention}[theorem]{Convention}

\newcommand{\qq}{\mathbb{Q}}

\newcommand{\zz}{\mathbb{Z}}
\newcommand{\rr}{\mathbb{R}}
\newcommand{\cc}{\mathbb{C}}

\newcommand{\nn}{\mathbb{N}}

\newcommand{\diam}{\mathrm{diam}}
\newcommand{\perim}{\mathrm{perim}}
\newcommand{\rds}{\mathrm{rds}}

\newcommand{\Mod}{\mathrm{Mod}}

\begin{document}

\pagestyle{plain}

\title{Comparison of Period Coordinates and Teichm\"uller Distances}
\author{Ian Frankel}

\maketitle

\begin{abstract}We show that when two unit area quadratic differentials are $\epsilon$-close with respect to good systems of period coordinates and lie over a compact subset $K$ of the moduli space of Riemann surfaces $\mathcal{M}_{g,n}$, then their underlying Riemann surfaces are $C\epsilon^{\alpha}$-close in the Teichm\"uller metric. Here, $\alpha$ depends only on the genus $g$ and the number of marked points, while $C$ depends on $K$.\end{abstract}

\tableofcontents

\section{Introduction and Statement of Main Result}

\subsection{Preliminaries}

\noindent The moduli space of compact Riemann surfaces of genus $g$ with $n$ unlabelled marked points (or $n$ deleted points), $\mathcal{M}_{g,n}$, is a complex variety and orbifold of dimension $3g - 3 + n$, with each point of $\mathcal{M}_{g,n}$ representing a biholomorphism class of compact genus Riemann surface with $n$ points deleted (or with $n$ marked points). The topology of $\mathcal{M}_{g,n}$ is induced by the Kobayashi metric. (The Kobayashi pseudo-metric on a complex analytic  $X$ space is the largest pseudo-metric for which all holomorphic maps of the hyperbolic disk into $X$ are non-expanding, and it is a non-degenerate metric in all cases we will consider.) Teichm\"uller constructed a metric whereby a $K$-quasiconformal homeomorphism (see section 3) between $X$ and $Y$ exists if and only if the distance between $(X;p_1,...,p_n)$ and $(Y;q_1,...,q_n)$ is at most $K$. Royden showed that this is the same as the Kobayashi metric. The orbifold universal cover of $\mathcal{M}_{g,n}$, the \emph{Teichm\"uller space} $\mathcal{T}_{g,n}$, is a contractible complex manifold of dimension $3g-3+n$, when $3g-3+n > 0$. In the same paper, Royden showed that the biholomorphism group of $\mathcal{T}_{g,n}$ is naturally isomorphic to the Teichm\"uller modular group, or mapping class group $\Mod(S_{g,n})$, of the surface $S_{g,n}$, and the quotient by this action is $\mathcal{M}_{g,n}$ \cite{Royden}.\\

\noindent Even though the Teichm\"uller-Kobayashi metric on $\mathcal{T}_{g,n}$ is not Riemannian, it has a geodesic flow, which is usually described as a dynamical system on the unit \emph{cotangent bundle}. For $X \in \mathcal{T}_{g,n}$, the cotangent space to $\mathcal{T}_{g,n}$ at $X$ consists of those meromorphic sections of the tensor square of the cotangent bundle to the complete Riemann surface corresponding to $X$ for which the only poles are simple and occur at the marked points of $X$. We refer to such sections as \emph{quadratic differentials}. We will write $QD(\mathcal{T}_{g,n})$ to denote the space of nonzero quadratic differentials on surfaces of genus $g$ with $n$ marked points.\\

\noindent Many analogies have been made between the geodesic flow for this metric and the geodesic flow on a closed negatively curved Riemannian manifold. Euclidean geometry is the main tool in the study of the Teichm\"uller geodesic flow; we describe briefly the connection here. Given a nonzero quadratic differential $q$ on $X$, there is an associated flat (Gaussian curvature = 0) metric on $X$ with a finite number of cone-type singularities, which we will call the \emph{$q$-metric}. The metric can be defined locally by taking the integral of the holomorphic $1$-form $\sqrt{q}$ to give an isometric chart to $\cc$, i.e. $$z \mapsto \int_{z_0}^z \sqrt{q}$$ gives an isometry between a neighborhood of $z_0 \in X$ with the $q$-metric and an open set in $\cc = \rr^2$ with the standard Euclidean metric, for each point $z_0$ where $q$ is holomorphic and nonvanishing.\\

\noindent Singularities of the $q$-metric occur at points where such charts cannot be defined - at zeros and poles of $q$. At these points, we have cone-type singularities with cone angle $(n+2)\pi$ at each zero of order $n \geq -1$. (If we allowed $q$ to have poles of higher order, they would be infinite distance away, so the metric cannot extend.) It should be noted that, since $\sqrt{q}$ is only defined up to sign, these charts to $\cc$ are unique up to the group $\{z \mapsto C \pm z: C \in \cc\}$. We will refer to a surface equipped with such a metric as a \emph{half-translation surface}. We note that it is also possible to take a surface with charts and transition maps lying in $\{z \mapsto C \pm z: C \in \cc\}$ and give it a metric for which those charts are isometric and recover a complex structure and a quadratic differential which is $dz \otimes dz$ in local coordinates for the system of charts. If, in addition, the resulting metric space can be completed by adding only finitely many points with cone-type singularities, the complex structure extends uniquely to the completion and the quadratic differential extends meromorphically. We also establish the following convention:\\

\begin{convention}Given $(X; p_1,...,p_n ;q)$ a Riemann surface with $n$ marked points $p_1,...,p_n$ and a nonzero meromorphic quadratic differential $q$ on $X$ whose only poles are simple and occur at a subset of $\{p_1,...p_n\}$ the set of \emph{singularities} of $(X; p_1,...,p_n ;q)$ are the zeros of $q$ on $X$ together with the collection of all marked points, regardless of whether or not $q$ has poles at those points.\end{convention}

\subsection{Period Coordinates and the Main Theorem}

\noindent We can integrate $\sqrt{q}$ along $q$-geodesic arcs $\gamma_i: [0,1] \to X$ chosen so that $\gamma_i((0,1))$ contains no singularities and $\gamma_i$ is injective on $(0,1)$, but $\gamma(0)$ and $\gamma(1)$ are (not necessarily distinct) singularities. Such arcs are called \emph{saddle connections}. A collection of such integrals forms a local holomorphic coordinate chart for $QD(\mathcal{T}_{g,n})$ in a neighborhood of $(X,q)$, provided that $q$ has $4g-4+n$ simple zeros and $n$ simple poles. We call such charts \emph{period coordinate} charts. A period coordinate chart can be extended to the boundary of any open precompact subset of $QD(\mathcal{T}_{g,n})$ on which it is well defined, and embeds into $\cc^{6g-6+2n}$ as a convex set. $QD(\mathcal{T}_{g,n})$ is locally finitely covered by the closures of such sets by Corollary \ref{CoroFinCoordSys}.\\

\noindent If a set $S \subset QD(\mathcal{T}_{g,n})$ is homeomorphic to a compact convex set with nonempty interior $K \subset \cc^{6g-6+2n}$, and a homeomorphism $f: S \to K$ whose restriction to the interior of $S$ is a period coordinate embedding, we say $S$ is a \emph{compact convex period coordinate patch}. In Appendix A we show that $QD(\mathcal{T}_{g,n})$ is locally finitely covered by compact convex period coordinate patches.\\

\noindent We will define a $Mod(S_{g,n})$-invariant path metric $d_\mathrm{Euclidean}$ (see Definition \ref{DefineEuclidean}) that has the property that every compact convex period coordinate patch is locally bi-Lipschitz to the corresponding subset of $\cc^{6g-6+2n}$.\\

\noindent We now state our main theorem:\\

\begin{theorem}\label{TheoremMain}
Let $g,n$ be nonnegative integers such that $3g - 3 + n > 0.$ Let $a_n = 1$ if $n = 0$ and $2$ if $n > 0$. Let $K$ be a compact subset of $\mathcal{T}_{g,n}$. Then for any unit area quadratic differentials $(X_1,q_1)$ and $(X_2,q_2)$ with in $X_1,X_2 \in K$, there is a constant $C_{K}$ such that $$d_{Teich}(X_1,X_2) < C_{K,s} d_{\mathrm{Euclidean}}(q_1,q_2)^{2/[2 + a_n(4g-4 + n)]}.$$
\end{theorem}

\noindent REMARK: Theorem \ref{TheoremMain} is true whether we restrict $d_\mathrm{Euclidean}$ to unit area quadratic differentials intrinsically or extrinsically (whether or not we consider paths that leave the space of unit area differentials). The reason for this is that $\|q\|$, the area of the $q$-metric, is locally Lipschitz as a function of $(X,q)$ with respect to the metric $d_\mathrm{Euclidean}$. Therefore if $(X_1,q_1)$ and $(X_2,q_2)$ are unit area and sufficiently close, the shortest path between them stays near the set of unit area quadratic differentials. A path in $QD(\mathcal{T}_{g,n})$ that starts and ends in $K$ can be projected onto the unit area subspace and if the path stays sufficiently close to $K$, this projection will only increase its length by a bounded factor.\\

\noindent The real content of the theorem is that it remains valid even when zeros of the quadratic differential are allowed to coincide with each other, or with marked points. If we required $K$ to be a compact subset of the \emph{principal stratum} of quadratic differentials, i.e. the space of quadratic differentials with $n$ simple poles and $4g-4+n$ simple zeros, then we could replace the H\"older exponent $2/[2 + a_n(4g-4 + n)]$ in Theorem \ref{TheoremMain} with $1$, by finding a common triangulation of $(X_1,q_1)$ and $(X_2,q_2)$ by saddle connections and mapping triangles to triangles.\\

\noindent We remark that no inequality in the other direction is possible, since it is possible that two quadratic differentials $q_1 \neq q_2$ may have the same underlying Riemann surface $X = X_1 = X_2$.\\

\noindent The \emph{strongly stable leaves} for the Teichm\"uller geodesic flow are the collections of unit area quadratic differentials that, on each chart given by period coordinates, can be obtained from each other by changing only the imaginary parts of period coordinates and changing systems of period coordinates. One might hope for a converse inequality if we make the additional assumption that $(X_1,q_1)$ and $(X_2,q_2)$ are on the same strongly stable leaf or strongly unstable leaf of the Teichm\"uller geodesic flow, for in such cases, $X_1,X_2$ are known to be different points in $T_{g,n}$, by the main result of \cite{HM}.\\

\subsection{Acknowledgments}

\noindent This research was conducted for the author's Ph.D. thesis at the University of Chicago. The author would like to thank Howard Masur for suggesting this problem. The author would also like to thank Howard Masur, Kasra Rafi, and especially Alex Eskin for helpful conversations.

\section{Sketch of Proof}

\noindent The proof is by construction of a quasiconformal map between nearby surfaces. The easiest way to build a quasiconformal map between $(X_1,q_1)$ and $(X_2, q_2)$ is by finding a common triangulation by saddle connections and applying a piecewise linear (PL) map. Then one can simply estimate the dilatation (the best possible constant $K$ such that the map is $K$-quasiconformal) on each triangle, and in so doing obtain a bound for the dilataion of the entire map.\\

\noindent Unfortunately, this strategy yields poor estimates unless all of the edge lengths and angles of the triangles can be bounded away from zero, and it usually can't be used to build maps between surfaces in different strata.\\

\noindent Our solution to this problem is to find a system of disks, which we call \emph{nearly regular right polygons}, or NRRP's that isolate clusters of singularities. We can triangulate the remainder of each of the surface, and build a PL map in the complement of these disks. As for the disks themselves, we need to find a quasiconformal extension of the boundary map to the NRRP, and estimate its dilatation. For this we use a Beurling-Ahlfors extension, which requires as input an estimate of the quasisymmetry of the boundary map. This is carried out section 6.\\

\noindent The boundary maps between disks are fixed since the map is PL outside of the disks. However, we need to understand the uniformization of these disks to the upper half-plane in order to estimate its dilatation. Our disks will be chosen so that, when doubled along the boundary, they have the isometry type of half-translation surfaces coming from conjugation-invariant quadratic differentials on $\hat{\cc}$. Almost all of the work consists of bounding the changes in the locations of the singularites as a function of the changes in period coordinates. Since the quadratic differentials are actually of the form $f(z)dz^2$ for some rational function $f \in \rr(z)$, it is actually easier to give a lower bound for the changes in period coordinates in terms of the changes of locations of the zeros and poles, for a sequence of perturbations tending to $0$. This is most of the work, and it is carried out in section 5.\\

\noindent The key to finding a lower bound is finding a scale on which perturbations do not have cancelling effects on period coordinates. For this we use a partial compactification of strata similar to the compactification used in \cite{5guys} to handle collisions of singularities. For us, the most important organizing information is which collections of singularities are colliding faster than others. The limit of a differential in this compactification is a stable curve with one component for each such cluster, and a meromorphic quadratic differential on each component. A colliding cluster of singularities will correspond to a component with a quadratic differential that has a higher order pole on $\infty$. The partial compactification is described in Definition \ref{DefCompactification}, and a convergent sequence is defined in Notation \ref{NotaCompactConverge}.\\

\noindent When a cluster of colliding singularities moves mostly in the same direction, we show that it essentially behaves as a point, and shadows a perturbation in a lower dimensional stratum. This is the content of Proposition \ref{PropClustermodels}. When a cluster of singularities move in what would appear to be cancelling directions according to Proposition \ref{PropClustermodels}, we find a leading order approximation for how the perturbation affects saddle connections in this cluster. If there are subclusters the leading term may again be zero, but we can break the perturbation into components, at least one of which must be detectable on some cluster. Appendix B is dedicated to showing that these kinds of perturbations do not cancel out. The end result of this argument is the key estimate, Lemma \ref{LemmaPerturbEff}, which bounds the change in the location of the singularities as a function of the change in the periods. (More precisely, we prove the contrapositive.) We finish section 5 by converting Lemma \ref{LemmaPerturbEff} to the form we need to apply to the specific disks, Lemma \ref{LemmaImplant}.

\section{Teichm\"uller Spaces and Quadratic Differentials}

\noindent In this section, we assemble the basic facts we need about Teichm\"uller spaces and coordinate systems on strata.

\subsection{Flat Metrics and Period Coordinates}

\noindent On a half-translation surface, the slope of a tangent vector at any nonsingular point is well-defined; in particular, vertical and horizontal are well-defined notions. Length and area are also well-defined in the $q$-metric, since the change of charts preserves them. It is also evident that if $X$ has the $q$-metric, then geodesics are polygonal arcs, and only change direction at cone points, and when they do, they turn by an angle of at least $\pi$ measured either way around the cone point. If $q$ is given in local coordinates by $(dx + i dy)^{\otimes 2}$, then the area of $X$ with the $q$ metric, given by $|q| = \int_X dx \wedge dy$ is the norm of $q$ in the Teichm\"uller cometric. (Note that this area form does not depend on which square root of $q$ is picked.)\\

\noindent The collection of all quadratic differentials whose associated metrics have finite area is a vector bundle over $\mathcal{T}_{g,n}$ (if we include the zero section), and the set of nonzero quadratic differentials admits a holomorphic stratification whereby each stratum consists of differentials whose metrics have the same number of singularities of each type; two singularities have the same type if they are both marked points and have the same cone angle, or they are both unmarked points and have the same cone angle. Two different strata can have the same underlying Riemann surfaces and quadratic differentials, but differ in which points are marked. The \emph{principal stratum} consists of those quadratic differentials with a simple pole at each of the $n$ marked points and $4g-4+n$ simple zeros. The complement of the principal stratum is an analytic subvariety, so the principal stratum sits inside $QD(\mathcal{T}_{g,n})$ as a dense set of full-measure (w.r.t. Lebesgue measure class).

\begin{definition} The \emph{universal half-translation surface} is a surface bundle over $T^*\mathcal{T}_{g,n}$ with fiber $S_g$, and the fiber over $(X,q)$ is a copy of $\overline{X}$ equipped with the $q$-metric.\end{definition}

\noindent We will give a more explicit description of the topology of this bundle later in this section.\\

\noindent With our convention, the set of all singularities of $q$-metrics is an \emph{incidence subvariety} of the universal half-translation surface. Locally they define a collection of sections of the universal half-translation surface, at least in a neighborhood of each point in the principal stratum. If $p_1,...,p_r$ are the singularities at a point $(X,q)$ in the principal stratum, any nearby $(X^\prime,q^\prime)$ will have singularities $p_1^\prime,...,p_r^\prime$, and each $p_i^\prime$ varies holomorphically with $(X^\prime,q^\prime)$ in a bundle over $T^*\mathcal{M}_{g,n}$ whose fibers are the Riemann surfaces represented. It is thus meaningful to talk about points being the same singularity on different surfaces. We may pick homotopy classes of arcs whose endpoints are the singularities, and the integrals $\int_{p_i^\prime}^{p_j^\prime} \sqrt{q^\prime}$ over such arcs will vary holomorphically. When the arcs are represented by saddle connections, we call the integrals $\int_{p_i^\prime}^{p_j^\prime} \sqrt{q^\prime}$ \emph{period coordinates}. In a neighborhood of any point in the principal stratum, some such collection of period coordinates associated to saddle connections will form a holomorphic coordinate chart to $\cc^{6g-6+2n}.$ On a dense subset of such coordinate charts, we can paste together Euclidean metrics to form a path metric on the moduli space of quadratic differentials, described in the next section.

\subsection{The Euclidean Metric on $QD(\mathcal{T}_{g,n})$}

\noindent A choice of norm on a vector space is equivalent to a choice of the closed unit ball, i.e. a compact, convex set with non-empty interior which is symmetric about the origin. Given a finite collection of norms we can simply take the convex hull of the unit of their unit balls to be the unit ball of the pasted metric. The resulting norm is of course equivalent (bi-Lipschitz) to any of the original Euclidean norms at each point. Call this the \emph{union convex hull operation}. It produces the smallest norm that is larger than a given set of norms. If we have an arc in a stratum, at each point of which some system of period coordinates has locally constant derivative, we describe a norm which is the speed of the arc at almost every point. For each real $L > 0$, there are locally only finitely many coordinate systems represented by saddle connections of length $\leq L$ (see Appendix A), in the following sense: for any point $(X,q) \in QD(\mathcal{T}_{g,n})$ with $q \neq 0$, there is a neighborhood of $(X,q)$ in which only finitely many such systems exist, \emph{even if $(X,q)$ is not in the principal stratum}.\\

\noindent Recall the definition of compact convex period coordinate patch in the discussion immediately before Theorem \ref{TheoremMain}. At each point $(X,q)$, we consider all good embeddings of sets containing $(X,q)$ in which $(X,q)$ maps to $\{z: |z| \leq \Theta(X,q)\}^n$, for some continuous proper function $\Theta$ which is large enough to ensure that the set of such good embeddings is not empty for any $(X,q)$. From Appendix A it is clear such a function exists. For example, $4$ times the $q$-diameter of $(X,q)$ is such a function.\\

\noindent We now make the following definition:

\begin{definition}\label{DefineEuclidean}
Fix a continuous function $\Theta: QD(\mathcal{M}_{g,n}) \to (0,\infty)$ as above. We define the \emph{Euclidean metric} $d_\mathrm{Euclidean}$ on the space of half-translation surfaces to be the path metric obtained by applying the union convex hull operation to the norms coming from compact convex period coordinate patches containing $(X,q)$, whose defining saddle connections have length less than $\Theta((X,q))$.
\end{definition}

\noindent It should be clear from the definition that $d_\mathrm{Euclidean}$ is the largest path metric for such that for all $(X,q)$, all compact convex period coordinate patches containing $(X,q)$ that take the point $(X,q)$ into $\{v: |v|_\infty < \Theta(X,q)\}^{6g-6+2n}$ are noncontracting in a neighborhood of $(X,q)$.\\

\noindent REMARK: The local Lipschitz class of the metric $d_\mathrm{Euclidean}$ depends only on the fact that the choice of it comes from a locally finite collection of period coordinates. Theorem \ref{TheoremMain} and any other statements that depend only on the local Lipschitz class this remain true if we apply define $d_\mathrm{Euclidean}$ in terms of any locally finite covering by compact convex period coordinate patches.\\

\noindent We note that the Euclidean metric extends to a path metric on the whole of $QD(\mathcal{T}_{g,n})$, but in the absence of the lower bound for the area, arbitrarily short paths with the Euclidean metric move arbitrarily far in Teichm\"uller space. Also, when we are referring to a particular stratum or stratum closure, the Euclidean metric is understood to be the \emph{intrinsic} metric on the stratum or stratum closure, rather than the Euclidean metric in the ambient moduli space.

\subsection{Teichm\"uller's Metric}

\noindent The Teichm\"uller space of Riemann surfaces of genus $g$ with $n$ marked points is the set of marked complex structures on a Riemann surface of genus $g$ with $n$ points deleted. One way to define this more precisely is, given a smooth oriented surface $S$ diffeomorphic to a surface of genus $g$ with $n$ points deleted, we look at all orientation-preserving homeomorphisms of complete finite volume hyperbolic Riemann surfaces into $S$, subject to the following equivalence relation: For $f_1: M_1 \to S$ and $f_2: M_2 \to S$, we say $(M_1,f_1)$ and $(M_2,f_2)$ are equivalent if $f_1^{-1} \circ f_2$ is homotopic to a biholomorphism. The set of equivalence classes of pairs $(M,f)$ is the \emph{Teichm\"uller} space of genus $g$ surfaces with $n$ punctures, which we will call $\mathcal{T}_{g,n}$. When it is convenient, we will sometimes view the punctures as marked points on a compact surface instead of deleted points.\\

\noindent $\mathcal{T}_{g,n}$ carries several metrics; the one we will be concerned with is the \emph{Teichm\"uller} metric. We will briefly summarize the geometric properties of this metric by characterizing its geodesics. Proofs of these can be found in, for instance, \cite{Hubbard} or \cite{FarbMarg}.

\begin{itemize}
\item Start with any nonzero holomorphic quadratic differential $\alpha$, i.e. a nonzero holomorphic section of the tensor square of the cotangent bundle of a Riemann surface. In local coordinates, this looks like $\alpha$, with $f$ holomorphic; at each point at which $f$ does not vanish we can find a holomorphic coordinate $z$ for which $\alpha = dz^2$.\\

\item In any simply connected region on such a chart, the level sets of the real and imaginary parts of $z$ give a pair of transverse \emph{measured foliations}, that is to say, foliations equipped with a transverse measure on the local leaf space, and the transverse measure is invariant under transition maps of a system of charts defining the foliation. (See \cite {FLP}, section 1.3 for a precise construction of the space of measured foliations.) Given any local coordinates $z = x + iy$, the integrals $\int_C |dx|$ and $\int_C |dy|$ are well-defined for any smooth contour of integration $C$. There is also an invariant volume form, and in the case that our surface is has punctures we will only consider those quadratic differentials for which the area form $dx \otimes dy$ assigns finite measure. This is equivalent to assuming the quadratic differential extends meromorphically to a compact Riemann surface with at most simple poles at each of finitely many new points.\\

\item Now, fix a Riemann surface $M_0$ and a pair of transverse measured foliations associated to a holomorphic quadratic differential $\alpha$ with finite area. Then, for each $\lambda \in \rr$, there is a Riemann surface $M_\lambda$, with a homeomorphism $f: M_0 \to M_\lambda$, smooth away from the zeros of $\alpha$, with holomorphic charts such that the push-forwards of $|dx|$ and $|dy|$ are $|e^\lambda dx|$ and $|e^{-\lambda} dy|$, and $\lambda \mapsto M_\lambda$ is a unit speed geodesic. The map $f$ is called a \emph{Teichm\"uller map}.\\

\item All geodesics in $\mathcal{T}_{g,n}$ arise in this way, and the length of a geodesic between two points in $\mathcal{T}_{g,n}$ is equal to the distance between those points in the Teichm\"uller metric. Any two distinct points in $\mathcal{T}_{g,n}$ lie on a unique geodesic. \end{itemize}

\noindent Suppose $f: U \to V$ is an a.e. differentiable orientation-preserving homeomorphism of bounded domains in $\cc$, with first order distributional derivatives in $L^2(U)$. Then the partial derivatives $$f_z := \frac{f_x - i f_y}{2}, f_{\bar{z}} := \frac{f_x + i f_y}{2}$$ are defined for a.e. $z$ and $|f_{\bar{z}}| < |f(z)|$ for a.e. $z \in U$.

\begin{definition}If $$K_{z_0}(f) := \left|\frac{|f_z(z_0)| + |f_{\bar{z}}(z_0)|}{|f_z(z_0)| - |f_{\bar{z}}(z_0)|}\right|$$ satisfies $K_{z_0}(f) \leq K$ for a.e. $z_0$ with respect to the Lebesgue measure class, then we say that $f$ is $K$-quasiconformal.\end{definition}

\noindent Since the quantity $K_{z_0}(f)$ is invariant under holomorphic changes of coordinates on $U$ and $V$, we say that a homeomorphism between Riemann surfaces is $K$-quasiconformal if it is $K$-quasiconformal with respect to a choice of holomorphic coordinate charts; by the above discussion the choice of charts does not matter. It is a theorem that $1$-quasiconformal maps are actually conformal.\\

\noindent The Teichm\"uller metric also has the following characterization:

\begin{theorem} (Teichm\"{u}ller) If $(Y_1,f_1)$ and $(Y_2,f_2)$ represent two points in $\mathcal{T}_{g,n}$ then their distance is at most $\frac{\log(K)}{2}$ iff there is a $K$-quasiconformal homeomorphism $g: Y_1 \to Y_2$ such that $f_2 \circ g \circ f_1^{-1}$ is homotopic to the identity on $S$. Teichm\"{u}ller maps are the unique dilatation-minimizing maps in their homotopy classes (that is, they are $K$-quasiconformal for the smallest possible $K$).\end{theorem}

\noindent With this metric, $\mathcal{T}_{g,n}$ is homeomorphic to $\rr^{6g-6+2n}$. This manifold admits a complex structure, and $QD(\mathcal{T}_{g,n})$ is a holomorphic vector bundle over $\mathcal{T}_{g,n}$ with the zero section deleted (as we have defined it). However we will not need to use this complex structure.\\

\noindent In addition to quadratic differentials that we used to characterize the Teichm\"uller metric, we will sometimes consider quadratic differentials on $\cc$ of the form $f(z)dz^2$, where $f$ is a rational function with at most one pole in $\cc$. They extend to $\hat{\cc}$ with a higher order pole at $\infty$. The induced flat metrics on $\cc$ have infinte area, and the distance to $\infty$ is infinte in such metrics.\\

\noindent We can use Teichm\"uller maps to create various bundles over Teichm\"uller space. We will pick a base point $X \in \mathcal{T}_{g,n}$. If $q$ is a quadratic differential over the Riemann surface $X$, we can associate to $q$ the Teichm\"uller map $\phi_q$ that moves $X$ along the geodesic described by $q$ by a distance equal to the area of $(X,q)$. This identifies Teichm\"uller space with the space of quadratic differentials over $X$, and we can use coordinates $(x,q)$ to refer to the point $\phi_q(x)$.\\

\noindent We can also trivialize the bundles of quadratic differentials and half-translation surface structures using Teichm\"uller maps. To do so, we need the following definition from \cite{FLP} and theorem from \cite{HM}. In the definition below, our base surface $S_{g,n}$ and all Riemann surfaces are understood to have $n$ points \emph{deleted}.

\begin{definition}Let $q_1$ and $q_2$ be vertical foliations of quadratic differentials on two Riemann surfaces $X_1, X_2 \in \mathcal{T}_{g,n}$, with vertical measured foliations $\mathcal{F}_v(q_1)$ and $\mathcal{F}_v(q_2)$, respectively. For a simple closed curve $C$ on the base surface $S_{g,n}$ with the $n$ points \emph{deleted}, let $i(\mathcal{F}_v(q_j),C)$ denote the inf of the transverse measure evaluated on curves homotopic to $C$ in the $q_j$-metric. If $i(\mathcal{F}_v(q_1),C) = i(\mathcal{F}_v(q_2),C)$ for all $C$, we say that $\mathcal{F}_v(q_1)$ and $\mathcal{F}_v(q_2)$ are equivalent and we write $\mathcal{F}_v(q_1) = \mathcal{F}_v(q_2)$.
\end{definition}

\begin{definition}The set of all equivalence classes of measured foliations is the space $\mathcal{MF}$. Its quotient by the action of $\rr^{>0}$ by scalars is \emph{PMF}.\end{definition}

\noindent The values of the numbers $i(\mathcal{F},C)$, where $C$ ranges over all (homotopy classes of) essential simple closed curves, give a weak-$*$ topology on the set of measured foliations, which give $\mathcal{MF}$ the homeomorphism type of the product of $\rr^{6g - 6 + 2n} \setminus \{ 0 \}.$ (See expos\'{e} 6 of \cite{FLP}). Moreover, if we quotient by the obvious $\rr^{>0}$ action by scalars, the result is homeomorphic to a sphere of dimension $6g-7+2n$.\\

\begin{theorem}\cite{HM} The map from $QD(\mathcal{T}_{g,n})$ to $\mathcal{T}_{g,n} \times \mathcal{MF}$ defined by $$(X,q) \mapsto (X,\mathcal{F}_v(q))$$ is a homeomorphism.\end{theorem}

\begin{theorem}\cite{HM} The map from the set of unit area quadratic differentials in $QD(\mathcal{T}_{g,n})$ to $\mathcal{T}_{g,n} \times \mathcal{PMF}$ defined by $$(X,q) \mapsto (X,\rr^{>0}\mathcal{F}_v(q))$$ is a homeomorphism.\end{theorem}

\begin{corollary} \cite{HM} Pick $(X,q) \in QD(\mathcal{T}_{g,n}).$ For each $Y$ in $\mathcal{T}_{g,n}$, there is a unique $q^\prime \in T_Y^*\mathcal{T}_{g,n}$ such that $\mathcal{F}_v(q) = \mathcal{F}_v(q^\prime)$. Moreover, given any base point $Y \in \mathcal{T}_{g,n}$, there is a trivialization of the vector bundle $T^*QD(\mathcal{T}_{g,n})$ over $QD(\mathcal{T}_{g,n})$ that sends each $(X,q)$ to $(X,q^\prime)$, where $q^\prime$ is the unique quadratic differential on $Y$ whose vertical foliation is equivalent to that of $q$. \end{corollary}

\noindent We note that the foliation of $QD(\mathcal{T}_{g,n})$ whose leaves are the sets of quadratic differentials with equivalent vertical foliations have the following property: if we fix a system of period coordinates, then moving along a leaf of this foliation only changes the imaginary parts of the period coordinates.\\

\noindent Now, we describe a trivialization of another bundle. Let $q$ be a quadratic differential over the base point surface $X$, let $q^\prime$ be a nonzero quadratic differential over $X$, and let $x \in \overline{X}$, the completion of $X$. We can use coordinates $(x,q,q^\prime)$ to refer to the point $\phi_q(x)$ on the surface $\phi_q(\overline{X})$ equipped with the metric given by the unit area quadratic differential on $\phi_q(\overline{X})$ whose vertical foliation is equivalent to $\mathcal{F}_v(q^\prime)$. The union of all such points $(x,q,q^\prime)$ is the \emph{universal half-translation surface of type} $S_{g,n}$. Similarly, we can take the universal cover of this bundle, which results in replacing each fiber with its universal cover in our coordinates (Teichm\"uller maps lift to universal covers). The coordinates are not particularly important, but the induced manifold topologies will be used for various compactness statements - for instance, we may refer to a sequence of saddle connections on quadratic differentials $(X_i,q_i)$ converging to a saddle connection on $(X,q)$ - this means that the convergence is in the Hausdorff metric on the space of compact subsets of the bundle. In the case when the surfaces in question have genus $0$ and are explicitly uniformized to $\hat{\cc}$, we have an alternative but topologically equivalent trivialization of the same bundles, so the notion of Hausdorff convergence is the same.

\begin{definition}Given a quadratic differential $(X,q)$, we may form a branched double cover as follows: Let $(\tilde{X},\tilde{q})$ be the metric completion of a metric double cover of the set of nonsingular points of $X$ given by $\{(X,v): v \mathrm{~is~a~vertical~unit~vector~w.r.t.~}q\}$. We call this the \emph{orienting double cover} of $(X,q)$, and it admits a degree $2$ map to $X$ that is branched over all points where the cone angle is of the form $(2n-1)\pi, n \in \nn$.\end{definition}

\noindent The orienting double cover is the square of a holomorphic one-form on $\tilde{X}$, and $\tilde{X}$ is connected only if $q$ is not the square of a holomorphic one-form. We can see that the orienting double cover comes with a flat metric and a unit vertical vector field at all non-singular points. Moreover, the map $i:X \to X$ defined by $(x,v) \mapsto (x,-v)$ is an involution fixing the ramification points. If $\Sigma$ is the set of singular points of $X$, and $\tilde{\Sigma}$ is the preimage of the set of singular points in the double cover. Now, on $\tilde{X}$ we can find a holomorphic $1$-form $\alpha$ such that $f^*{q} = \alpha^{\otimes 2}$, and $\alpha$ is unique up to sign.

\subsection{Triangulations and Degenerations}

\noindent Every nonzero quadratic differential admits a cellular decomposition whose 1-cells are saddle connections and whose open 2-cells are isometric to the interiors of triangles in Euclidean space. We shall refer to such decompositions as \emph{triangulations}, even though not every edge has distinct vertices.\\

\noindent There is one triangulation, due to \cite{Delaunay}, which can be constructed as follows: if we delete all the poles of a nonzero meromorphic quadratic differential, $(X,\alpha)$ we can take the metric universal cover of $X \setminus \{x_1,...x_n\}$ where the $x_i$ are the marked points; let $\Gamma$ be the deck group of this covering map. This makes sense because on $X \setminus \{x_1,...,x_n\}$ the flat metric associated to $\alpha$ is still a path metric. Then, complete the resulting space $\tilde{X}$. (The new points added will be points of cone angle $\infty$). Then, for each point in $\tilde{X}$, we consider the set of singularities and cone points of angle $\infty$ that are closest to $X$. If for some point $p$ there are 3 or more such points distance $R_p$ from $p$, then there is an embedded Euclidean disk of radius $R > 0$ around $p$. Since the length of the shortest saddle connection is nonzero, there are only finitely many singularities on the boundary of this disk, and they have a cyclic ordering. For each such $p$ draw the chords that are the boundary of the convex hull of these points in the disk, and label them as $1$-cells. We now have a $\Gamma$-invariant cell decomposition. Now quotient by $\Gamma$. The result is a cell decomposition in which the $0$-cells are singularities, the $1$-cells are saddle connections, and the $2$-cells are convex polygons, which can then be triangulated by diagonals.

\begin{definition}We refer to the above construction as the \emph{Delaunay Triangulation}. \end{definition}

\noindent The Delaunay triangulation has a number of useful properties. It is invariant under scaling the quadratic differential by nonzero complex numbers, and for a typical surface it does not involve a choice, and if it does, only finitely many choices are possible. (By a typical surface, we mean a full measure set with respect to a Lebesgue measure, described below). One of them is that the lengths of the saddle connections are bounded by twice the diameter of the flat metric, and another is that all of the angles of the triangles are bounded away from zero given an upper bound on the diameter of the surface and a lower bound on the length of the shortest saddle connection.\\

\noindent A triangulation by saddle connections, together with the collection of cone angles, imposes relations between the periods of the saddle connections. Given a triangulation, there is some number of periods that determines all of the others.\\

\noindent In \cite{ietsandmfs} the following complex manifold structure on strata in $QD(\mathcal{T}_{g,n})$ was described, using the homology of the orienting double cover:

\begin{definition}Let $(X,q)$ be a quadratic differential with double cover $(\tilde{X},\tilde{q})$, let $\Sigma$ be the set of singularities of $X$ with inverse image $\tilde{\Sigma}$, and let $\iota$ be the involution of $\tilde{X}$. Then $$H_1^{odd}(\tilde{X},\tilde{\Sigma}; \cc) := \{\varphi \in H_1(\tilde{X},\tilde{\Sigma}; \cc): \iota_*(\varphi) = -\varphi\}.$$ $$H_{odd}^1(\tilde{X}, \cc) := \{\alpha \in H_1(\tilde{X},\tilde{\Sigma}; \cc): \iota^*(\alpha) = -\alpha\}.$$ Similarly, define $H_1^{even}$ and $H_{even}^1$ to be the eigenspaces of $\iota_*$ and $\iota^*$ with eigenvalue $1$. \end{definition}

\noindent Since $\iota^2, \iota_*^2, (\iota^*)^2$ are the identity map on their domains, we have

$$H_1^{odd}(\tilde{X},\tilde{\Sigma}, \cc) \oplus H_1^{even}(\tilde{X},\tilde{\Sigma}; \cc) = H_1(\tilde{X},\tilde{\Sigma}; \cc).$$ and 
$$H_{odd}^1(\tilde{X},\tilde{\Sigma}; \cc) \oplus H_{even}^1(\tilde{X},\tilde{\Sigma}; \cc) = H_1(\tilde{X},\tilde{\Sigma}; \cc).$$

\noindent We also have the usual duality between homology and cohomology groups since we are using field coefficients.\\

\noindent Fix a subset $S$ of a stratum of $QD(\mathcal{T}_{g,n})$ sharing a triangulation $T$ by saddle connections in fixed homotopy classes relative to the singularities, and each saddle connection $\gamma$ on any $(X,q) \in S$ has two lifts $\gamma^\prime$ and $\gamma^{\prime\prime}$ to $\tilde{X}$. Fix $T$ as an oriented graph, and assume the edges have names. On the orienting double cover of a $(X,q)$, there is a triangulation $\tilde{T}$ that maps to $T$ as an oriented graph. For each $\gamma$, fix a choice of $\gamma^\prime$ and $\gamma^{\prime\prime}$. Then $\tilde{T}$ comes equipped with this choice of names as well. In other words, if we say two surfaces have the same triangulation, then an isomorphism of the two as simplicial complexes determines a name-preserving isomorphism of their orienting double covers as simplicial complexes, which is unique up to composition with the involution. We require that one of the two possible isomorphisms of the double cover preserves the edge names.

\begin{definition}A triangulation $T$ by saddle connections \emph{degenerates} to a cell decomposition $T^\prime$ if there is a sequence of unit area half-translation surfaces $\{X_i\}_{i = 1}^\infty$, all with the same marked triangulation $T$, that converge in the Hausdorff topology in the universal curve, each cell of $T$ in the sequence of surfaces $X_i$ converges in the Hausdorff topology to a cell or union of cells of a surface $X$ with the cell decomposition $T^\prime$. If $v$ is a vertex of $T$ in each $X_i$, we say that $v$ \emph{limits to} $v^\prime$. If a collection of vertices all limit to the same $v \in v^\prime$ in $X$, we say that those vertices \emph{collide}.\end{definition}

\begin{definition}Let $T$ be a triangulation of $S_g$, together with a prescribed cone angle for each vertex of $T$. Assume $T$ has $4g-4 + n$ vertices with prescribed cone angle $3\pi$ (the zeros of $T$) and $n$ vertices with cone angle $\pi$ (the poles of $T$). We say a simply connected subset $A$ of $QD(\mathcal{T}_{g,n})$ is \emph{$T$-convex} if each $(X,q) \in A$ has a triangulation by saddle connections isotopic to $T$ or a degeneration $T^\prime$ of $T$ such that

\begin{itemize}
\item The order of vanishing of $q$ at each vertex $v$ of $T^\prime$ is the sum of the orders of vanishing of all vertices $v_k$ of a sequence $(X_i,q_i)$ triangulated by $T$ and degenerating to $v$ in $(X,q)$. (Given a maximal collection $D$ of vertices joined by degenerate edges of $T$ whose union is connected and contains only nullhomotopic closed curves, the order of vanishing of $q$ at the degenerate singularity is equal to the number of zeros minus the number of poles in $D$.)
\item For each oriented edge $\gamma$ of $T$, there is a half-plane $H_\gamma \subset \cc$ not containing $0$ in its interior, such that if we vary $(X,q)$ continuously, we always have $$\int_{\gamma^\prime}\sqrt{\tilde{q}} \in H_\gamma.$$
\item For $(X,q)$ let $v(X,q)$ be the vector of periods of saddle connections of $(X,q)$ that belong to $T$, such that each edge $\gamma$ the coordinate is in $H_\gamma$. Then for any $t \in (0,1)$ $(X_1,q_1),(X_2,q_2) \in A$ there is some $(X_3,q_3) \in A$ such that $$v(X_3,q_3) = tv(X_1,q_1) + (1-t)v(X_2,q_2).$$
\end{itemize}\end{definition}

\noindent We will show in Appendix A that there is a locally finite collection of triangulations $T_i$ on $QD(\mathcal{T}_{g,n})$ such that $QD(\mathcal{T}_{g,n})$ is a locally finite union of $T_i$-convex sets, and the sets are invariant under scaling by $\rr^{>0}$. The purpose of this definition is to allow us to analyze collisions of singularities while letting them retain their individual identities.\\

\noindent A $T$-convex subset of a stratum admits an embedding into $H_{odd}^1(\tilde{X},\tilde{\Sigma}; \cc)$ by integrating a choice of square root of $\tilde{q}$ along relative cycles in $H_{odd}^1(\tilde{X},\tilde{\Sigma}; \cc).$ We always choose the choice of square root so that 

$$\int_{\gamma^\prime}\sqrt{\tilde{q}} \in H_\gamma.$$

\noindent The choice of triangulation gives a local trivialization of the vector space $H_{odd}^1(\tilde{X},\tilde{\Sigma};\cc)$, so this gives us a system of charts into $\cc^n$, where $n$ is the dimension of $H_{odd}^1(\tilde{X},\tilde{\Sigma};\cc)$; in fact they form an atlas of charts, and the dimension of the connected component of a stratum is the dimension of $H_{odd}^1(\tilde{X},\tilde{\Sigma};\cc)$.\\

\noindent There is a canonical identification $$H_{odd}^1(\tilde{X},\tilde{\Sigma};\cc) = \mathrm{Hom}(H_1^{odd}(\tilde{X},\tilde{\Sigma};\zz),\cc)$$ which gives rise to a natural Lebesgue measure on each stratum, by picking a basis for the free $\zz$-module $H_1^{odd}(\tilde{X},\tilde{\Sigma};\zz)$, we can write any quadratic differential in terms of this basis as an element of $\cc^d$ for some $d$ by taking the integral over this basis, and the natural Lebesgue measure on $\cc^d = \rr^{2d}$ is invariant under the adjoint action by $\mathrm{GL}_n(\zz)$. We refer to this measure on any stratum of $QD(\mathcal{T}_{g,n})$ as \emph{Masur-Veech} measure. We can obtain a measure on the space of unit area differentials by a standard cone construction; the coned measure of a set of unit area differentials is defined to be the measure of the union of their multiples by scalars in $(0,1)$.\\

\noindent It is well know that these measures are finite on the whole stratum, invariant under the Teichm\"uller geodesic flow. They were constructed by Masur and Veech in \cite{ietsandmfs} and \cite{Ve86} in order to apply ergodic theory to the study of the Teichm\"uller flow. In this note it is of use to us that the condition of having a vertical or horizontal saddle connection is Masur-Veech measure $0$. If $\mu$ and $\nu$ are strata and $\nu \subset \bar{\mu}$, then $\nu$ is nowhere dense and measure $0$ with respect to Masur-Veech measure on $\bar{\mu}$, in the sense that there are neighborhoods of $\nu$ whose intersections with $\mu$ have arbitrarily small measure.

\section{Quadratic Differentials on the Sphere}

\noindent In this section, we collect basic facts about collisions of singularities in quadratic differentials. We show that every cluster of singularities is biholomorphic to a cluster of singularities on a meromorphic quadratic differential on $\hat{\cc}$, and describe a partial compactification of strata in terms of this uniformization.

\subsection{Singularities and $\delta$-Clusters}

\begin{notation}
Let $A,B$ be two positive real quantities that depend on a common variable. We write $A \dot{\prec} B$ or $B \dot{\succ} A$ if there is a positive constant $c$ such that $A \leq cB$. We write $A \dot{\asymp} B$ if $A \dot{\prec} B \dot{\prec} A$.
\end{notation}

\noindent We will use the following version of Mumford's compactness criterion.

\begin{theorem}Let $K$ be a closed set in the moduli space of unit area quadratic differentials $QD^1(\mathcal{M}_{g,n})$. For $(X,q) \in K$, let $S(q)$ denote the infimum of the $q$-lengths of simple closed curves in $X$ which are not nullhomotopic or homotopic to a loop around a puncture. Then $S(q) \dot{\succ} 1$ as $(X,q)$ ranges over $K$ if and only if $K$ is compact.\end{theorem}

\noindent This is really a combination of Mumford's original theorem \cite{Mumf}, which was about lengths of the curves in hyperbolic metrics, with any of several other works that can be used to compare $q$-metrics to hyperbolic metrics. See for instance \cite{ThickThin} or \cite{Maskit}.\\

\noindent There are two concrete, equivalent descriptions of complex manifold structures on strata of $QD(\mathcal{T}_{0,n})$. Given a quadratic differential, we may assume three of its poles are at $0,1,\infty$, since it necessarily has four more poles than zeros (counted with multiplicity) and all of the poles are simple. Once this is fixed, it follows that each stratum can be represented locally as a space of differentials of the form $\frac{P(z)}{z(z-1)Q(z)}dz^2$, where the roots of $P$ have prescribed multiplicities, the roots of $Q$ do not repeat and are distinct from $0$ and $1$, and the number of roots of $P$ of each multiplicity that are also roots of $Q$ is fixed. The locations of the roots of $P$ and $Q$, together with the ratio of leading coefficients, give local coordinate charts.\\

\noindent It follows easily that, once we restrict to any fixed triangulation, the periods are defined by holomorphic functions of the coefficients on each stratum, or equivalently, by the locations of the singularities in $\cc \setminus \{0,1\}$. Periods vary holomorphically with respect to these coordinates: if the locations of singularities vary holomorphically in the plane, we can assume that the endpoints of any saddle connection remain fixed at $a$ and $b$, and one other singularity remains fixed at $\infty$ by applying a holomorphic choice of M\"obius transformation, and along a saddle connection, $\sqrt{q}$ varies in a $1$-parameter family $h_t$. We can compute the derivative as $$\frac{d}{dt}\int_{a}^b h_t(z)dz$$ which we can simply differentiate under the integral over a fixed contour of integration. In fact, even if the endpoints do not remain fixed, differentiation under the integral remains valid when the endpoints are zeros of the differential. If we never allow our singularities to come close to colliding, this is an effective way to estimate the order of magnitude of $d_\mathrm{Euclidean}$ as a function of the change of the locations of the singularities. When collections of singularities collide or come close, we can make estimates from two perspectives: one, that they never really collided (we zoom in) or two, that they stay together (we zoom out). To aid us in analyzing these clusters of singularities we make the following definition:

\begin{definition}Given $0 < \delta < 1,$ a metric space $X$, and a discrete set $S \subset X$, we say that $D \subset S$ is a \emph{$\delta$-cluster} in $S$ if $D$ contains at least two points, and the distance between any two points in $D$ is at most $\delta$ times the distance from any point in $D$ to any point in $S \setminus D$. \end{definition}

\noindent We will first show that when a subset of the singularities forms a $\delta$-cluster $D$, and the diameter of $D$ is small compared to the lengths of all nontrivial simple closed curves, there are enough saddle connections connecting singularities in $D$.

\begin{definition}Given a graph $G$ and a subset $W$ of the vertices, the \emph{induced} subgraph with respect to $W$ is the maximal subgraph whose vertex set is $W$.\end{definition}

\noindent The Delaunay triangulation has the following property:

\begin{lemma}\label{LemmaClusterSaddles} Suppose $(X,(x_1,...,x_n),q) \in QD(\mathcal{T}_{g,n})$ or $QD(\mathcal{M}_{g,n})$. Let $S(q)$ denote the infimum length of simple closed curves on $X \setminus (x_1,...,x_n)$ that are not homotopic to punctures or constant maps, where lengths are taken in the $q$-metric, and let $\Delta(q)$ be the diameter of $X$ in the $q$-metric. Then for each $\epsilon > 0$ there is a number $\delta_0 > 0$ such that whenever $\frac{S(q)}{\Delta(q)} > \epsilon$, and $\delta < \delta_0$, for every $\delta$-cluster $D$ in the singularity set of $q$, the induced subgraph on the vertex set $D$ from the Delaunay triangulation is a connected graph. Moreover, $D$ contains at most one marked point and all cycles in this induced subgraph are nullhomotopic in $X$.\end{lemma}

\noindent Proof: We may scale the metric by a real number so that $\Delta(q) = 1$, and we may pick $\delta$ to be less than $S(q)/8$. We may also assume $\delta < 1/8.$ We note that the distance between two distinct poles is at least half of $S(q)$, so no $\delta$-cluster contains 2 distinct marked points.\\

\noindent Now suppose that $T_1$ and $T_2$ are disjoint subsets of $T$ with no edges between them in the Delaunay triangulation, and $T_1 \cup T_2  = T$. We may consider the completion of the metric universal cover of $X \setminus \{x_1,...,x_n\}$. Let $v$ be a saddle connection connecting points that project into $T_1$ and $T_2$ whose length is minimal. Indeed, the shortest path from $T_1$ to $T_2$ must be a single saddle connection, since if it were a sequence of more than one saddle connection, any intermediate singularity would be closer to the start and end than the start and end are to each other, and must therefore belong to any $\delta$-cluster containing the start and end. Pick a lift of this saddle connection and, starting at its midpoint, move along the locus of points equidistant from its ends until a point is reached that is equidistant from a third singularity. It may be the case that we do not have to move at all. We eventually meet a point $p$ which is equidistant from our two original singularities plus some more, which is the center of a 2-cell in the decomposition into convex polygons, a triangulation of which forms the Delaunay triangulation. Moving around the boundary of the circle in one direction or the other, it must be the case that all singularities we hit are at least as close to one of the ends of the saddle connection we started on the midpoint of as to each other, and therefore they are edges in the Delaunay triangulation (no matter how it is chosen). Thus they are in the original $\delta$-cluster, and the sequence of chords connecting them projects down to sequence of edges in the Delaunay triangulation that connect $T_1$ and $T_2$, all of which have both ends in $T$. $\Box$

\begin{notation}Let $(X,q)$ be a quadratic differential, and let $\gamma$ be an (oriented) arc such that the endpoints of $\gamma$ are the only singularities of the $q$-metric contained in $\gamma$. Let $\tilde{\gamma}^1,\tilde{\gamma}^2$ be the two oriented lifts of $\gamma$ in the orienting double cover $(\tilde{X},\tilde{q}).$\end{notation}

\begin{lemma}Given a quadratic differential $(X,q) \in QD(\mathcal{T}_{0,n})$ with a simple pole at $\infty$, and a collection of saddle connections forming a tree $B$ whose vertex set is the set of singularities of $(X,q)$ that are not $\infty$. Let $\gamma_1,...,\gamma_r$ be oriented saddle connections in $X$, and let $\gamma_i^{\prime\prime},\gamma_i^\prime$ be the lifts of $\gamma_i$ in the orienting double cover of $(X,q)$. Then $\{[\gamma_i^{2}] - [\gamma_i^1]: 1 \leq i \leq r\}$ is a basis for $H_1^{odd}(\tilde{X},\tilde{\Sigma}; \cc)$. \end{lemma}

\noindent Proof: In $\tilde{X}$, the inverse image of $B$ is a graph $\tilde{B}$, and $\tilde{X} \setminus \tilde{B}$ is a ramified double cover of $\hat{\cc} \setminus B$, whose only branching is of order 2 over $\infty$; call the ramification point $\tilde{\infty}$. This is biholomorphic to a disk. Thus $\tilde{B}$ is the $1$-skeleton of a $CW$ decomposition of $\tilde{X}$ with a unique 2-cell, which we may call $U$. The boundary of $U$ is zero (which can be seen since $H^2(X,\zz) \neq 0$ and there is only one 2-cell, or alternatively because its boundary wraps twice around the graph, on opposite sheets and hence $\partial U$ traverses each 1-cell twice in each direction). So $\tilde{B} \hookrightarrow X$ induces isomorphisms $H_i(\tilde{B};\cc) \to H_i(\tilde{X};\cc)$ for $i \leq 1$. Since $\tilde{B}$ is invariant with respect to the involution, the isomorphism on the first homology respects the even-odd decomposition. So all absolute cycles in $H_1^{odd}(\tilde{X},\tilde{\Sigma};\cc)$ are homologous to cycles in $H_1^{odd}(\tilde{B},\tilde{\Sigma} \setminus \{\tilde{\infty}\};\cc)$. It is clear that all of the absolute cycles are in the odd part, for both $\tilde{B}$ and $\tilde{X}$, since $\hat{\cc}$ and $B$ have no homology in dimension $1$.\\

\noindent From the long exact sequences of the pairs $(\tilde{X},\tilde{\Sigma})$ and $(\tilde{B},\tilde{\Sigma} \setminus \{\tilde{\infty}\})$ we get natural inclusions $$H_1(\tilde{X};\cc) \hookrightarrow H_1(\tilde{X},\tilde{\Sigma};\cc) \mathrm{~and~}$$ $$H_1(\tilde{B};\cc) \hookrightarrow H_1(\tilde{B},\Sigma \setminus \{\tilde{\infty}\};\cc).$$ An element of either quotient group is exactly determined by the image of any cycle representing it under the connecting map into $H_0(\tilde{\Sigma})$ or $H_0(\tilde{\Sigma} \setminus \{\tilde{\infty}\})$. If such an element is odd, it must be the case that its boundary is a linear combination terms of the form $\{p_i - \iota(p_i)\}, p_i \in \Sigma.$ Since $\tilde{\infty}$ is invariant under $\iota$ the exact same relative cycle classes are realizable as well. $\Box$

\begin{corollary}The periods of $B$ form a holomorphic coordinate chart. $\Box$\end{corollary}

\noindent A similar argument shows that, when $\tilde{\cdot}$ denotes the preimage of $\cdot$ in the orienting double cover,

\begin{corollary}\label{CoroAnybasis} If $\gamma$ is a simple closed curve containing no singularities, which bounds a disk $\mathbb{D} \subset  X$, and the singularities in $D$ form a set $S$, then the periods of any tree with vertex set $S$ and all edges saddle connections contained in $\mathbb{D}$ form a basis for $H^1(\tilde{\mathbb{D}},\tilde{S}; \cc)$ as a vector space over $\cc$.\end{corollary}

\begin{lemma}Let $(X,q) \in QD(\mathcal{T}_{0,n})$. If $\Sigma$ is its singularity set and $p$ is a pole in $\Sigma$, then any Delaunay triangulation has a connected induced subgraph with respect to $\Sigma \setminus \{p\}.$\end{lemma}

\noindent REMARK: This property is shared with the $L^\infty$ Delaunay triangulations considered in the appendix.\\

\noindent Proof: There is only one ray of each slope emanating from $p$ so there cannot be a saddle connection that starts and ends at $p$. The edges emanating from $p$ are cyclically ordered counterclockwise, and there is an edge connecting each consecutive pair of them, so for any pair of neighbors of $p$ there is a path in the graph avoiding $p$. Therefore, $p$ is not a cut vertex of the $1$-skeleton. $\Box$\\

\noindent On compact sets, any two holomorphic coordinate systems are bi-Lipschitz, so the main difficulty will be in dealing with degenerations to lower-dimensional strata. Our main object will be to examine periods of clusters of singularities close to $0$.\\

\noindent In the sequel, we will often have to deal with a meromorphic quadratic differential $q$ defined on $\cc$, and we will often have to speak of distances between points, diameters of sets, etc. in two different metrics: the usual Euclidean metric on $\cc = \rr^2$ (which does not depend on the choice of $q$) and the singular metric that depends on $q$.

\begin{notation}$d_\cc$ and $\diam_\cc$ shall denote the distance and diameter in the usual metric on $\cc$ and $d_q$ and $\diam_q$ shall be used to denote distance and diameter in the singular metric defined by $q$. $\perim_\cc$ and $\perim_q$ will be used to denote perimeters of regions. $\rds_\cc(A)$ and $\rds_q(A)$ will denote the maximum radius of a ball contained entirely in $A$ with respect to each of $d_\cc, d_q$. \end{notation}

\noindent We now discuss a certain class of quadratic differential on $\cc$ which we use to model clusters of singularities. These differentials will extend meromorphically to $\hat{\cc}$, but will have higher order poles at $\infty$.\\

\begin{definition}We say $q$ is a \emph{cluster differential} if it satisfies the following:

$q$ is a quadratic differential on $X = \cc$ which takes one of the following forms: $q = p(z)dz^2$, where $p$ is a monic polynomial of degree at least $2$ whose roots sum to $0$ or $q = \frac{p(z)}{z}dz^2$. If $p$ is of the second type, we say that $0$ is a \emph{marked point}.\end{definition}

\noindent We consider two cluster differentials to be distinct if one has a marked point and the other does not, even though they may have the same underlying quadratic differential on $\cc$. That is to say, $0$ is \emph{allowed} to be a root of $p$.\\

\noindent The set of cluster differentials with a fixed degree and number of marked points are complex manifolds, with the coefficients of the polynomial as coordinates. The complex dimension of the space of cluster differentials is $\deg(p) - 1$ if the differentials are of the form $p(z)$ and there is no marked point, and $\deg(p)$ if there is a marked point.\\

\noindent Spaces of cluster differentials are also stratified by the numbers and types of singularities, and the periods of any spanning tree of saddle connections forms a local holomorphic coordinate system.\\

\noindent In fact, the periods of the saddle connections in the tree and the cone angles determine the $q$-metric entirely, since the metric has constant curvature $0$ away from the tree. This determines $(\cc,q)$ as a metric space, from which the conformal structure can be recovered; this determines the differential and the locations of its singularities up to a M\"obius transformation fixing $\infty$. Since we know that $p$ is monic and we know which direction is vertical, we can recover $p$ up to multiplying all of its singularities by a root of unity. If the singularities have names then locally there is a unique correct choice.) We may write $q$ as $q(z_1,...,z_n,w)$ where $n$ is the number of zeros of $q$ and $w$ is the marked point. We are of course constrained to the hyperplane $\sum(z_i) = 0$ if there is no marked point.\\

\noindent If a cluster differential has no poles, then its infinitesimal metric is locally nonpositively curved in the sense of Alexandrov. If a quadratic differential induces a complete Alexandrov nonpositively curved metric on a simply connected Riemann surface, then there is a unique length-minimizing geodesic between every pair of points. In particular, a cluster differential has only finitely many saddle connections if it has no poles. We also see that if we take a branched double cover of a cluster differential branched only over the pole, and pull back the differential, we get a cluster differential with no poles and a marked point. All saddle connections on our original cluster differential lift to one of finitely many saddle connections upstairs so there were only finitely many saddle connections on the original cluster differential, and the number of saddle connections is bounded by the number of pairs of singularities in the double cover.\\

\noindent There is an analogous definition of $T$-convex for cluster differentials - we may say that a set $A$ is $\Gamma$-\emph{convex} if $\Gamma$ is the graph whose vertices are singularities and whose edges are saddle connections, the periods of each saddle connection that is an edge of $\Gamma$ vary continuously in a half plane, and we can take convex combinations in period coordinates without leaving $A$.

\begin{proposition}\label{PropScale}The effect of multiplying the zeros and poles of a cluster differential $q$ by $t$ is to scale all periods by $t^{(m+2)/2}$ where $m$ is the rational function degree of $q$, and angles between saddle connections meeting at a singularity are preserved.\end{proposition}

\noindent Proof: We give the proof for $q = p(z)dz^2$ without the marked point since the proof with a marked point is identical. The first statement is a straightforward change of variables: $$\int_{ta}^{tb} [t^mp(z/t)]^{1/2} dz = \int_a^b t^{m/2}[p(z)]^{1/2} t dz.$$

\noindent The left hand side and the right are the values of the two periods in question, which clearly differ by $t^{(m+2)/2}$. The only way that the cone angles can differ is by integral multiples of $\pi$, but they don't differ at all since we can vary $t$ continuously and cone angles vary continuously with $t$.\\

\noindent A similar change of variables tells us the following:

\begin{proposition}\label{PropDerivScale}Fix a stratum $\mathcal{Q}$ of quadratic differential with $n$ singularities in $\cc$, and let $m$ be the total number of zeros minus number of poles (counted with multiplicity). For an $n$-tuple $p$ of distinct points in $\cc$ which can be the singularity set of an element in $\mathcal{Q}$, let $q_p$ be the cluster differential with those singularities (it is assumed that for each coordinate, the singularity has a prescribed type). Then for any  $v \in \cc^n$ and any period $P_i$ of any saddle connection we have $$\frac{d}{dt}P_i(q_{sp + tv})|_{t = 0} = s^{m/2}\frac{d}{dt}P_i(q_{p + tv})|_{t = 0}.$$\end{proposition}

\begin{corollary}\label{CoroClusterSizes}Let $q$ be a cluster differential of rational function degree $m$. Let $\Sigma$ be the set of singularities of $q$. Then $$\diam_q(\Sigma) \dot{\asymp} \diam_\cc(\Sigma)^{(2 + m)/2}.$$ \end{corollary}

\noindent Proof: By Proposition \ref{PropScale} this reduces to the following claim: If $\diam_{\cc}(\Sigma) = 1,$ then $\diam_{q}(\Sigma) \dot{\asymp} 1$.\\

\noindent We will only deal with the case in which there is a marked point since that is more difficult. What we will do is prove an upper bound on the $q$-diameter of the $d_\cc$ ball of radius $1$ about the marked point and a lower bound on the $q$-length of an arc whose endpoints belong to $\Sigma$ and are distance $1$ apart with respect to $d_\cc$.\\

\noindent Indeed, suppose that the marked point is $0$. We will give an upper bound for the $q$-length of any radius of the unit circle $|z| = 1$, since any two points in $\Sigma$ can be joined by a path contained in a pair of two such radii. By rotating our coordinate system, we may assume that the radius is the interval $[0,1] \subset \rr$.\\

\noindent If $q = \frac{f(z)}{z}dz^2$ and $f$ has degreee $m+1$ then for all points $x \in [0,1]$ we have $|f(x)/x| \leq (x+1)^n/x$. Therefore, the $q$-length of the line segment $[0,1]$ is bounded by

$$d_q(0,1) \leq \int_0^1 x^{-1/2} dx.$$ This gives us the upper bound.\\

\noindent For a lower bound, note that a rectifiable arc $\gamma$ from $a$ to $b$ with $d_\cc(a,b) = 1$ must travel Euclidean distance at least $1/2$ outside of the following union of $m+2$ disks: $$\cup_{z_i \in \Sigma}\left\{z: |z - z_i| < \frac{1}{4(m+2)}\right\}.$$

\noindent As before, assume that $q = f(z)/z dz^2$ and $z_1,...,z_{m+1}$ are the roots of $f$. we have $$|f(z)/z| \leq \left|\frac{z - z_{m+1}}{z}\right| \prod_{j = 1}^m |z - z_j| \geq \frac{|z - z_{m+1}|}{1 + |z - z_{m+1}|}\cdot \prod_{j = 1}^m |z - z_j| \geq \frac{1/(4m+8)}{1 + 1/(4m+8)}\cdot \frac{1}{(4m+8)^m}.$$

\noindent The length of $\gamma$ is therefore bounded below by $$\int_\gamma |\sqrt{q}| \geq \frac{1}{2}\left| \frac{1/(4m+8)}{1 + 1/(4m+8)} \cdot \frac{1}{(4m+8)^m}\right|^{1/2}. \Box$$

\begin{proposition}\label{PropPolyScale}Let $f$ be a monic polynomial of degree $d$ with $f(z) \neq 0$ whenever $|z| < 1$, and assume $r < 1$. Then there is a constant $C_{r,d}$ depending on $r$ and $d$ but not on $f$ such that on the disk $\{|z < r|\}$, we have $C_{r,d}^{-1} \leq |f(z)/f(0)| \leq C_{r,d}$. Moreover, if $r$ is small enough we may take $C$ arbitrarily close to $1$.\end{proposition}

\noindent Proof: For a monic polynomial $z - w$, this follows from the fact that whenever $|w| \geq 1$ and $|z_1|,|z_2| < r,$ we have $$\left|\frac{z_1-w}{z_2-w} - 1\right| \leq \frac{|z_1| + |z_2|}{|w| - r} \leq \frac{2r}{1-r}. \Box$$

\begin{notation}For the purpose of Corollary \ref{CoroBallsizes}, let $B_s$ denote the disk $\{z: z < s\}$.\end{notation}

\begin{corollary}\label{CoroBallsizes} Let $r < 1/4$. Let $p(z)dz^2$ be a cluster differential with a pole of order $m+4$ at $\infty$ and no singularities in $\cc$ outside $B_r$. Let $f$ and $g$ be polynomials with fixed degrees and no zeros in $B_1$. $q = \frac{f(z)}{g(z)}p(z)dz^2$.

Then $$d_q(0,\partial B_r) \dot{\asymp} \left|\frac{f(0)}{g(0)}\right|^{1/2}r^{(2+m)/2}, \rds_q(B_r) \dot{\asymp} \left|\frac{f(0)}{g(0)}\right|^{1/2}r^{(2+m)/2},$$ $$\diam_q(B_r) \dot{\asymp} \left|\frac{f(0)}{g(0)}\right|^{1/2}r^{(2+m)/2}\mathrm{~and~}\perim_q(B_r) \dot{\asymp} \left|\frac{f(0)}{g(0)}\right|^{1/2}r^{(2+m)/2}.$$

The implied constants depend on $m$ and the degrees of $f$ and $g$.\end{corollary}

\noindent Proof: For the base case where $f,g$ are both constant, we can apply Corollary \ref{CoroClusterSizes}. For other cases, we may simply apply Proposition \ref{PropPolyScale}. $\Box$\\

\begin{corollary}\label{CoroDontleave} Let $q$ be as in Corollary \ref{CoroBallsizes}. Then, for each $C \in (0,1)$ there is a constant $c > 0$ such that if all zeros and poles of $p$ are contained in $B_{cR}$, then any length minimizing path in the $d_q$-metric between points in $B_{cr}$ is contained in $B_{Cr}$.\end{corollary}

\noindent Proof: By the estimates in Corollary \ref{CoroBallsizes}, if $c/C$ is small enough, then the diameter of $B_{cr}$ is less than the $d_q$ distance from $\partial B_{cr}$ to $\partial B_{Cr}. \Box$\\

\noindent The following lemma says that while a faraway singularity may change the size of a $\delta$-cluster, it does little to change the shape:\\

\begin{lemma}\label{LemOutsiders} Let $f,g,B_s$ be notations as in Corollary \ref{CoroBallsizes}. Let $q = p(z)dz^2$ be a cluster differential with all singularities in $B_r$, and let $q^{\prime} = \frac{f(z)}{g(z)}p(z)dz^2$. Then for each $\epsilon > 0$ there is some $r_0 > 0$ such that whenever $r < r_0$ and $a,b,c,d \in B_r$, $$1 - \epsilon < \frac{d_q(a,b)d_{q^\prime}(c,d)}{d_q(c,d)d_{q^\prime}(a,b)} < 1+\epsilon.$$ Moreover, if $\Gamma$ is a tree whose vertices are the singularities of $p$ and whose edges are saddle connections whose lengths with respect to $d_q$ are as small as possible, then there is a tree $\Gamma^\prime$ on whose endpoints are saddle connections of $q^\prime$ and whose edges are homotopic (rel endpoints) to those of $\Gamma$ whenever $r < r_0$ is sufficently small, and angles between corresponding edges differ by less than $\epsilon$. \end{lemma}

\noindent Proof: The first claim is immediate from the pointwise estimates Proposition \ref{PropPolyScale}, since the metric is scaled pointwise by a near constant.\\

\noindent For the second claim, we need to show that if we start with our embedding of the graph $\Gamma \subset \cc$ and modify it by making its edges to become geodesic with respect to $d_{q^\prime}$, the edges remain saddle connections. $\Gamma$ is obtained by a greedy algorithm: pick saddle connections one at a time, repeatedly picking the shortest saddle connection that does not join two vertices belonging to the same connected component of the graph with edges already picked. (This is because for any saddle connection $e$ not used, the fundamental cycle of $\Gamma \cup e$ with $e$ must consist of edges no longer than $e$.) Therefore, whenever an edge $e$ is drawn, it must be the case that there is no singularity of distance less than the length of $e$ away from both ends of $e$, and there are two equilateral triangles sharing the edge $e$ with no singularities in the interior of either triangle. Also there are no other singularities on $e$. It follows that the angle formed by $e$ and any previously drawn edge is at least $\pi/3$.\\

\noindent The estimate of Proposition \ref{PropPolyScale} then implies that the lengths and angles of all tangent vectors are scaled by nearly the same constant, so for any points $a,b$ joined by a $q$-geodesic with no singularities, we have $$\frac{\int_a^b \sqrt{q^\prime}}{\int_a^b\sqrt{q}} = \sqrt{f(0)/g(0)}[1 + o(1)] ~ \mathrm{as} ~ r_0 \to 0.$$

\noindent We may scale the entire differential by a constant without changing the geodesics, so we may assume $\sqrt{f(0)/g(0)} = 1$. Now, if we have a map from a region bounded by the union of two equilateral triangles into $\rr^2$ sharing an edge, and its derivative is close enough to $\left( \begin{array}{rl} 1 & 0 \\ 0 & 1 \end{array} \right)$, then the geodesic in $\rr^2$ joining the endpoints is contained in the image of the region. The result follows.$\Box$

\begin{proposition}Let $f,g$ be polynomials such that $g$ has no repeated roots, and $\deg(f) \geq \deg(g) - 3$. Let $q =\frac{f(z)}{g(z)}dz^2$ be a meromorphic quadratic differential on $\hat{\cc}$ whose poles are all simple, except perhaps for a higher order pole at $\infty$. Assume that the roots of $f$ and $g$ are all contained in $B_{r_0}$ for some $r_0 > 0$. Then for each $r > 0$ there is a number $R$ depending only on $r, r_0$ and the degrees of $f$ and $g$, such that the shortest path between any points in $B_r$ with respect to the metric $d_q$ is contained in $B_R$.\end{proposition}

\noindent Proof: We will assume $f,g$ are monic, since scaling the metric by a constant doesn't change whether or not a path is length-minimizing. Let $m = \deg(f) - \deg(g)$. As $|z| \to \infty$ we have $|\frac{f(z)}{g(z)}|^{1/2}/|z^{m/2}| \to 1$ uniformly for all $q$ satisfying our hypotheses. Integrating this pointwise bound gives us the following:\\

\noindent The distance from the circle of radius $r$ to the circle of radius $r^2$ therefore has asymptotics $$\int\limits_r^{r^2} |\sqrt{q}| \sim \int\limits_r^{r^2} |z^{m/2}dz| \sim \left\{ \begin{array}{rl}
2 r^{-1/2} & \mathrm{~if~} m = -3 \\
\log(r) & \mathrm{~if~} m = -2 \\
\frac{2r^{m+2}}{m+2} & \mathrm{~if~} m \geq -1
\end{array} \right. $$

\noindent As $r \to \infty$, the ratio between the actual $q$-distance between $\partial B_r$ and $\partial B_{r^2}$ and the asymptotic value converges to $1$ at a rate that depends only on $r_0$, and the degrees of $f$ and $g$.\\

\noindent The $q$-lengths of semicircles on $\partial B_r$ are asymptotic to $\pi r\cdot r^{m/2}$. Thus, the distance between two points on $\partial B_r$ is at most $[1 + o(1)]\pi/4$ times the distance from $\partial B_r$ to $\partial B_{r^2}$ and back when $m = -3$, and $o(1)$ times the distance from $\partial B_r$ to $\partial B_{r^2}$ and back for $m \geq -3$. $\Box$\\

\noindent We could have used $B_{cr}$ in place of $B_{r^2}$ and gotten the same conclusion in the case $m = -3$, for any sufficiently large $C$. By performing the change of coordinates $z \mapsto 1/z$, $c = 1/C$, we conclude the following:

\begin{proposition} If a quadratic differential $q$ on $\hat{\cc}$ has only simple poles, one occurring at $z = 0$, and no other poles or zeros in the ball $B_r,$ then there is some $c \in (0,1)$ depending only the number of zeros and poles of $q$ of each multiplicity, such that no length-minimizing path from two points outside of $B_r$ passes through $B_{cr}$. $\Box$\end{proposition}

\begin{proposition}\label{PropContinDiam}Fix integers $e_1,...,e_k \geq -1$ with $\sum\limits_{i = 1}^k e_i \geq -3$. Let $q(\alpha,z_1,...,z_k)$ denote the quadratic differential $\left[\prod\limits_{i=1}^k (z-z_i)^{e_i}\right]dz^2$ on $\cc$. (It will be meromorphic on the Riemann sphere with, a pole at $\infty$, which is simple if and only if $\sum e_i = -3$.) Assume the set of tuples $(\alpha,z_1,...,z_k)$ is restricted to the set where $z_i \neq z_j$ if $e_i = e_j = -1$. Then for all $i,j \in \{1,...,k\}$ $d_{q(\alpha,z_1,...,z_k)}(z_i,z_j)$ varies continuously with respect to $(\alpha,z_1,...,z_k)$.\end{proposition}

\noindent Proof: Let $\epsilon > 0$. By the estimates Corollary \ref{CoroBallsizes} the $q$-diameter, radius, perimeter, and distance from the boundary to the point $z_i$ of $$B_{r_i}(z_i) := \{z: |z-z_i| <  r_i\}$$ are all $O(r_i^{1/2})$ as $r \to 0$. We pick  around all singularities small enough so that each has diameter less than $\epsilon$, and also a ball $B_\infty$ about $\infty$ such that no length-minimizing path between $z_i$ and $z_j$ can enter. We can pick the balls much smaller than the distance between any two points $z_i,z_j$ unless $z_i = z_j$. In this case we pick $r_i = r_j$. In a neighborhood of $(\alpha,z_1.,,,z_k)$. We can fix $r_1,...,r_k$ and $\alpha$ so that each of these balls has diameter, radius, and perimeter less than $\epsilon$ in the $q$-metric in a neighborhood of $(\alpha,z_1,...,z_k)$. Then, as $(\beta,y_1,...,y_k)$ converges to $(\alpha,z_1,...,z_k)$ the sizes of the balls $B_i(z_i)$ still satisfy $$\diam_{q(\beta,y_1,...,y_k)}(B_{r_i}(z_i)) \dot{\prec} \epsilon.$$ The distances $d_{q(\beta,y_1,...,y_k)}(B_{r_i}(z_i),B_{r_j}(z_j))$ vary continuously, since the metric varies smoothly on $\hat{\cc} \setminus (B_\infty \cup \bigcup\limits_{i = 1}^k B_{r_i}(z_i)).$ Since length-minimizing paths from $B_{r_i}(z_i)$ to $B_{r_j}(z_j)$ stay in this region, it follows that $d_{q(\beta,y_1,...,y_k)}(y_i,y_j),$ viewed as a function of $(\beta,y_1,...,y_k,)$ can be written as a sum of a continuous function and a function taking values in $[-\epsilon,\epsilon].$ Since $\epsilon$ was arbitrary the proposition follows.$\Box$

\begin{definition} A \emph{M\"obius normalization} of $\mathcal{T}_{0,n}$ or $QD(\mathcal{T}_{0,n})$ is a collection of uniformization maps of the underlying Riemann surfaces to $\hat{\cc}$ such that the points mapped to $0,1,$ and $\infty$ define three continuous sections of the universal curve or universal half-translation surface, and all each of these three sections intersects each fiber at a marked point. A \emph{M\"obius normalized} collection $A$ of quadratic differentials on $\mathcal{T}_{0,n}$ is the set of quadratic differentials on $\hat{\cc}$ that pull back to elements of $A$ across the uniformization maps.\end{definition}

\begin{proposition}\label{PropJointContDist} $d_{q(\alpha,z_1,...,z_k)}(a,b)$ is jointly continuous in $(\alpha,z_1,...,z_k,a,b)$. \end{proposition}

\noindent Proof: We can just pretend $a$ and $b$ are singularities $z_{k+1}$ and $z_{k+2}$ with $e_{k+1} = e_{k+2} = 0$ and apply Proposition \ref{PropContinDiam}.$\Box$

\begin{corollary}\label{CoroEquicont} Fix a M\"obius normalization of $QD(\mathcal{T}_{0,n})$ and a compact set $K \subset QD(\mathcal{T}_{0,n})$ for which $\infty$ is always a pole. Then for any compact subset $V \subset \cc$ and any $\epsilon > 0$, there is a number $\delta > 0$ such that for all $q \in K$, $d_q(x,y) < \epsilon$ whenever $x,y \in V$ and $d_\cc(x,y) < \delta$. The same is true if we allow $K$ to vary over a compact set of cluster differentials. \end{corollary}

\noindent Proof: $\{(x,y,q): x,y \in V, d_q(x,y) \geq \epsilon\}$ is compact, so $d_\cc(x,y)$ attains a minimum there, which we can take to be $\delta$.

\begin{corollary}\label{CoroHausCont} Suppose $\{t_m\} \to t_\infty$ and $\{u_m\} \to u_\infty$ are convergent sequences in $\cc$ and $\{q_m\}$ to $q_\infty$ is a M\"obius normalized convergent sequence in $QD(\mathcal{T}_{0,n})$ or a space of cluster differentials. Then any sequence $\gamma_m: [0,1] \to \cc$ of $d_{q_m}$-geodesics from $t_m$ to $u_m$ which are length minimizing has a subsequence that converges uniformly (as $\cc$-valued functions) to a constant speed length-minimizing geodesic in $q_m$ along a subsequence.\end{corollary}

\noindent Proof: By Proposition \ref{PropJointContDist} it is clear that any subsequential limit of constant speed length-minimizing geodesics must be a constant speed length-minimizing geodesic. By a diagonalization argument, we can pass to a a subsequence such that $\gamma_m(t)$ converges for all $t \in [0,1] \cap \qq.$ Uniform convergence follows from Corollary \ref{CoroEquicont}. $\Box$\\

\begin{corollary}\label{CoroC1}Suppose that $q_n \to q_\infty$ in $QD(\mathcal{T}_{g,n})$ or a space of cluster differentials. Let $\gamma_n: [0,1] \to \cc$ be a geodesic arc for the $q_n$-metric, parametrized to be constant speed, converging to a gedoesic arc $\gamma_\infty$ in the $q_\infty$-metric that does not pass through any singularities. $\gamma_\infty: [0,1] \to \cc$ does not pass through a singularity of $q_\infty$, the sequence of maps $\gamma_n$ converge to $\gamma_\infty$ in the $C^1$ topology on $C^1([0,1])$. \end{corollary}

\noindent Proof: We already know that the lengths converge. It is easy to see that the directions converge as well, since the arcs formed by concatenating the $q_\infty$ geodesic from $\gamma_\infty(0)$ to $\gamma_n(0)$, the $\gamma_n$-geodesic from $\gamma_n(a)$ to $\gamma_n(b)$, and the $q_\infty$-geodesic from $\gamma_n(1)$ to $\gamma_\infty(1)$ are homotopic rel endpoints, via homotopies that pass through no singularities, to $\gamma_\infty$ for all sufficiently large $n$. The integrands of $\sqrt{q_n}$ converge to $\sqrt{q_\infty}$ on an open set containing all of these arcs. It thus follows that, if we take the correct branch of the square root, $$\int_{\gamma_n(0)}^{\gamma_n(1)}\sqrt{q_n} \to \int_{\gamma_\infty(0)}^{\gamma_\infty(1)}\sqrt{q_\infty}.$$

\noindent This implies that the complex lengths $\ell_n$ of the segments $\gamma_n([0,1])$ converge to some $\ell_\infty$. By Corollary \ref{CoroHausCont} and the fact that the tangent vector that maps $\sqrt{q}$ to $1$ varies continuously in $(q,z)$ away from singularities, we conclude that $\ell_i/\sqrt{q}$ the lengths and directions of the tangent vectors $\{\gamma_n^\prime(t): t \in [0,1]\}$ converge uniformly in $T\cc$.$\Box$\\

\noindent Fix a space of cluster differentials, or else fix $c > 0,$ compact subset $K$ of $QD(\mathcal{T}_{0,n})$, and assume that for each $(X,q) \in K$, we pick a marking such that $\infty$ is a pole and such that for all singularities $p$, we have $d_q(p,\infty) > c\Delta(q)$. (Recall $\Delta(q)$ is the diameter of the $q$-metric). This choice is invariant under the group of M\"obius transformations fixing $\infty$. Then we have the following:

\begin{proposition}\label{PropVanishInEither} The notions of $\delta$-cluster with respect to $\cc$-metric and $q$-metric are comparable on $K$ in the following sense:
\begin{itemize}
\item For each sufficiently small $\delta > 0$ there exists $\delta^\prime > 0$, depending only on $\delta$ and $K$, such that whenever $S$ is a $\delta$-cluster of singularity set of $q \in K$ with respect to the $\cc$-metric, $S$ is a $\delta^\prime$-cluster of singularities with respect to the $q$-metric.\\
\item Conversely, for each sufficiently small $\delta > 0$ we can find $\delta^\prime > 0$, depending only on $\delta$ and $K$, such that such that whenever $S$ is a $\delta$-cluster of singularity set of $q \in K$ with respect to the $q$-metric, $S$ is a $\delta^\prime$-cluster of singularities with respect to the $\cc$-metric.\\
\item The two statements above are true for strata of cluster differentials in place of differentials in $QD(\mathcal{T}_{0,n})$.
\end{itemize}
\end{proposition}

\noindent In other words, we can tell when a collection of singularities is collapsing strictly faster than any proper superset using either metric.\\

\noindent Proof: By scale invariance of the ratios of distances in both metrics, we will assume all quadratic differentials in $K$ are of the form $\frac{P(z)}{Q(z)}dz^2$ with $P,Q$ monic of fixed degrees. $P$ is allowed to have multiple roots and $Q$ is not, and $P$ and $Q$ may have roots in common. By applying M\"obius translations fixing $\infty$, we can also assume that the set $\Sigma$ of singularities of each $q$ and is contained in a ball of radius $R$ about $0$ but not in a ball of radius $r$, for some $R > r > 0$. A closed subset of the space of differentials with these restrictions is compact if and only if there is a lower bound on the distance between any two roots of $Q$; conversely, for any $K$ we can find some $R$ and $r$ constraining the singularities of $q$ for all $q \in K$, after a M\"obius transformation.\\

\noindent Let $\{q_n\}_{n = 1}^\infty$ be a sequence in $K$, uniformized to $\hat{\cc}$ so as to satisfy the constraints in the above paragraph, and let $S_n$ be a collection of singularities in the $q_n$-metrics. Then by Corollary \ref{CoroBallsizes} we see that $\diam_\cc(S_n) \to 0$ if and only if $\diam_{q_n}(S_n) \to 0$. We may rephrase this as $\diam_{q_n}(S_n)/\Delta(q_n) \to 0$ if for each $\delta > 0$, $S_n$ is eventually contained in a $\delta$-cluster with respect to $d_\cc$, if and only if for each $\delta > 0$, eventually $S_n$ is eventually contained in a $\delta$-cluster with respect to $d_{q_n}$.\\

\noindent Now let $S_n \subsetneq S_n^\prime$ be sets of singularities such that $S_n^\prime$ is a $\delta_n$ cluster in $q_n$, and $\diam_{q_n}(S_n) \to 0$. By Lemma \ref{LemOutsiders} we can ignore all singularities outside of $S_n^\prime$ and use the above argument deduce that $\frac{\diam_\cc(S_n)}{\diam_\cc(S_n^\prime)} \to 0$ if and only if $\frac{\diam_{q_n}(S_n)}{\diam_{q_n}(S_n^\prime)} \to 0$.\\

\noindent The proposition then follows by induction on the total number of singularities. $\Box$

\begin{definition}Let $\delta < 1/10.$ We say that a saddle connection is \emph{internal} to a $\delta$-cluster $D$ in of the $q$-metric if it stays in the $2\diam_q(D)$ neighborhood of $D$ (in the $q$-metric). (In particular, this includes any saddle connection that is the shortest path from one endpoint to the other.)\end{definition}

\begin{definition}Given a meromorphic quadratic differential $\frac{P(z)}{Q(z)}dz^2$ on $\cc$ whose poles are all simple, we say that a $\delta$-cluster of singularities (where $\delta$-cluster is with respect to the flat metric) is \emph{shrunk} if it contains more zeros than poles (counted with multiplicity). Given a shrunk $\delta$-cluster of singularities, we say that the \emph{center} of the $\delta$-cluster is the weighted average of the singularities, where each point is counted with weight equal to the order of vanishing of $\frac{P(z)}{Q(z)}$ at that point (simple poles count as $-1$.)\end{definition}

\noindent These definitions are motivated by the fact that the diameter in the singular $q$-metric shrinks much faster than the diameter in the nonsingular metric as a cluster converges to a point and all other singularities remain fixed, and all other periods in a good coordinate system are well-approximated by the periods of the quadratic differential that is obtained when we replace the $\delta$-cluster with a single singularity at its center, which is a zero of order equal to the sum of the orders of vanishing of $P/Q$ at the singularities being replaced.\\

\noindent The equivalence of complex coordinate systems then follows because a holomorphic local homeomorphism of complex manifolds is a local biholomorphism - that is, the locations of the singularities in $\cc$ are holomorphic functions of the periods.

\subsection{Modelling Collisions of Singularities by Cluster Differentials}

\noindent We claim that cluster differentials form a model for a neighborhood of any disk without poles:

\begin{proposition}If $(X,q)$ is a simply connected Riemann surface equipped with a holomorphic quadratic differential with finitely many zeros, and the metric $d_q$ is complete, then $(X,q)$ is determined up to isometry by the isometry type of its geodesic convex hull of its zeros and the exterior angles at the vertices of the convex hull.\end{proposition}

\noindent Proof: Assume that a ball of radius $r$ about a singularity $p$ contains all singularities. By Alexandrov non-positive curvature, there is a unique geodesic between any pair of points. Moreover, geodesic rays from a fixed point $p$ have a circular order from their initial direction, and the directions they exit singularities if they enter the same singularity. We can therefore construct all geodesic rays from $p$, and reconstruct the complement of the geodesic convex hull of the zeros as a union of sectors of disks of infinite radii. $\Box$

\begin{proposition}\label{PropEverythingIsCluster}If $(X,q)$ is a (noncompact) simply connected (noncompact) Riemann surface with a complete metric $d_q$ coming from a holomorphic quadratic differential $q$, and $q$ has finitely many zeros and no poles, then $(X,q)$ is a cluster differential. \end{proposition}

\noindent Proof: Given $(X,q)$, we will construct a sequence of quadratic meromorphic quadratic differentials $\{q_n\}$ on $\hat{\cc}$ that converge to a holomorphic quadratic differential $q_\infty$ on every compact subset of $\cc$. Consider the balls $B_R(p)$ in $(X,q)$, for some fixed singularity $p \in X$. Now, the circle $c_R(p)$ of radius $R$ about $p$ has geodesic curvature $1/(R-s)$ at $z$ if the geodesic from $p$ to $z$ last changed direction at a cone point distance $s$ from $z$. In particular, the geodesic curvature of $c_R(p)$ is between $\frac{1}{R}$ and $\frac{1}{R-r}$ if $B_r(p)$ contains all singularities of $d_q$, except at finitely many points where the curvature is undefined, and never changes sign where the curvature is discontinuous.\\

\noindent STEP 1: CONSTRUCTION OF $\{q_n\}$. There are geodesic tangent lines to $c_R(p)$, and a finite number of these lines are vertical or horizontal. For large $R$, the number of such lines remains fixed, their points of tangency vary continuously in $R$, and they alternate vertical and horizontal. If $V(R)$ and $H(R)$ are consecutive vertical and horizontal tangent lines with respect to the circular order, then $V$ intersects $H$. The convex hull of all such intersection points is a piecewise geodesic arc with right-angles at its singularities, and it bounds a disk $D_R$. We will take $(X_n,q_n)$ to be a sequence of quadratic differentials to be the union of disks $D_{R_n}$ with their Schwarz reflections about the boundary (as metric spaces). Note that this Schwarz reflection gives us a canonical way to extend the vertical and horizontal foliations, and creates cone points of angle $\pi$ at the corners of $\partial D_R$. However, we have not yet described the sequence of uniformization maps on $X_n$ that make $q_n$ convergent on compact sets.\\

\noindent STEP 2: PICKING COORDINATES FOR $\{X_n\}$: We note that if we rescale $q_n$ to have unit area, then the $d_q$ distances between poles remain bounded below. (The area of $D_R$ is a quadratic polynomial in $R$ for large $R$, and the distances between poles grow linearly.) Thus we can apply Mumford's compactness criterion. By the equivalence of $\delta$-clusters in $d_{q_n}$ and $d_\cc$, for any sufficiently small positive $\delta$ there is an $N < \infty$ such that set of zeros on $D_{R_n}$ forms a $\delta$-cluster in $\cc$ for all $n > N$ in both $d_\cc$ and $d_q$. However, it is also true that for any sufficiently small $\delta > 0$ no proper subset of the zeros in $D_{R_n}$ forms a $\delta$-cluster, for any sufficiently small $\delta$ and large $n$. Fix any pair of vertical and horizontal tangent lines $V$ and $H$ to have intersection equal to $\infty$ under choice of uniformizations $X_n \to \hat{\cc}$. Then, there is a sequence of M\"obius transformations $T_n(z) = a_nz + b_n$ such that the collection of zeros converges in the configuration space $\mathrm{Conf_n(\cc)},$ along a subsequence, since ratios of $\cc$-distances between zeros in $D_{R_n}$ remain bounded above and below. Moreover, the locations of all other singularities go to $\infty$.\\

\noindent STEP 3: PROVING UNIFORM CONVERGENCE ON COMPACT SETS: it follows that on $\cc$, each $q_n$ takes the form $p_n(z)h_n(z)dz^2$, where we may assume (passing to a subsequence if necessary) that $p_n(z)$ has convergent coefficients, and $h_n$ is a rational function whose zeros and poles tend to $\infty$ uniformly in $n$. Therefore, there are numbers $C_n$ such that $h_n(z) = C_n(1 + o(1))$ on any ball as $n \to \infty$. Moreover, since the periods of saddle connections converge, we see that $h_n$ must converge uniformly to a \emph{fixed constant} on each compact subset of $\cc$. Moreover, judicious choices of $a_n$ and $b_n$ allow us to assume that $p_n$ is monic and its roots sum to $0$. $\Box$

\begin{corollary}Let $B$ be the closed ball of $d_q$-radius $1$ in a quadratic differential $(X,q)$ about a point $p \in X$. If $B$ has compact closure in $X$ and is homeomorphic to a closed disk and $q$ has no poles except possibly at $p$, then the interior of $B$ is isometric to a ball in a cluster differential. \end{corollary}

\noindent Proof: The case of no poles is obvious, and the case of one pole follows from taking a double cover branched only over the pole. $\Box$\\

\begin{definition}\label{DefCompactification} Let $\mu$ be a set whose elements correspond to singularities of a quadratic differential in $\mathcal{T}_{g,n}$ or a space of cluster differentials. To each element of $\mu$ we associate an integer that no less than $-1$ corresponding to order of vanishing, and a boolean variable which is $1$ for marked points and $0$ for unmarked points. Let $\mathcal{Q}(\mu)$ be the corresponding stratum in the Teichm\"uller space of half-translation surfaces.\\

\noindent For each stratum $\mathcal{Q}(\mu)$ of quadratic differential or cluster differential with at least one zero and at least two singularities, there is a partial compactification $\bar{\mathcal{Q}}(\mu)$, each of whose elements consist of the following data, up to an equivalence relation:

\begin{itemize}
\item A collection of subsets $S_\ell$ of $\mu$ of size at least $2$, such that each $S_\ell$ contains at most one singularity corresponding to a marked point. We require that the collection of subsets be totally ordered with respect to $\subsetneq$. We also require $S_0 = \mu$ to be one of the subsets. Let $m_\ell$ be the sum of the orders of vanishing of singularities in $S_\ell$ (poles count as $-1$).
\item A directed tree $T$ whose vertices consist of the sets $S_\ell$. Let the edges be such that there is a directed path from $v_{\ell_1}$ to $v_{\ell_2}$ if and only if $S_{\ell_1} \subset S_{\ell_2}$. For each vertex $S_\ell \neq S_0$, let $\phi(\ell)$ be the unique subset of $\mu$ such that there is a directed edge from $S_\ell$ to $S_{\phi(\ell)}$.
\item A quadratic differential $(X_\ell,q_\ell)$ with one singularity for each element $p$ of $\mu$ in $S_\ell$ not in any $S_j \subsetneq S_\ell$, and one singularity for each $S_j \subsetneq S_\ell$. Each $p$ not in any $S_\ell$ corresponds to the same type of singularity as $p$ in $\mu$, and for each $S_j$ there is a singularity whose order of vanishing is the sum $m_j$ of orders of vanishing of the singularities in $S_j$, and which is a marked point if and only if $S_j$ contains a marked point.
\item We require that $(X_0,q_0)$ have genus $g$ with $n$ marked points. For $\ell \neq 0$ we require that $(X_\ell,q_\ell)$ be a cluster differential with a pole of order $m_\ell+4$ at $\infty$, with a marked point if and only if $S_\ell$ has a marked point. We also require that the $d_\cc$ diameter of the set of singularities of $S_\ell$ is $1$.
\item For each oriented saddle connection $\gamma_i$ of the quadratic differential corresponding to vertex $\phi(S_\ell)$ that starts at the singularitiy corresponding to $S_\ell$, we associate an angle $\theta(\gamma_i) \in \rr/2\pi\zz$ such that the direction of $\gamma_i$ is counterclockwise from the direction of any other such saddle connection $\gamma_k$ by a (cone) angle of $[\theta(\gamma_i) - \theta(\gamma_k)](m_\ell+2)/2$.
\end{itemize}
The equivalence relation is generated by the condition that if $q_\ell = f(z)dz^2$, then $f(e^{i \alpha} z)[d(e^{i \alpha}z)]^2$ and the collection of angles $\theta(\gamma_k) - \alpha$.
\end{definition}

\noindent This is similar to the compactification of strata in \cite{5guys}, except that we do not consider quadratic differentials that escape to $\infty$ in the moduli space of quadratic differentials, and we consider the limiting objects to be different if its components are rotated relative to each other, whereas we are less specific with regard to the relative sizes of different components.\\

\noindent We will not need to construct the topology of $\bar{\mathcal{Q}}(\mu)$, though it can be constructed real analytically by repeatedly taking blowups (over $\rr$) of the closure of $\mathcal{Q}(\mu)$ along the loci where singularities collide and taking finite branched covers of the blown up subspace. However, we describe when a sequence in $\mathcal{Q}(\mu)$ converges to an element of $\bar{\mathcal{Q}}(\mu)$.

\begin{definition}Let $\{(Y_N,q_N)\}_{N = 1}^\infty$ be a sequence of quadratic differentials in $\mathcal{Q}(\mu)$. We say the sequence converges to a given element of $\bar{\mathcal{Q}}(\mu)$ if for every $\delta \in (0,1)$, there exists some $M < \infty$ such that
\begin{itemize}
    \item $d_\mathrm{Euclidean}((Y_N,q_N),(X_0,q_0)) < \delta$ whenever $N > M$.
    \item $\delta$-clusters in the set of singularities of $(Y_m,q_m)$ correspond exactly to the sets $S_\ell$ whenever $M > N$. That is, they have the same number of singularities of each type, same nesting, graph of inclusions, etc. Call these clusters $S_\ell(N)$. For each $\ell$ we will refer to the sequence $S_\ell(N)$ as a \emph{vanishing cluster}.
    \item For each $S_\ell$, $\ell \neq 0$, there is a disk $\mathbb{D}_{\ell}(N) \subset Y_N$ containing only the singularities in $S_\ell(N)$. For some $\lambda_N \in \cc^\times$, $(\mathbb{D}_{\ell}(N), d_{\lambda_N q_N})$ is locally isometric to the disk $\{z: z < 1/\delta\}$ in the $q_\ell^\prime$-metric on $\cc$, for some cluster differential $(\cc,q_\ell^N)$, where $d_\mathrm{Euclidean}(X_\ell^\prime,q_\ell^N) < \delta$.
    \item For each oriented saddle connection $\gamma_k$, if we pick a sequence of saddle connections $\gamma_k^N$ coverging to $\gamma_k$ in the universal half-translation surface of containing $q_{\phi(\ell)}$, and the choice of disks and isometries above, the $d_{q_i^N}$ ray in $\cc$ that has the same initial location and direction as the image of $\gamma_k$, when parametrized by arc-length, is of the form $\rr(t) e^{i\alpha(t)}, t \in [0,\infty)$ and $|\theta(\gamma_k) - \lim\limits_{t \to \infty} \alpha(t)| < \epsilon (mod ~ 2\pi/zz).$
\end{itemize}
\end{definition}

\begin{proposition} Every convergent sequence in $QD(\mathcal{T}_{g,n}$ consisting of elements of $\mathcal{Q}(\mu)$ has a subsequence that converges to an element of $\bar{\mathcal{Q}}(\mu)$.\end{proposition}

\noindent Proof: This is just a simple matter of extracting subsequences, since the space of cluster differentials with $\diam_\cc = 1$ is compact for each stratum of cluster differential, as well as the fact that the space of directions $\theta$ is compact. $\Box$\\

\noindent We remark that there is an equivalent characterization of the convergence to an element of $\bar{\mathcal{Q}}(\mu)$ when the underlying Riemann surface is $\hat{\cc}$ and the differentials are M\"obius normalized. In this case, there is a sequence of M\"obius transformations $T_i(N),$ all fixing $\infty$, such that the the images of singularities in $S_\ell(N)$ under $T_i(N)$ converge to those of a suitable cluster differential. By Proposition \ref{PropPolyScale}, Lemma \ref{LemOutsiders}, and Proposition \ref{PropScale} it follows that this sequence of cluster differentials can be rescaled to converge in the Euclidean metric, so that their limit is isometric to $(\cc,d_{q_i})$. It should also be noted that if $q_i$ has no poles, then for some numbers $\lambda_N \in \cc,$ $\{\lambda_N (T_N^*)^{-1}(q_N) dz^{-2}\}_{N = 1}^\infty$ is a uniformly convergent sequence of holomorphic functions on each compact subset of $\cc$. For each vanishing cluster, let $z_{j,N}$ be a sequence of elements corresponding to a singularity in $S_\ell(N)$. After possibly permuting singularities of the same type, we conclude that there is a sequence of numbers $a_N \in \rr$ and $b_n \in \cc$ such that $T_N(z) = a_N(z) + b_N$, with $a_N \to \infty$, such that the singularities $T_N(z_{j,N})$ all converge to the singularities of a the differential $q_i$. In summary, we have the following:

\begin{proposition}\label{PropMobiusClusters}Let $\{q_N\}_{N = 1}^\infty$ be a M\"obius normalized sequence of quadratic differentials on $\hat{\cc}$ in a stratum of cluster differentials or in $QD(\mathcal{T}_{0,n})$, converging to an element of $\bar{\mathcal{Q}}(\mu)$ with a pole at $\infty$. If $S_\ell$ is a vanishing cluster, then after permuting singularities $z_{j,N}$ of the same type in $S_\ell(N)$, there exist M\"obius transformations $T_N(z) = a_N z + b_N$, $a_N \in \rr^+$, $b_N \in \cc$ with $a_N \to \infty$, such that the singularities $T_N(z_{j,N})$ converge to the singularities of the differential $q_i$. There are scalars $\lambda_N \in \cc$ such that $\lambda_N ((T_N)^*)^{-1}$ converges uniformly on compact subsets of $\cc$ if $S_\ell$ does not contain a pole, and on compact subsets of $\cc \setminus \{0\}$ if $S_\ell$ contains a pole. $\Box$\end{proposition}

\noindent Because of the convergence of the rescaled quadratic differentials, we may also conclude the following version of Corollaries \ref{CoroHausCont} and \ref{CoroC1}:

\begin{corollary}\label{CoroC1forclusters}Let $\lambda_N,T_n,q_N,S_\ell(N)$ be as in Proposition \ref{PropMobiusClusters}. If $\gamma_N: [0,1] \to \cc$ is a sequence of constant speed saddle connections internal to $S_\ell(N)$ then $T_N(\gamma_N)$ has a Hausdorff convergent subsequence. If the limit is a saddle connection in $q_\ell$ then the convergence is $C^1$ on compact subsets of $(0,1). \Box$\end{corollary}

\noindent We conclude this section by collecting several simple facts about cluster differentials whose proofs we defer to Appendix B. These results are needed for the proof of Lemma \ref{LemmaPerturbEff}.

\begin{definition}Let $\delta \in (0,1)$. A \emph{possibly collapsed} $\delta$-cluster is a $\delta$-cluster or a single point.\end{definition}

\begin{proposition}\label{PropPersistentCluster} Let $K$ be a compact subset of $QD^1(\mathcal{T}_{g,n})$ and let $(X,q) \in K$ have singularity set $\Sigma$. Let $\gamma: [0,1] \to QD(\mathcal{T}_{g,n})$ be a rectifiable path with respect to the Euclidean metric on $QD(\mathcal{T}_{g,n})$ starting at $\gamma(0) = (X,q) \in K$. We write $\gamma(t) = (X(t),q(t))$, and we let $\Sigma(t)$ be the singularity set of $\gamma(t)$.\\

Let $\Sigma^\prime \subset \Sigma$ let $\Sigma^\prime(t)$ be a finite subset of $X(t)$ that varies continuously in the Hausdorff topology along $\gamma(t)$, and assume that $\Sigma^\prime(0)$ is a possibly collapsed $\delta$-cluster for some $\delta > 1$.\\

Then for any $\delta_1 \in (0,1)$ there are positive real numbers $C,\delta_2 > 0$, depending only on $K$ and $\delta_1$, such that the following holds: if $\Sigma^\prime(0)$ is a $\delta_2$-cluster of Euclidean length at most $C\diam_q(\Sigma)/\delta_2$ in $QD(\mathcal{T}_{g,n})$, there is a unique way to choose $\Sigma^\prime(t)$ for all $t$ such that the number of marked points in $\Sigma(t)$ and the number of zeros minus poles (counted with multiplicity) remain constant along $\gamma$, even if the cardinality of $\Sigma^\prime(t)$ does not remain constant. \end{proposition}

\begin{proposition}\label{PropClusterProj} Let $K \subset QD(\mathcal{T}_{g,n})$ be compact. Let $U$ be a contractible open subset of $QD(\mathcal{T}_{g,n})$ with closure in $K$. Suppose that there is a continuous choice of possibly collapsed $\delta$-cluster in $U$, that is to say, a continuous function $h$ from $U$ to the space of compact sets in the universal half-translation surface of type $g,n$ (endowed with the topology of the Hausdorff metric) that maps each element $(X,q) \in U$ to a finite collection of singularities of $(X,q)$, such that $h(X,q)$ is a possibly collapsed $\delta$-cluster in the singularity set of $(X,q)$ and moreover, the number of marked points and sum of the orders of vanishing of quadratic differentials at points in the image of $h(X,q)$ is constant. (By proposition \ref{PropPersistentCluster}, if $\delta$ is small enough and the diameter of $U$ is small enough, then any such function $h$ is uniquely determined by its value at any point in $U$.)

Let $M$ be the moduli space of cluster differentials with the same number of marked points and total order of vanishing as $h(X,q)$. Then there is a number $\delta_K$ such that if $\delta < \delta_K$, and $U$ and $M$ are given Euclidean metrics in the sense of Definition \ref{DefineEuclidean}, there is a Lipschitz map $F: U \to M$ with the following property: a disk containing all of the singularities of $F(X,q)$ maps isometrically into $(X,q)$ such that singularities of $F(X,q)$ map to $h(X,q)$ and the vertical foliation is preserved. The Lipschitz constant depends only on $K$.\end{proposition}

\noindent Proof: This is immediate from local finiteness of period coordinate systems coming from bounded length saddle connections, i.e. from Proposition \ref{PropLocallyFinite}. $\Box$\\

\begin{lemma}\label{LemAsympClust}Let $\{q_N\}$ be a sequence in $QD(\mathcal{T}_{0,n})$ that converges to an element of $\bar{\mathcal{Q}}$, and Fix a vanishing cluster $\{S_\ell(N)\}$ for the sequence $q_N$. Assume $q_N$ is given the normalization such that one pole, which is not part of a vanishing cluster, is at $\infty$, the center of $S_\ell(N)$ is always $0$, and $q(z)dz^{-2}$ is a quotient of monic polynomials. Assume the sum of the orders of vanishing of $q_N$ on $S_\ell(N)$ is $m$.

\noindent Let $z_{j,N}$ range over the non-infinite singularities of $q_N$, with $q_N$ vanishing at $z_{j,N}$ to order $e_{j,N}$, and let $$t_N = \prod\limits_{z_{j,N} \notin S_\ell(N)} (-z_{j,N})^{e_{j,N}},~\mathrm{and}$$ $$\alpha_N = \prod\limits_{z_{j,N} \in S_\ell(N)} \left(z - t_N^{1/(m+2)}z_{j,N}\right)^{e_j} dz^2.$$

Let the map $F_N$ be the locally defined map $F$ for the cluster $S_\ell(N)$ to the space of cluster differentials associated to $S_\ell(N)$ as defined in proposition Proposition \ref{PropClusterProj}. Then, for an appropriate choice of the $(m+2)^{\mathrm{nd}}$ root of $t_N$, $d_\mathrm{Euclidean}(\alpha_N, F_N(q_N)) = o(\diam_{q_N}(S_\ell(N))).$ \end{lemma}

\noindent Proof: This is deferred to Appendix B. $\Box$\\

\noindent The following proposition says that if we perturb a quadratic differential but preserve the isometry type all $\delta$-clusters for some sufficiently small $\delta$, then for a reasonable period coordinate chart, the change in the periods is comparable to the Euclidean distance.\\

\begin{proposition}\label{PropNoCancellation}Let $K$ be a compact subset of $QD(\mathcal{T}_{g,n})$ or a space of cluster differentials, and let $L > 0$, and assume that on $K$ the systems of saddle connections associated to the Euclidean metric on $K$ consists of saddle connections of length less than $L$. Then there is some $\delta_0 > 0$ such that for all $\delta < \delta_0$, the following holds: as $X$ varies in $K$ with singularity set $\Sigma$ and orienting double cover $(\tilde{X},\tilde{\Sigma})$, let $H_\delta^1(\tilde{X},\tilde{\Sigma}; \cc)$ denote the subspace of $H^1(\tilde{X},\Sigma,\cc)$ that vanishes on all saddle connections internal to $\delta$-clusters of singularities.

Then, if $\{X_m\}_{m = 1}^\infty$ is a convergent sequence in $K$ with singularities $\{\Sigma_m\}$ and $Y_m$ is another such sequence lying on a common convex $T$-convex period coordinate chart $U \subset \cc^n$ with $X_m$ for some $T$, and $U$ is generated by periods of saddle connections of length less than $L$, and $X_m$ is connected to $Y_m$ by a line segument of length $\epsilon_m \to 0$ in $U$, and such that, in local period coordinates, $X_m - Y_m \in H_\delta^1(\tilde{X}_m,\tilde{\Sigma}_m; \cc) \cap H_{odd}^1(\tilde{X}_m,{\Sigma}_m,\cc)$, then $d_\mathrm{Euclidean}(X_m,Y_m) \dot{\asymp} \epsilon_m$. \end{proposition}

\noindent Proof: Deferred to Appendix B. $\Box$

\section{Perturbing Quadratic Differentials on the Sphere}

\noindent We now establish a setup to discuss saddle connections in convergent sequences of quadratic differentials on $\hat{\cc}$. We will consider all strata of quadratic differentials with a pole fixed at $\infty$, and such that all poles in $\cc$ are simple. We simply want to know that our quadratic differentials have enough saddle connections for us to apply Corollary \ref{CoroAnybasis}. The construction in Definition \ref{DefCompactification} also applies to such differentials.\\

\subsection{Derivatives of Period Coordinates}

\begin{notation}\label{NotaModelnotation}
Fix integers $e_1,...e_r \geq -1$. Assume $\sum\limits_{m = 0}^r e_m \geq -3$. For any $r$-tuple of distinct complex numbers $(z_1,...z_r)$, let $q(z_1,...,z_r) = \prod\limits_{j = 1}^r (z - z_j)^{e_j}dz^2$.

Given a finite system of saddle connections $\gamma_i$ for $q(w_1,...w_r)$ with $\gamma_i$ having endpoints $w_{i_1},w_{i_2}$, there is a neighborhood $U$ of $(w_1,...,w_r)$, such that for all $(z_1,...,z_r) \in U$ the saddle connection in $(\cc,q(z_1,...,z_r)$ with endpoints $z_{i_1}$ and $z_{i_2}$ can be chosen Hausdorff continuously for all $\ell \in I$, such that $\gamma_i$ is the collection of saddle connections at $(w_1,...,w_r)$. The periods vary holomorphically with respect to $(z_1,...,z_r)$; call these periods $P_i$, and their partial derivatives $\frac{\partial P_i}{\partial z_j}$.
\end{notation}

\noindent Differentiation under the integral gives us the following basic formula:

\begin{proposition}\label{PropDifunderint}1. If $\gamma_i$ is a saddle connection of the quadratic differential $q(z) = f(z)dz^2$ with period $P_i$, and either $e_j \geq 1$ or $v_j$ is not an endpoint of $\gamma_i$, then we have $$\frac{\partial P_i}{\partial z_j} = \int_{\gamma_i} \frac{e_j\sqrt{q(z)}}{2(z_j - z)}.$$

\noindent 2. In all other cases, we can recover the partial derivative $\frac{\partial P_i}{\partial z_j}$ by performing a 1-parameter family of M\"obius transformations so that the endpoints of $\gamma_i$ remain fixed and pushing forward the family of differentials, and differentiating under the integral.\end{proposition}

\noindent Proof: The only cases in which the formula for $\frac{\partial P_i}{\partial z_j}$ is not the immediate result of differentiation under the integral are those in which on endpoint of $\gamma_i$ is $z_j$ and $e_j \geq 1$. Suppose $a$ and $b$ are the endpoints of $\gamma_i$ and $b = z_j$. Assume $q_h(z) = q(z)\left(\frac{z - b - h}{z - b}\right)^{e_j}.$ Fix a contour of integration from $a$ to $b$ to be the saddle connection $\gamma_i$ for the $q$-metric. Then, $$\frac{1}{h}\left[\int_a^{b+h} \sqrt{q_h(z)} - \int_a^b \sqrt{q(z)}\right] = \int_a^b \frac{1}{h}[\sqrt{q_h(z)} - \sqrt{q}] + \int_b^{b+h}\frac{1}{h}q_h(z).$$

\noindent We claim that the last integral on the right hand side is $O(h^{e_j/2})$ because we can take our contour of integration from $b$ to $b+h$ to have length $h$ and the integrand is $O(h^{e_j})$ along the entire contour. The remaining term limits to the desired formula.\\

\noindent In other cases, the integrals need not converge. To prove our remaining claim, we simply note that any holomorphic 1-parameter family of quadratic differentials $\prod(z - z_j(t))^{e_i}dz^2$ admits a holomorphic 1-parameter family of M\"obius transformations sending $z_{j_1}(t),z_{j_2}(t),\infty$ to $z_{j_1}(0)$ and $z_{j_2}(0)$. The push-forwards of $q(t)$ will be of the form $(z-z_{j_1})$. $\Box$

\subsection{The Limit of the Matrix $\frac{\partial P_i}{\partial z_j}$}

\begin{notation}\label{NotaCompactConverge}
Fix a sequence of quadratic differentials $q_N = q(z_{1,N},...,z_{r,N})$ that converge to an element of the partial compactifcation of Definition \ref{DefCompactification}. Assume that the length-minimizing spanning trees of saddle connections are isomorphic to some fixed $\Gamma$ as ribbon graphs, each $z_{j,N}$ corresponds to a fixed vertex $v_j$ of $\Gamma$ for all $N$. Also assume that for each vanishing cluster $S_\ell$ such that $z_{j,N} \in S_\ell(N)$ there is some sequence $T_n$ as in Proposition \ref{PropMobiusClusters} such that $T_N(z_{j,N})$ converges for all $j$ such that $z_{j,N} \in S_\ell(N)$ to the corresponding singularity of the cluster differential $q_\ell$ corresponding to $S_\ell$.

For this sequence, let $\gamma_i$ vary over the edges of $\Gamma$ and have period $P_i$. Let $\Gamma_{i,N}$ be the saddle connection corresponding to $\gamma$ in $(\cc,q_N)$. Assume we have picked $\sqrt{q_N}$ so that for some sequence of positive real numbers $t_N$, the sequence $\int_{\gamma_{i,N}} t_N \sqrt{q_N}$ converges to a nonzero value, and let $P_i$ be the holomorphic function associated to $\gamma_i$.

Finally, let $M_N$ be matrix $(M_N)_{ij} = \left.\frac{\partial P_i}{\partial z_j}\right|_{(z_1,...,z_r) = (z_{1,N},...,z_{r,N})}.$ Finally, let $M_N^\prime$ be obtained from $M_N$ by the following operation: for each maximal vanishing cluster $S_\ell$, delete the rows corresponding to the saddle connections internal to vanishing clusters, and replace the columns corresponding to singularities in vanishing clusters by a single column that is the sum of the columns corresponding to singularities in $S_\ell(N).$\end{notation}

\begin{proposition}\label{PropClustermodels} Let $\{q_N\}_{N = 1}^\infty = \{q(z_{1,N},...,z_{r,N})\}_{N=1}^\infty$ be a sequence satisfying the hypotheses of Notation \ref{NotaCompactConverge}, and the root of the associated tree is $q_0$. Let $m$ be such that $q$ has a pole of order $m+4$ at $\infty$. Assume no vanishing cluster contains a pole, and there is at most one singularity of cone angle $2\pi$ in each vanishing cluster. Then we have the following:\begin{enumerate}
\item For each $i$, $P_i(z_{1,N},...,z_{r,N})$ converges as $N \to \infty$.
\item $M_N$ converges to a matrix $M_\infty$.
\item If $(z_j)_\infty = (z_k)_\infty$, then $e_j(M_\infty)_{ik} = e_k (M_\infty){ij}$.
\item $M_N^\prime$ converges to a matrix $M_\infty^\prime$
\item The kernel of $M_\infty^\prime$ has dimension 1 and it is spanned by the all $1$ vector.
\item $M_\infty^\prime$ is the matrix of partial derivatives for the periods of the saddle connections and vertices of the degeneration of $T$ that occurs as $q_n \to q_\infty$.\\
\item If $\gamma_i$ corresponds to a sequence of saddle connections internal to a vanishing cluster, then $\frac{\partial P_i}{\partial z_j} \to 0$ for all $j$.
\end{enumerate}

More concisely, if $\cdot^\prime$ denotes the object that $\cdot$ degenerates to as $N \to \infty$, $$\frac{\partial P_i}{\partial z_j} \to \frac{e_j}{e_j^\prime}\frac{\partial P_i^\prime}{\partial z_j^\prime}.$$ If $e_j = e_j^\prime = 0$ then $\frac{e_j}{e_j^\prime} = 1$. (In particular, the period of a saddle connection that degenerates to a point has partial derivatives tending to $0$.) \end{proposition}

\noindent Proof: The first claim follows easily from Corollaries \ref{CoroHausCont}, \ref{CoroC1}, and \ref{CoroBallsizes}.\\

\noindent Items 2-7 determine the value of $M$, up to the value of the derivative of the saddle connection of a non-shrunk vanishing cluster.\\

\noindent Suppose there are no vanishing clusters. Then the period coordinates $P_i$ are a holomorphic coordinate system for our stratum of cluster differential in a neighborhood of $q_0$. Indeed, the periods of the $r-1$ edges of $\Gamma$ form a local coordinate system for the stratum of the differentials $q(z_{1,N},...,z_{r,N})$. It is also true that the values of $r-1$ elements of $\{z_1,...,z_r\}$ determine a local coordinate system if the remaining $z_j$ remains fixed, since $q(z_1,...,z_r)$ is biholomorphically equivalent to $q(z_1 + \alpha,...,z_r + \alpha)$ for any $\alpha \in \cc$ (each is the pullback of the other via translation). Since any injective holomorphic map from a domain in $\cc^n$ to a domain in $\cc^n$ is a biholomorphism onto its image, it follows that the kernel of $M$ consists only of perturbations of $(z_0,...,z_r)$ that preserve the isomorphism type of the quadratic differential $q_0$. This subspace is precisely the span of $\left[\begin{array}{c} 1 \\ ... \\ 1 \end{array}\right]$. In this case, the dimension of the stratum of the limiting dimension equals the claimed rank of $M$. It follows that all claims hold in this case.\\

\noindent Now, recall that the sequences of saddle connections $\gamma_{i,N}: [0,1] \to \cc$ each have the property that for some sequence of M\"obius transformations $T_{i,N}(z) = t_{i,N}(z) + u_{i,N}$ fixing $\infty$ as in Proposition \ref{PropMobiusClusters}, $T_{i,N} \circ \gamma_{i,N}$ converges uniformly and $c_1$ on compact subsets of $(0,1)$ to a saddle connection in some cluster differential $q_\ell$.\\

\noindent Assume that $T_{i,N}(z) = t_{i,N} z + u_{i,N}$ is chosen so that $T_{i,N}(z)\gamma_i(q_N)$ converges to the saddle connection corresponding to $\gamma_i$ on one of the quadratic differentials $q_\ell$ associated to the limit of the sequence $\{q_N\}$.\\

\noindent Now, we would like to apply Proposition \ref{PropDifunderint} with the change of coordinates $T_{i,N}(z) = \zeta$, so let $$\alpha_{i,N}(\zeta) = t_{i,N}^{(m+2)}\prod\limits_{j = 1}^r(\zeta - \zeta_j)^{e_j}d\zeta^2 = (T_{i,N})_*(q_N).$$ If $\gamma_{i,N}$ has endpoints $a_{i,N}$ and $b_{i,N}$ for the differential $q_N$, we then have $$\frac{\partial Q_i}{\partial z_j}|_{(z_1,...,z_r)=(z_{1,N},...,z_{r,N})} = \int_{a_{i,N}}^{b_{i,N}} \frac{e_j}{2} \frac{\sqrt{q_N(z)}}{z_{j,N} - z} = \int_{T_{i,N}(a_{i,N})}^{T_{i,N}(b_{i,N})} t_{i,N}^{-m/2}\frac{e_j}{2}\frac{\sqrt{\alpha_{i,N}(\zeta)}}{\zeta_{j,N} - \zeta}.$$

\noindent We would like to apply the Lebesgue dominated convergence theorem to the right-hand side, parametrizing the integrals to be constant speed $d_{(T_{i,N})_*(q_N)}$-geodesic maps from $[0,1]$ to $\cc$. To do this, we need the following:

\begin{claim}A dominating function for the sequence $\int_{T_{i,N}(a_{i,N})}^{T_{i,N}(b_{i,N})} t_{i,N}^{-m/2}\frac{e_j}{2}\frac{\sqrt{\alpha_{i,N}(\zeta)}}{\zeta_{j,N} - \zeta}$, where the contours are parametrized as constant speed geodesics from $[0,1]$ to $\cc$ with respect to $d_{T_*{q_N}}$, exists whenever $z_{j,N}$ is a zero of $q_N$ or not an endpoint of $\gamma_{i,N}$.\end{claim}

\noindent Proof: Suppose that $\{\gamma_{i,N}\}$ is not internal to a vanishing cluster.\\

\noindent Since we are assuming that our contours of integration belong to length-minimizing spanning trees, For each $z$ on our contour of integration, we have $$d_{q_N}(z,(z_{j,N}) \geq  \min(d_{q_n}(z,(a_{i,N})),d_{q_N}(z,b_{i,N})).$$ If this were not the case, we would have $$d_{q_N}(a_{i,N},b_{i,N}) > \min(d_{q_N}(a_{i,N},z_{j,N}),d_{q_N}(b_{i,N},z_{j,N}),$$ and a saddle connection from $a_{i,N}$ to $b_{i,N}$ could not belong to a length-minimizing spanning tree.\\

\noindent All singularities of $q_N$ remain bounded, the poles are all bounded away from $z_{j,N}$, and $z_{j,N}$ is a zero of $q_N$, so we therefore have $$d_{q_N}(z, z_{j,N}) \dot{\prec} \int_0^{|z - z_{j,N}|} t^{1/2} dt = 2|z - z_{j,N}|^{3/2}/3.$$

\noindent We thus have, along our contour of integration, $$|z - z_{j,N}|^{-1} \dot{\prec} d_{q_N}(z,z_{j,N})^{-2/3} \leq \min (d_{q_N}(z, T_{i,N}(a_{i,N})),d_{q_N}(z, T_{i,N}(b_{i,N})))^{-2/3}.$$

\noindent We have parametrized our arcs to be $d_{q_N}$ geodesics instead of $d_{\alpha_n}$ geodesics, but applying Lemma \ref{LemOutsiders} to the proofs of Corollary \ref{CoroHausCont} and Corollary \ref{CoroC1} tells us that the contours of integration converge (Hausdorff, in total length, and $C^1$ on compact subsets of $(0,1)$) to geodesics with respect to $d_{\alpha_\infty}$, where $\alpha_\infty = \lim\limits_{n \to \infty} \alpha_n$. Therefore, for some $C > 0$, our sequence of integrals is dominated by $\int_0^1 C[x(1-x)]^{-2/3}dx.$ This completes the proof of the claim in this case. $\Box$

\subsection{The Derivative of Period Coordinates, Rescaled}

\noindent Now, we may do the same for the case of a saddle connection internal to a vanishing cluster. By translating our original sequence $q_N$ so that one endpoint of $\gamma_{i,N}$ is $0$ and not equal to $z_j$ we can assume that $\zeta = T_{i,N}(z) = t_{i,N}z$ with $t_{i,N} \to \infty$. If $\{S_\ell(N)\}_{N = 1}^\infty$ is the minimal vanishing cluster to which $\gamma_{i,N}$ is internal write $q_N = f_{i,N}(z)g_{i,N}(z)dz^2$ with $g(z) = \prod\limits_{j \in S_\ell(N)}(z - z_{j,N}^{e_j}).$ Then the sequence of integrals becomes $$\begin{array}{rcl} \left. \frac{\partial P_i}{\partial z_j}\right|_{(z_1,...,z_r) = (z_{1,N},...,z_{r,N})} & = & \int_{a_{i,N}}^{b_{i,N}} \sqrt{f(z)g(z)} \frac{e_j dz}{2(z - z_{j,N})}\\
{} & = & \int_{\alpha_{i,N}}^{\beta_{i,N}} [f_N(\zeta/t_{i,N})]^{1/2} [g_N(\zeta/t_{i,N})]^{1/2} \frac{e_j d\zeta}{2(\zeta - \zeta_{j,N})}. \end{array}$$

\noindent Now, the contour of integration converges, and $f_N(\zeta/t_{i,N}) = f_N(z)$ converges to a constant, possibly zero. If $m_\ell$ is the degree of $g(z)$ then $[t_{i,N}^{-m}h_N(\zeta)]$ converges as a function of $\zeta$. By Lemma \ref{LemOutsiders} and the argument for the saddle connections not internal to vanishing clusters, for $x \in [0,1]$ we have $$\left|\frac{h_N(\zeta)}{\zeta - \zeta_{j,N}} \circ T_{i,N} \gamma_{i,N}(x)\right| \leq C[x(1-x)]^{-2/3}.$$ We therefore see that not only does a dominating function exist, but we can take the dominating function to be $t_{i,N}^{-m/2}C[x(1-x)]^{-2/3}$, $m > 0$, if we are willing to start the sequence at some large $N$. This proves item 7. The rest of the claims follow from Proposition \ref{PropDifunderint} and the dominated convergence theorem.

\begin{corollary}\label{CoroLogmodelsprep}Let $q_N = \prod\limits_{j = 1}^s (z-z_{j,N})^{e_j}dz^2$ be a sequence of cluster differentials with no poles satisfying the hypotheses of Proposition \ref{PropClustermodels}, and assume that at least one of the sequences $\{z_{j,N}\}$ does not converge to $0$. Fix integers $e_{s+1},...,e_r$ and let $(z_{s+1,N},...,z_{r,N})$ be a sequence of $(s-r)$-tuples of complex numbers converging to $(\infty,\infty,...,\infty)$ in $\hat{\cc}^{r-s}$. Assume $\{\gamma_{i,N}\}$ is a sequence of saddle connections of the differential $q_N$ that converges (Hausdorff), and its period $P_i$ is locally a holomorphic function of $(z_1,...,z_s)$. Let $Q_{i,N}$ be the periods of a sequence saddle connections in the metric associated to the quadratic differential $\prod_{j = 1}^r (z - z_{j,N})^{e_j}$ which converge (Hausdorff) to same limit as $\gamma_{i,N}$, which are locally given by the holomorphic function $Q_i(z_1,...,z_r)$.\\

\noindent If $\gamma_{i,N}$ does not converge to a point, for $1 \leq j \leq s$, the following limits exist and are equal:

$$\lim\limits_{N \to \infty} \frac{\partial Q_{i,N}}{\partial z_j} = \frac{\partial \log Q_{i,N}}{\partial z_j} = \lim\limits_{N \to \infty} \frac{\partial \log P_i}{\partial z_j}.$$

\noindent If $\gamma_{i,N}$ converges to a point, then let $\gamma_{i^\prime,N}$ be another convergent sequence of saddle connections that does not converge to a point, and let them have periods $P_{i^\prime,N}$ and $Q_{i^\prime,N}$.

$$0 = \lim\limits_{N \to \infty} \frac{\partial Q_{i,N}}{\partial z_j} = \frac{\partial Q_{i,N}}{\partial z_j}\frac{1}{Q_{i^\prime,N}} = \lim\limits_{N \to \infty} \frac{\partial P_i}{\partial z_j}\frac{1}{P_{i^\prime,N}}.$$

\noindent For $s + 1 \leq j \leq r$,

$$\lim\limits_{N \to \infty} \frac{\partial Q_{i,N}}{\partial z_j} = \frac{\partial \log Q_{i,N}}{\partial z_j} = 0.$$ \end{corollary}

\noindent Proof: An identical dominated convergence argument. We can treat the faraway singularities (those outside $S_\ell{N}$) as scalars in the limit, and by Lemma \ref{LemOutsiders} we can use the same system of saddle connections $\gamma_i$ we would if there were no finite singularities outside of $S_\ell(N)$. $\Box$

\begin{corollary}\label{CoroLogmodels} Given a sequence $\{q_N\}$ converging as in Proposition \ref{PropClustermodels}, let $M_{\ell,N}$ be the submatrix of $M_N$ corresponding to saddle connections internal to the vanishing cluster $\{S_\ell(N)\}$. Then we have the following as $N \to \infty$:
\begin{itemize}\item $\left. \frac{\partial P_{i}}{\partial z_j} \right|_{(z_1,...,z_r) = (z_{1,N},...,z_{r,N})} = o(|M_{\ell,N}|)$ if $j \notin S_\ell$.
\item If $\gamma_{i,N}$ corresponds to a row of $M_{\ell,N}$ then $\sum\limits_{j \in S_\ell} \left. \frac{\partial P_{i}}{\partial z_j} \right|_{(z_1,...,z_r) = (z_{1,N},...,z_{r,N})} = o(|M_{\ell,N}|).$
\item If $S_l \subsetneq S_\ell$ is a vanishing cluster and $\sum\limits_{j \in S_l} a_j e_j = 0$, then $$\sum_{j \in S_l} a_j \left. \frac{\partial P_{i}}{\partial z_j} \right|_{(z_1,...,z_r) = (z_{1,N},...,z_{r,N})} = o(|M_{\ell,N}|).$$\end{itemize}\end{corollary}

\noindent Proof: This is immediate from Corollary \ref{CoroLogmodelsprep} and Proposition \ref{PropDerivScale}. $\Box$

\begin{corollary}\label{CoroFirstOrderOnClusterMagnitude}Let $\{q_N\}, M_{\ell,N}$ be as above, and assume the center of $S_{\ell,N}$ is $0$ for all $N$. Let $t_N = \prod\limits_{z_{j,N} \notin S_\ell(N)} z_j (-z_{j,N})^{e_{j,N}}.$ we have $|M_{\ell,N}| \dot{\asymp} |t_N|^{1/2} \diam_\cc(S_\ell(N))^{m/2},$ where $m$ is the number of zeros in $S_\ell(N),$ counted with multiplicity. In fact, the sequence of matrices $t_N^{-1/2}\diam_\cc(S_\ell(N))^{m/2}M_{\ell,N}$ converges to a matrix whose cokernel corresponds to the space spanned by periods of saddle connections internal to proper vanishing subclusters of $\{S_\ell(N)\}$. \end{corollary}

\noindent Proof: We are free to consider the same system of saddle connections we would use in the absence of faraway singularities by Lemma \ref{LemOutsiders}.\\

\noindent Assume $z_j$ a zero of order $z_j$ degenerates to $z_{j^\prime}$ a zero of order $e_j^\prime$ times the corresponding column of the matrix of partial derivatives for the limiting differential $q_\ell$. If we change coordinates on $\cc$ by the positive real scalar $\diam_\cc(S_\ell(N))$ so that $S_\ell(N)$ has diameter $1$, then Propositions \ref{PropPolyScale}, \ref{PropDerivScale} and \ref{PropClustermodels} imply that column $j$ of the matrix $|t_N|^{-1/2}\diam_\cc(S_\ell(N))^{-m/2}M_{\ell,N}$ converges to $\frac{e_j}{e_{j^\prime}}$, together with additional zeros for the saddle connections internal to proper subclusters of $q_\ell$. Since the matrix of partial derivatives for $q_\ell$ has a nonzero limit, so does $|t_N|^{-1/2}\diam_\cc(S_\ell(N))^{-m/2}M_{\ell,N}$, and the result follows. $\Box$\\

\noindent We would have liked to extend Proposition \ref{PropClustermodels} and its corollaries to clusters that include a single pole. Unfortunately, the partial derivatives of a period need not remain bounded when a zero collides with a pole. Consider a the family of quadratic differentials $\frac{z(z-1)dz^2}{z + t}$ where $t$ varies. If $t$ is a negative real number, then the segment $[0,1]$ of the real line is a vertical saddle connection; call its period $P(t)$. Then, by Proposition \ref{PropDifunderint} we have $$P^\prime(t) = \frac{1}{2}\int_0^1 \sqrt{\frac{x(x-1)}{(x+t)^3}}dx.$$

\noindent The integrand is purely imaginary and does not change sign. Thus, by the monotone convergence theorem, $$\lim\limits_{t \to 0^-}|P^\prime(t)| = \int_0^1 \frac{(1-x)^{1/2}}{x}dx \geq \int_0^{1/2}\frac{dx}{2x} = +\infty.$$

\noindent However, for the estimates we need, it is possible to take a double cover branched over the pole and pull back the quadratic differential; this converts the pole into a singularity with cone angle $2\pi$ at the cost of duplicating other singularities.\\

\subsection{The Main Estimate}

\begin{notation}\label{NotaFixpoles} Fix a stratum of quadratic differentials in $QD(\mathcal{M}_{0,n})$, with no zeros of order 2 or higher. Assume all quadratic differentials are normalized to have simple poles at $\infty$ and two other fixed points $A$ and $B$. Then each has the form $\lambda_q[(z - A)(z - B)]^{-1}\prod_{j = 1}^s (z - z_j)^{e_j}dz^2$. Fix integers $\{e_j\}_{j = 1}^s$ with $-1 \leq e_j \leq 1$ for all $j$ and $e_1 + ... + e_s = -1$ and denote such a differential as $q(\lambda,z_1,...,z_r)$. Let $$d_\mathrm{Sym}(q_1,q_2) = \inf_\phi\sum_w|\phi(w) - w|+|\lambda_{q_1} - \lambda_{q_2}|$$ where $w$ ranges over all non-infinite singularities of $q_1$, and $\phi$ ranges over all bijections between the singularities of $q_1$ and $q_2$ that preserve the type of singularity.\end{notation}

\begin{proposition}\label{PropStableClust}Suppose that $\{q_N\}$ and $\{r_N\}$ are sequences converging to possibly distinct elements of the partial compactification $\bar{\mathcal{Q}}$, but such that $d_\mathrm{Sym}(q_N,r_N) \to 0$ with respect to a common M\"obius normalization. Let $\phi_N$ be the type-preserving bijection between the singularities of $q_N$ and the singularities of $r_N$ that minimizes $\sum |\phi(z_{j,N}) - z_{j,N}|$ for this M\"obius normalization. Let $\{S_\ell(N)\}$ be a vanishing cluster for the sequence $\{q_N\}$. Suppose that as $N \to \infty$, $\frac{d_{\mathrm{Sym}(q_N,r_N)}}{\diam_\cc(S_{m}(N))} \to 0$ for every vanishing cluster $\{S_m(N)\}$ in $q_N$ with $S_{m}(N) \supsetneq S_{\ell}(N)$ for all $N$. Then $\phi_N(S_\ell(N))$ is a vanishing cluster for $r_N$.\end{proposition}

\noindent Proof: It suffices to consider the minimal $S_m$ that properly contains $S_\ell$. It is clear from the definition of vanishing cluster that for any sequence of sets of singularities $S(N)$, $S_m(N) \supseteq S(N) \supsetneq S_\ell(N)$ we must have $\diam_\cc(T(N)) \dot{\asymp} \diam_\cc(S_m(N))$, since otherwise we could pass to a subsequence along which $T(N)$ belonged to a vanishing cluster with $S_m(N) \supsetneq S(N) \supsetneq S_\ell(N)$. By minimality of $S_m$ this can't happen, since the nesting of the vanishing clusters is determined by the limit in $\bar{\mathcal{Q}}$. $\Box$

\begin{lemma}\label{LemmaPerturbEff} Suppose $\{q_N\}_{N = 1}^\infty = \{q((\lambda_N,z_{1,N},...,z_{r,N})\}_{N = 1}^\infty$ is a sequence in $QD(\mathcal{T}_{0,n})$, normalized in the sense of Notation \ref{NotaFixpoles} that converges in the sense of Definition \ref{DefCompactification}, such that no vanishing cluster contains a pole and every vanishing cluster contains a zero. Suppose that in addition, no vanishing cluster for the sequence $\{q_N\}$ has more than $k$ zeros (counted with multiplicity), $k \geq 0.$ Let $\{r_N\}$ be another sequence with the same M\"obius normalization, converging to a point in the partial compactification of Definition \ref{DefCompactification}, and assume that for each singularity type $m$ and each $p \in \cc$, the same number of singularities of type $m$ that limit to $p$ along the sequences $\{q_N\}$ and $\{r_N\}$ is the same. (In particular, $\{q_N\}$ and $\{r_N\}$ have the same limit in $QD(\mathcal{T}_{0,n}).$) Then $$d_{\mathrm{Euclidean}}(q_N,r_N) \dot{\succ} d_\mathrm{Sym}(q_N,r_N)^{2/(2+k)} ~\mathrm{as}~N \to \infty.$$\end{lemma}

\noindent Proof: The proof is divided into two steps. In the first step we fix a subsequence of potential counterexamples with specified combinatorics of the singularities belonging to vanishing clusters, and isolate some piece of the difference between $q_N$ and $r_N$. In the second we show that this piece causes a non-trivial change in period coordinates that is not cancelled by any other components.\\

\noindent STEP 1: Each maximal vanishing cluster of $\{q_N\}$ is also a maximal vanishing cluster for $\{r_N\}$, and we will only consider bijections that preserve maximal vanishing clusters, since $d_\mathrm{Sym}$ is only realized by such a bijection for large $N$. In fact, we can take this further. For each vanishing cluster $\{S_\ell(N)\}$ of $\{q_N\}$ we can pass to a subsequence whereby for each vanishing cluster $S_\ell(N)$, $$d_\mathrm{Sym}(q_N,r_N) \dot{*}_\ell \diam_\cc(S_\ell(N))$$ where $\dot{*}_\ell$ depends the choice of $\ell$ but is either $\dot{\prec}$, $\dot{\succ}$, or $\dot{\asymp}$.\\

Assume $$q_N = \frac{\lambda_N}{(z-A)(z-B)}\prod_{j = 1}^r (z-z_{j,N})^{e_j}dz^2, r_N = \frac{\mu_N}{(z-A)(z-B)}\prod_{j = 1}^r (z-w_{j,N})^{e_j}dz^2,$$

\noindent and let $\phi_N$ be the bijection realizing $d_\mathrm{Sym}(q_N,r_N)$ takes $z_j$ to $w_j$ for all $j$. We will describe how to break the perturbation realized by $\phi_N$ into components which we will call $v_{\ell,N}$ and define below.\\

\noindent Let $\{S_\ell(N)\}$ be either a vanishing cluster of the sequence $\{q_N\}$ or the entire set of singularities (besides $\infty$). Let the sequences $\{v_N\}$ and $\{v_{\ell,N}\}$ satisfy the following:\\

\begin{itemize}
\item $v_{\ell,N}$ is an $(r+1)$-tuple of complex numbers, corresponding to the coordinates $(\lambda,z_1,...,z_r)$ and $(\mu,w_1,...,w_r)$.
\item If $\{S_\ell(N)\}$ is a vanishing cluster of $\{q_N\}$, then the support of $v_{\ell,N}$ is in the entries corresponding to singularities in $S_\ell$.
\item If $S_k \subsetneq S_\ell$ is a vanishing cluster of $\{q_N\}$, then $v_{\ell,N}$ is constant on the entries corresponding to singularities in $D^\prime$. (All points in each proper subcluster move the same distance and in the same direction.)
\item For each $N$, $\sum\limits_{\ell} v_{\ell,N} = (\mu - \lambda, w_1 - z_1,...,w_s - z_s) =: v_N$.
\item If $S_\ell$ is a vanishing cluster and $v_{\ell,N} = (0,a_1,...,a_s)$ then $\sum\limits_{j = 1}^s a_je_j = 0$. (Together with other conditions, this means that $v_{\ell,N}$ does not move the center of $S_{\ell,N}$.)
\end{itemize}

\noindent The purpose of this notation, and the idea of the rest of the proof, is that for some $\ell,$ the component $v_{\ell,N}$ can be detected either by the periods of saddle connections internal to $S_\ell(N)$, or by periods internal to the largest $S_k(N) \supset S_\ell(N)$ with $\diam_\cc(S_k(N)) \dot{\prec} \|v_{\ell,N}\|$. Then we show that this does not get cancelled by other components of the perturbation.\\

\noindent We may pass to a subsequence with the following property: for any, $S_k$ and $S_\ell$, say, not necessarily distinct, the following converge (possibly to $0$ or $\infty$):

\begin{itemize}
\item $\frac{\diam_\cc(S_{\ell}(N))}{\diam_\cc(S_k(N)))}$
\item $\frac{\|v_{\ell,N}\|}{\diam_\cc S_{k,N}}$
\item $\frac{v_{\ell,N}}{\|v_N\|}$
\end{itemize}

\noindent In other words, the relative sizes of the various clusters and components of $v_N$ converge, and if $v_{\ell,N}$ which is not vanishingly small relative to any $v_{k,N}$, the direction of $v_{\ell,N}$ converges. Moreover, the sizes of perturbations of vanishing clusters converge relative to the sizes of all vanishing clusters.\\

\noindent Now, pick $S_\ell$ to be maximal subject to the following condition: Either $|v_{\ell,N}| \dot{\asymp} |v_N|$, or $\sum\limits_{S_l \subseteq S_\ell} |v_{l,N}| \dot{\asymp} |v_N|$ and $\sum\limits_{S_l \subseteq S_\ell} |v_{\ell,N}| \dot{\succ} \diam_\cc(S_{\ell}(N)).$\\

\noindent In particular, this means that either $S_\ell(N)$ is always the entire set of singularities, or $S_\ell(N)$ and $\phi_N(S_\ell(N))$ are both vanishing clusters of the sequence $\{r_N\},$ by Proposition \ref{PropStableClust}.\\

\noindent In both cases, Let $T_N(z) = a_N z + b_N$ be a sequence of M\"obius transformations fixing $\infty$ with the property that $T_N(S_\ell(N))$ converges (Hausdorff) to the limiting differential associated to $S_\ell$; we may assume $a_N \in \rr$. Then $T_N(\phi_N(S_\ell(N)))$ also converges to the singularity set of a cluster differential along a subsequence.\\

\noindent STEP 2: Estimating the perturbation in period coordinates.\\

\noindent Now we have two cases. In both cases, we can detect a change in Euclidean distance from the periods internal to $S_\ell(N)$ by proposition \ref{PropClusterProj}.\\

\noindent CASE 1: $\frac{|v_N|}{\diam_\cc(q_N)} \dot{\succ} 1.$

\noindent In this case, for all large $N$ we can define the projection map $F$ associated to the clusters $S_\ell(N)$ and $\phi_N(S_\ell(N))$ in Proposition \ref{PropClusterProj} on a connected open set containing both, by Proposition \ref{PropStableClust}. By Lemma \ref{LemAsympClust} and Corollary \ref{CoroClusterSizes} it follows that $$d_\mathrm{Euclidean}(q_N,r_N) \dot{\succ} d_\mathrm{Euclidean}(F_N(q_n),F_N(r_N)) \dot{\asymp} \max\{\diam_{q_N}(S_\ell(N)),\diam_{r_N}\phi_N(S_\ell(N))\}.$$

\noindent The last inequality can be justified as follows: First of all, the ratios between the factors $t_N$ in Lemma \ref{LemAsympClust} for $q_N$ and $r_N$ differ by a ratio of $1 + o(1)$. Now, if two cluster differentials $q,r$ have $d_\cc$-diameter at most $1$, and one has diameter $1$, and their distance with respect to $d_\cc$ is bounded below, then $d_\mathrm{Euclidean}(q,r)$ is also bounded below by compactness. Our last inequality then comes from Lemma \ref{LemAsympClust} and the fact that scalars induce similarities on spaces of cluster differentials.\\

\noindent Now, by Corollary \ref{CoroClusterSizes} and Lemma \ref{LemAsympClust} we see that if $\max\{\diam_{\cc}(q_N),\diam_\cc(r_N)\} = \epsilon_N,$ then $$\max\{\diam_{q_N}(S_\ell(N)),\diam_{r_N}(\phi_N(S_\ell(N)))\} \dot{\succ} \epsilon_N^{(m+2)/2}, ~ \mathrm{and}$$ $$\max\{\diam_{q_N}(S_\ell(N)),\diam_{r_N}(\phi_N(S_\ell(N)))\} \dot{\asymp} \epsilon_N^{(m+2)/2}$$ only if $q_N$ contains all $m$ singularities.\\

\noindent CASE 2: $\frac{|v_N|}{\diam_\cc(q_N)} \to 0.$

\noindent In this case we have $|v_{\ell,N}| \dot{\asymp} |v_N|.$ Again we choose $F_N$ to be the Lipschitz projection associated to $S_\ell(N)$. Assume the center of $|v_{\ell,N}|$ is always zero and let $t_N = \prod\limits_{z_{j,N} \notin S_\ell(M)} (-z_{j,N})^{e_{j,N}}.$ Assume $m$ is the number of zeros in $S_\ell{N}$ (each zero has multiplicity $1$, so we don't need to specify that we are counting with multiplicity.)

\noindent If $M_{\ell,N}$ is as in Corollary \ref{CoroLogmodelsprep} then the sequence of matrices $$\{t_N^{-1/2}\diam_\cc(S_{\ell,N})^{-m/2} M_{\ell,N}\}_{N = 1}^\infty$$ converges (by Corollary \ref{CoroLogmodels} and Proposition \ref{PropDerivScale}), and its kernel is the span of those vectors that are allowed to be $v_{k,N}$ with $k \neq \ell$. What we mean by this is that when we defined $v_{k,N}$, we required that it belong to a certain subspace of $\cc^{r+1}$ and all vectors in this subspace are in the kernel of the limiting matrix. The periods of saddle connections internal to proper subclusters of $S_\ell(N)$ span the cokernel of the limiting matrix. The limiting matrix $$\lim\limits_{N \to \infty} t_N^{-1/2}\diam_\cc(S_{\ell,N})^{-m/2} M_{\ell,N}$$ has full rank when restricted to the space of vectors that are allowed to be $v_{\ell,N}$.\\

\noindent We may move from $q_N$ to $r_N$ along a piecewise-smooth path $\beta_N$, so that the motion of each singularities $z_{j,N}$ to $\phi(z_{j,N})$ is a straight line segment with constant speed along each of two segments of $\beta_N$. On the first piece of $\beta_N$, the total displacement of the vectors $(\lambda_N,z_{1,N},...,z_{r,N})$ will be $v_{\ell,N}$ and on the second piece, the total displacement will be the remainder.\\

\noindent The singularities starting at $S_\ell(N)$ will be in the stratum of $S_\ell(N)$ for all but finitely many points along each segment. For any such sequence of paths, if we pick $s_N$ to be one point in this stratum along each $\beta_N$, and construct the matrix $\hat{M}_{\ell,N}$ analogously to the matrix $M_{\ell,N}$ at $s_N$, we can pass to a subsequence that converges in the sense of Notation \ref{NotaCompactConverge}. We also note that $\frac{\partial P_i}{\partial \lambda} = o(M_{\ell,N})$, since $\frac{\partial \log P_i}{\partial z_j}$ is not bounded for all $j$ but $\frac{\partial P_i}{\partial \lambda}$ is. The limit of the sequence of matrices must be the same, since it only depended on the limiting differential for $S_\ell$, which is the same for any subsequence of the sequence $\{s_N\}$. We therefore conclude that the matrices $M_{\ell,N}$ converge \emph{uniformly after rescaling by constants $c_N$ along $\beta_N$} in the sense that the matrix $\hat{M}_{\ell,N}$ differs from $M_{\ell,N}$ at $q_N$ by $o(|M_\ell,N|)$. In particular, this does not depend on which choices of saddle connections that degenerate to the proper edges of the limiting differential for $S_\ell$ we picked at each $s_\ell$. Moreover, the displacement of periods of $S_\ell$ is $M_{\ell,N}(v_{\ell,N}) + o(|M_{\ell,N}(v_N)|)$ along the first segment of $\beta_N$ and the displacement is $o(|M_{\ell,N}(v_N - v_{\ell,N})|)$ along the second segment of $\beta_N$.\\

\noindent This, combined with Proposition \ref{PropNoCancellation} tells us that $d_\mathrm{Euclidean}(F_N(q_N),F_N(r_N))$ is comparable to the size of the leading order approximation $M_{\ell,N}v_{\ell,N}$. By Corollary \ref{CoroFirstOrderOnClusterMagnitude} we have $$|M_{\ell,N} v_N| \dot{\asymp} |t_N|^{1/2}|s_N|^{m/2}|v_{\ell,N}| \dot{\succ} \diam_\cc(S_{\ell,N})^{(k-m)/2}\cdot[\diam_\cc(S_\ell(N))]^{m/2}|v_N| \cdot |v_N|^{(2+k)/2}.$$

\noindent We explain the sources of the inequalities: the size of $t_N$ is at least comparable to the product of the zeros outside of $S_\ell(N)$, since all poles remain bounded away from $0$. Since we only care about a multiplicative constant, we may ignore all zeros and poles that do not belong to a vanishing cluster containing $\{S_\ell(N)\}$ and there are at most $k-m$ additional zeros inside any such vanishing cluster, each of which is larger than $\diam_\cc(S_\ell).$ Hence $|t_N| \dot{\succ} |v_N|^{(k-m)/2}.$ The remaining inequalities are clear from $\diam_\cc(S_\ell(N)) \dot{\succ} |v_\ell,N| \dot{\asymp} |v_N|.$\\

\noindent Then by Proposition \ref{PropClusterProj} we get $$d_\mathrm{Euclidean}(q_N,r_N) \dot{\succ} d_\mathrm{Euclidean}(F_N(q_N),F_N(r_N)) \dot{\succ} |v_N|^{(k+2)/2}. \Box$$

\begin{lemma}\label{LemmaImplant} Fix a compact subset $H$ in the open upper half-plane $\mathbb{H}$. Fix a real number $C > 1$ and an integer $m \geq 0$. Consider the space $K$ of quadratic differentials in $\mathcal{T}_{0,2n}$ with the following restrictions:

\begin{itemize}
\item they are of the form $\frac{F(z)}{G(z)}dz^2$ with $F,G \in \rr[z]$, $\deg(F) = 2m$, $\deg(G) = 2m - 3$.
\item $G$ may have simple roots at $i$ and $-i$, but all other roots of $G$ are real. (Either way, $G$ has real coefficients.) The roots of $G$ are the marked points of the underlying Riemann surface.
\item $G$ has no repeated roots, and the distance between any two consecutive real roots of $G$ is between $C$ and $1/C$.
\item $G$ has a root in $[-C,C]$ if $i$ and $-i$ are roots of $G$; otherwise, $0$ and $1$ are roots of $G$.
\item No root of $F(z)$ is real, and one root from each conjugate pair belongs to $H$.\\
\end{itemize}

Assume $q_1 = \frac{F_1(z)}{G_1(z)}dz^2$ and $q_2 = \frac{F_2(z)}{G_2(z)}dz^2$ are elements of $K$. Let $\phi$ range over all conjugation-invariant bijections of the roots of $F_1$ with the roots of $F_2$ that fix $i$ if $i$ is a root of $G$ and let $\psi$ be the bijection of $\psi$ of the roots of $G_1$ with the roots of $G_2$ that preserves the order of the real roots and fixes $i$ and $-i$, if they are roots of $G_1$ and $G_2$. Let

$$\epsilon = \min\limits_{\phi}\max\{\phi(w) - w| F_1(w) = 0\}\cup \{\psi(w) - w| G_1(w) = 0\} = $$

Then there are constants $M$, $\epsilon_0 > 0$ depending only on $K$, such that the following hold:\\

If $\pm i$ are not roots of $G$, then $d_\mathrm{Euclidean}(q_1,q_2) \dot{\succ} \epsilon^{(m+2)/2}.$\\

If $\pm i$ are roots of $G$, then $d_\mathrm{Euclidean}(q_1,q_2) \dot{\succ} \epsilon^{m+1}.$
\end{lemma}

\noindent Proof: First we will explain how both conclusions follow from Lemma \ref{LemmaPerturbEff}. If $\pm i$ are roots of $G_1$ and $G_2$ the we can take a double cover of each that preserves the real line and is branched at $\pm i$ via $z = h(w) = \frac{w^2 - 1}{2w}.$ Then $\pm i$ are marked points, but not poles on the pullbacks $h^*(q_1),h^*(q_2)$. Then $\iota^*(q)$ is a conjugation-invariant quadratic differential on $\hat{\cc}$ whose poles are all simple and occur on the projective real line. There are uniform upper and lower bounds on the distance between real poles, and the zeros with positive real part are in a fixed compact part of the upper half-plane.\\

\noindent Now, suppose that $z_1$ is a root of $F_1$ and $z_2$ is a root of $F_2$, and $|z_1 - z_2| = \epsilon.$ Then if $h(w_1) = z_1$ and $h(w_2) = z_2$, we have $$|w_1 - w_2| = |h(w_1) - h(w_2)| \cdot \left|\frac{2w_1w_2}{1 + w_1w_2}\right| \dot{\succ} \epsilon.$$

\noindent We have a similar estimate for the inverse images of the roots of $G_1$.\\

\noindent Now, we can apply Lemma \ref{LemmaPerturbEff} to both. $\Box$

\section{Building the Quasiconformal Map}

\noindent In this section, we build a quasiconformal map with the desired properties to prove Theorem \ref{TheoremMain}.

\subsection{Cutting into Triangles and Nearly Regular Right Polygons}

\begin{definition}\label{DefNRRP} A \emph{nearly regular right polygon}, or \emph{NRRP}, is a subset $P$ of a half-translation surface $(X,q)$ with the following properties:

\begin{itemize}
\item $P$ is homeomorphic to a closed disk
\item $P$ contains at least one singularity in its interior, and at most one pole, but has no singularities on its boundary
\item $\partial P$ is piecewise geodesic and all pieces are vertical or horizontal segments
\item every interior angle of $\partial P$ is $\pi/2$
\item If $P$ is doubled along its boundary to form a sphere $P \cup \overline{P}$, where $\overline{P}$ is the Schwarz reflection of $P$, then there is a choice of holomorphic coordinates, $\partial P$ with the real projective line, and the quadratic differential on $P \cup \partial P$ has one of the two forms in Lemma \ref{LemmaImplant}.
\end{itemize}

Let $R$ be the minimum $q$-distance from $\partial P$ to a singularity in $P \setminus \partial P$. We call $R$ the \emph{radius} of the NRRP.
\end{definition}

\noindent On a quadratic differential it is always possible to find an NRRP around each singularity with all sides of the same length. If we change the differential by a small enough amount, and in so doing break up such a singularity into multiple singularities, we can still find an NRRP that stays Hausdorff close to our original NRRP. The motivation for this definition is that if one wants to build a quasiconformal map between two half-translation surfaces, piecewise-affine maps which are affine on some triangulation will tend to have large dilatation even when two Riemann surfaces are nearly conformal, if the triangulation is nearly degenerate. Given an edgewise-linear map between boundaries of similarly shaped NRRP's in two quadratic differentials that are close in $QD(\mathcal{T}_{g,n})$, we can extend to a quasiconformal map NRRP's, and explicitly estimate its quasiconformal dilatation.

\begin{proposition} \label{PropDecompControl} On any compact subset $K$ of $QD(\mathcal{T}_{g,n})$, there are constants $C = C(K) > 0$ and $\delta = \delta(K)$ such that for every surface $X$ in $K$, we can partition the set of singularities of $X$ into $\delta$-clusters and singletons (for some $\delta \in (0,1)$) and associate to each $\delta$-cluster or singleton an NRRP containing it such that the following hold:

\begin{itemize}
\item Each NRRP contains only the singularities associated to it.
\item Each NRRP has radius and all boundary side lengths at least $C(K)$.
\item If we consider all lifts of the NRRPs in the universal cover of the underlying compact Riemann surface, the distance between any two NRRPs is more than twice as large as the side length of any NRRP.
\end{itemize}
\end{proposition}

\noindent Proof: It is clearly possible to choose a constant $C(X)$ at each $X \in K$ such that the properties hold for $C(X)$, by picking a small NRRP around each singularity. (NRRPs of a given radius persist under sufficiently small perturbation). For each $X$, we can find a neighborhood of $X$ the conditions all hold with the constant $C(X)/2$ instead of $C(X)$, and $K$ is finitely covered by such neighborhoods. It follows that we can take $C(K)$ to be half the minimum value of $C(X)$ used in a finite subcover. $\Box$\\

\noindent Let $Y$ be the closure of the complement of the union of a system of NRRPs satisfying the hypotheses of Proposition \ref{PropDecompControl} in a half-translation surface $X$, and assume every singularity of $X$ is contained in some NRRP. If we double $Y$ along $\partial Y$, the foliations by vertical and horizontal segments extend by reflection, and $Y$ acquires the structure of a quadratic differential, with singularities of cone angle $3\pi$ at the vertices of $\partial Y$.\\

\noindent Moreover, if we we take the Delaunay triangulation of the resulting surface, the midpoint of any boundary edge is strictly closer to the endpoints of that edge than to any other singularity, so all edges of $\partial Y$ belong to the Delaunay triangulation of the double cover.\\

\noindent Therefore, \emph{it makes sense simply to speak of the Delaunay triangulation of} $Y$. That is to say, if $\Sigma$ is the set of vertices of $Y$ (which are all on the boundary of $Y$), then we have proved:

\begin{proposition} $Y$ has a triangulation with the property that the circumcenter of each triangle is the boundary of a standard Euclidean disk of radius equal to the distance to any vertex of the triangle. The boundary edges of $Y$ are all edges. $\Box$ \end{proposition}

\noindent Now, we describe how to build a quasi-conformal map between two nearby quadratic differentials, provided they are within some distance $c(K)$ that depends only on $K$.\\

\noindent Let $X_1$ and $X_2$ be two such Riemann surfaces, and let $Y_1$ be the complement of a system of NRRPs for $X_1$ satisfying the hypotheses of Proposition \ref{PropDecompControl}.\\

\noindent The circumradii of Delaunay triangles on any half-translation surface are bounded above by twice the diameter of the surface, which is bounded on $K$. The side lengths of $Y_1$ are bounded below by $C(K)$. Thus by the law of sines, the angles of Delaunay triangles are bounded below. So we could have picked $c(K)$ small enough that when we change all of the period coordinates by at most $c(K)$, all angles of the Delaunay triangulation remain bounded away from zero. So the Delaunay triangulation for $Y$ is uniformly bounded away from degeneration, i.e. all of the side lengths and angles are uniformly bounded away from $0$. If $d_\mathrm{Euclidean}(X_1,X_2) < c(K)$, then along a path from $X_1$ to $X_2$ we can make a choice of NRRP decomposition so that the radii, edges of NRRPs, and edges of Delaunay triangles that start as edges of $Y_1$ vary in a Lipschitz manner with respect to arc length, and the Lipschitz constant depends only on $K$. Let $Y_2$ be the complement of the NRRPs in $X_2$, as chosen along the path.\\

\noindent By Proposition \ref{PropMultStable} the triangulation of $Y_1$ persists to a triangulation $Y_2$, and this determines an affine map on triangles of $Y_1$ to a triangulation of $Y_2$. Since no edge gets too short, and no angle gets too close to $0$, it follows that if the length of the path is $\epsilon$, then as $\epsilon \to 0$ this part of the map is $1 + O(\epsilon)$-quasiconformal on these triangles, and the implied constants in the $O(\epsilon)$ depend only on $K$.

\subsection{The Beurling-Ahlfors Extension}

\noindent The final step is to show how to build a quasiconformal map between NRRP's and estimate its dilatation. We will do this by means of a Beurling-Ahlfors extension. Given an orientation preserving self-homeomorphism $h$ of $\rr^+$, we have the following:

\begin{definition}The \emph{Beurling-Ahlfors extension} of $h$ to the upper half plane \emph{with parameter} $r > 0$ is the function $f_r$ given by\end{definition}

$$f_r(x + iy) = \frac{1}{2}\int_0^1 [h(x + yt) + h(x - yt)]dt + \frac{ir}{2}\int_0^1 [h(x + yt) - h(x - yt)]dt.$$

\begin{definition}Let $\rho \in [1,\infty).$ We say that an orientation preserving homeomorphism $h: \rr \to \rr$ is $\rho$-\emph{quasisymmetric} on $\rr$ if for all $x\in \rr, t > 0$ we have $$\frac{1}{\rho} \leq \frac{h(x + t) - h(x)}{h(x) - h(x - t)} \leq \rho.$$\end{definition}

\noindent Beurling and Ahlfors \cite{BA} proved

\begin{theorem}If $h$ is $\rho$-quasisymmetric then for some explicit choice of $r$, the Ahlfors Beurling-Ahlfors extension $f_r$ is $\rho^2$-quasiconformal. Moreover, $r = 2 + O(\rho - 1)$ as $\rho \to 1^+$.\end{theorem}

\noindent Note that $r = 2$ extends a M\"obius transformation of $\mathbb{RP}^1$ to a M\"obius transformation of $\mathbb{CP}^1$. Any choice of $r$ yields a quasiconformal extension, but possibly with worse dilatation.\\

\noindent In the event that a NRRP does not contain a marked point, the Beurling-Ahlfors extension, together with the piecewise-affine maps with compatible boundary conditions completes our construction of a quasiconformal map. In the case with marked points, we will have to compose with a quasi-conformal map on each NRRP which is trivial on the boundary and moves the marked point to the correct location. For this we will uniformize our map to the upper half-plane and apply the unique $\rr$-linear map on the complex plane fixing the real line pointwise and taking our marked point to its desired location.\\

\noindent We have established uniform control over the NRRP's. If the only interior singularity is a pole, then our map of NRRP's is in fact already an isometry in the singular metric, hence conformal. In all other cases, we have at least four poles on the boundary, and we can assume three of them are $0,1,\infty$. Subject to this normalization, if we double each NRRP along its boundary, the singularities in the interior of each NRRP can be chosen to belong to some fixed compact subset of the upper half-plane.\\

\noindent This motivates the following:\\

\begin{lemma}Suppose that $P_1,P_2$ are NRRP's. Let $h_q: \partial P_1 \to \partial P_2$ be a piecewise affine identification of sides with respect to the singular metrics on $P_1$ and $P_2$, affine on each pair of corresponding sides and assume that when we uniformize each of $P_1$ and $P_2$ to the upper half-plane, three pairs of corresponding corners are sent to $0,1,\infty$, and the distance between any two finite singularities on the boundary is between $1/B$ and $B$, $B > 1$, and all interior singularities belong to some fixed compact subset $K$ of the open upper half-plane. For all sufficiently small $\epsilon$, if the singularities of $P_1$ and $P_2$ can be put in bijection so that corresponding singularities are distance at most $\epsilon$ apart, then the quasisymmetry constant $\rho$ of the boundary map $h: \rr \to \rr$ induced by $h_q$ and our chosen uniformizations is at most $C\epsilon$, for some $C$ depending only on $K$, the number of sides of $P_1$ and $P_2$, and $B$.\end{lemma}

\noindent Proof: The boundary map $h$ is continuous, and it is differentiable except possibly at points corresponding to vertices of the NRRP. It is enough to show that the derivative of the boundary map satisfies $h^\prime(x) - 1 = O(\epsilon)$ uniformly on the complement of these points. \\

\noindent Suppose that one pair of corresponding sides is $h([a_1,b_1]) = [a_2,b_2]$ where $a_1,a_2,b_1,b_2$ are all finite. (We will deal with the intervals with an endpoint at $\infty$ later.) Assume the quadratic differentials are $\frac{p_i(z)}{(z - a_i)(z - b_i)q_i(z)}dz^2$ on $[a_i,b_i]$ where $p_i$ and $q_i$ are polynomials. Then for $a_1 < c < b_1$ our boundary map $h$ is defined so that

$$\frac{\int_{a_1}^c \left(\frac{p_1(z)}{(z - a_1)(z - b_1)q_1(z)}\right)^{1/2}dz}{\int_{a_1}^{b_1} \left(\frac{p_i(z)}{(z - a_i)(z - b_i)q_i(z)}\right)^{1/2}dz} = \frac{\int_{a_2}^{h(c)} \left(\frac{p_2(z)}{(z - a_2)(z - b_2)q_2(z)}\right)^{1/2}dz}{\int_{a_2}^{b_2} \left(\frac{p_2(z)}{(z - a_2)(z - b_2)q_2(z)}\right)^{1/2}dz}.$$

\noindent We differentiate both sides of the above equation in $c$ and solve for $h^\prime(c)$:

$$h^\prime(c) = \frac{\int_{a_2}^{b_2} \left(\frac{p_2(z)}{(z - a_2)(z - b_2)q_2(z)}\right)^{1/2}dz}{\int_{a_1}^{b_1} \left(\frac{p_1(z)}{(z - a_1)(z - b_1)q_1(z)}\right)^{1/2}dz}\left(\frac{p_1(c)q_2(h(c))(h(c) - a_2)(h(c) - b_2)}{p_2(h(c))q_1(c)(c - a_1)(c - b_1)}\right)^{1/2}.$$

\noindent By applying changes of coordinates consisting of translations by real numbers to each NRRP, we can make a few simplifying assumptions: At the cost of increasing $\epsilon$ to $2\epsilon$ we can do a translation to assume $a_1 = a_2 = 0$. After this simplifying assumption, our argument can be rephrased: let $$f_i(z) = \frac{\left(\frac{p_i(z)}{z(z - b_i)q_i(z)}\right)^{1/2}}{\int_{0}^{b_i} \left(\frac{p_i(t)}{t(t - b_i)q_i(t)}\right)^{1/2}dt}.$$ Then $h^\prime(c) = \frac{f_1(c)}{f_2(h(c))}.$ So we need to show $\log \left[\frac{f_1(c)}{f_2(h(c))}\right] = O(\epsilon).$

\noindent We will actually do this under the assumption that $0 < c \leq 2b_i/3.$ A basically identical argument will cover the cases $b_i/3 < c < b$, and a final argument will extend this to cover the two pairs of sides of $P_1$ and $P_2$ that have $\infty$ as an endpoint after we uniformize.\\

\noindent In what follows, $C$ will be used to denote various positive constants that depend only on choices we have already made - its meaning may vary from one line to the next. The following will be immediate from our compactness assumptions. To simplify matters, we first prove 

\begin{claim} $|\log(h(c)/c)| \leq C\epsilon$.\end{claim}

\noindent The idea is that the $f_i$ are equal to $z^{-1/2}$ times two very similar functions. Recall that $h$ is defined so that $\int_0^c f_1(z)dz = \int_0^{h(c)}f_2(z)dz.$\\

\noindent Define $g_i(z) = z^{1/2}f_i(z).$\\

\noindent We know the following:

\begin{itemize}
\item The $g_i$ are either both purely imaginary or both purely real, and never change sign.
\item $-C \leq \log|g_i(z)| \leq C$ and $|g_i^\prime(z)| \leq C$ for $z \in [0,3b_i/4]$.
\item $|\log(f_1(z)/f_2(z)) = |\log (g_1(z)/g_2(z))| < C\epsilon$ for $z \in (0,3b_i/4]$.
\item $-C \leq \log [z^{1/2}\int_0^z f_i(t)dt] \leq C$
\item $\frac{d}{dz}\log\left|\int_0^z t^{-1/2}g_i(t)dt\right| \geq z^{-1/2}/(Cz^{1/2}) = 1/(Cz)$ for $z \in (0,3b_i/4]$.
\end{itemize}

\noindent From the last observation we see that the amount that we need to perturb $z$ to change the value of $\log \int_0^z t^{-1/2}g_i(t)dt$ by $\epsilon$ is at most $Cz\epsilon$, for all sufficiently small $\epsilon$. In particular, for $i = 2$, we can start with $\int_0^c f_2(t)dt$, and by moving the upper limit of integration most $C\epsilon c$ away from $c$ we can adjust the value of the integral of $f_2$ by a factor of at least $1 + C\epsilon$. The intermediate value theorem implies that the adjustment of the upper limit of integration needed to get from $c$ to $h(c)$ by is at most $C\epsilon c$. This establishes the claim.\\

\noindent Now that we we have the claim, the main case of the lemma follows quite easily: the logarithmic derivatives of all factors of $f_i$ are bounded by some constant $C$ on $[c,h(c)]$ or $[h(c),c]$, and $\log(f_1/f_2)$ is uniformly close to $(0,3b/4)$. So the claim, plus boundedness of logarithmic derivatives of the remaining factors gives us $g_1(c)/g_2(h(c)) = 1 + O(\epsilon)$, and by the claim $c/h(c) = 1 + O(\epsilon)$. This proves the main case of the lemma. Interchanging the roles of $a_1$ and $b_1$ handles the case in which $c$ is in the right half of the interval $[a_1,b_1]$.\\

\noindent Now we must handle the infinite intervals for $1 < c < \infty$, we can conjugate $h(z)$ by the map $z \to 1/z$ to a function $H(c)$ satisfying an estimate of the type we originally had. The lemma and its proof hold for $H$; in particular $H(1/c)/(1/c) = 1 + O(\epsilon)$. Finally, we have $$h^\prime(c) = \frac{d}{dc}\left(\frac{1}{H(1/c)}\right) = H^{\prime}(1/c) \frac{1}{c^2 H(1/c)^2} = 1 + O(\epsilon).$$

\noindent This completes the proof. $\Box$\\

\noindent Proof of Theorem \ref{TheoremMain}: \noindent Now, we describe how to build a quasi-conformal map between two nearby quadratic differentials, provided they are within some distance $c(K)$ that depends only on a compact set $K \subset QD(\mathcal{T}_{g,n})$.\\

\noindent Let $X_1$ and $X_2$ be two such Riemann surfaces, and let $Y_1$ be the complement of a system of NRRPs for $X_1$ satisfying the hypotheses of Proposition \ref{PropDecompControl}.\\

\noindent The circumradii of Delaunay triangles on any half-translation surface are bounded above by twice the diameter of the surface, which is bounded on $K$. The side lengths of $Y_1$ are bounded below by $C(K)$. Thus by the law of sines, the angles of Delaunay triangles are bounded below. So we could have picked $c(K)$ small enough that when we change all of the period coordinates by at most $c(K)$, all angles of the Delaunay triangulation remain bounded away from zero. So the Delaunay triangulation for $Y$ is uniformly bounded away from degeneration, i.e. all of the side lengths and angles are uniformly bounded away from $0$. If $d_\mathrm{Euclidean}(X_1,X_2) < c(K)$, then along a path from $X_1$ to $X_2$ we can make a choice of NRRP decomposition so that the radii, edges of NRRPs, and edges of Delaunay triangles that start as edges of $Y_1$ vary in a Lipschitz manner with respect to arc length, and the Lipschitz constant depends only on $K$. Let $Y_2$ be the complement of the NRRPs in $X_2$, as chosen along the path.\\

\noindent Since the triangulation persists, this determines an affine map on triangles of $Y_1$ to a triangulation of $Y_2$. Since no edge gets too short, and no angle gets too close to $0$, it follows that if the length of the path is $\epsilon$, then as $\epsilon \to 0$ this part of the map is $1 + O(\epsilon)$-quasiconformal on these triangles, and the implied constants in the $O(\epsilon)$ depend only on $K$. We need only extend this map to the interiors of the disks.\\

\noindent Now, the NRRP's themselves are represented by quadratic differentials on disks, which can be uniformized to the closure of the upper half-plane in $\hat{\cc}$; we may assume three pairs of corresponding vertices map to $0,1,\infty$ under this uniformization. The period coordinates of the quadratic differentials on these disks differ by $O(\epsilon)$, and so the locations of the zeros and poles in uniformized coordinates differ by $O(\epsilon^{\alpha_{g,n}})$, where $\alpha_{g,n} = {2/[2 + a_n(4g-4 + n)]}$, by Lemma \ref{LemmaPerturbEff}. This gives us the estimate we need for the dilatation of the map on the NRRP's.\\

\noindent Our final task is to modify the map, if necessary, to make sure that any marked point gets mapped to the correct marked point. However, since the marked point (there is at most one) is in a fixed compact subset $K$ of the upper half-plane, our boundary map differs from the identity on any compact subset of the real line by $O(\epsilon^{\alpha_{g,n}})$, it follows that the Beurling-Ahlfors extension moves points in $K$ by $O(\epsilon^{\alpha_{g,n}})$ as well. The $\rr$-linear map that is the identity on the boundary and moves a fixed point in $K$ by $O(\epsilon^{\alpha_{g,n}})$ also has constant dilatation that is $1 + O(\epsilon^{\alpha_g,n})$. We have constructed a homeomorphism that satisfies the dilatation bound at almost every point, as its first partial derivatives are defined off of a measure $0$ set, they are clearly bounded away from every point except possibly the singularities and vertices of the NRRP's, since there is a conformal metric in which the partial derivatives converge piecewise, which shows that if we delete a finite set the first partial derivatives are in $L_{loc}^2.$ To finish, we simply recall that a homeomorphism which is $M$-quasiconformal on the complement of a finite set is $M$-quasiconformal, see e.g. \cite{Hubbard}. $\Box$

\appendix

\section{Local Finiteness of Period Coordinate Systems}

In this appendix we show the following, leaving the proof to the end of the section:

\begin{proposition}\label{PropLocallyFinite} There are only finitely many systems of period coordinates represented by saddle connections of length at most $D$ in any compact subset $K$ of $QD^1(\mathcal{T}_{g,n})$. \end{proposition}

\noindent To this end, we will show that every such period coordinate system is related to a basis of the edges of the Delaunay triangulation by one of finitely many transition matrices, and there are only finitely many choices Delaunay triangulation on $K$. We will always assume that our triangulations are \emph{labelled} and \emph{marked}, that is, the vertices and edges are distinguishable, we know which marked points (if any) on the base surface correspond to which singularities, and we know which homotopy classes of curves on the surface are represented by which homotopy classes of curves on the graph, and finally, we will also include in the data of the Delaunay triangulation the sign of the slope of each saddle connection.

\begin{proposition}Up to the action of the mapping class group, there are only finitely many triangulations that can occur as Delaunay triangulations in $QD^1(\mathcal{T}_{g,n})$.\end{proposition}

\noindent Proof: The number of triangles is a function only of the stratum, since the sum of the angles of a triangle is $\pi$ and the sum of the angles of all triangles is the sum of the cone angles of the singularities. (There is one triangle for each vertical separatrix emanating from a singularity.) $\Box$

\begin{proposition}Only finitely many Delaunay triangulations occur in any compact $K \subset QD^1(\mathcal{T}_{g,n})$.\end{proposition}

\noindent Proof: If not, then infinitely many such triangulations occur in a single orbit of the $\Mod(S_{g,n}).$ Given a fixed Delaunay triangulation $T$, The space of complex lengths we can assign to the edges and still have a unit area quadratic differential with systole $> \epsilon$ and such that the triangulation remains Delaunay (but possibly becomes degenerate) is compact. (If there are no short simple closed curves, Mumford's compactness criterion implies the diameter of the surface, and hence the diameter of each Delaunay triangle, is bounded. The inequalities that guarantee that the triangulation be a Delaunay triangulation are closed conditions.) This contradicts the proper discontinuity of the action of $\Mod(S_g,n)$ on $QD^1(\mathcal{T}_{g,n}).$\\

\begin{proposition}There is a uniform bound on the number of saddle connections of length at most $D$ on $K$.\end{proposition}

\noindent Proof: We first reduce to the case in which $q \in K$ has no marked points. If $q$ has marked points, we simply find $q^\prime$ belonging to a compact $K^\prime$ in a moduli space of higher genus surfaces, by taking any branched cover branched over all marked points, and possibly one other point that does not belong to any saddle connections, as close as we like to the point on $q$ that is the greatest distance from any singularity, and give it the metric associated to the quadratic differential that is half the pullback of $q$ (This is so that the area is 1). The resulting surface is an Alexandrov non-positively curved space since it is locally Euclidean with all cone angles greater than $2\pi$. Since there is a lower bound on the length of the shortest simple closed curve it belongs to a compact part of some $\mathcal{TQD}_{g^\prime}$, and every saddle connection of $q$ is the image of a saddle connection of the branched cover.\\

\noindent Now, to deal with surfaces without poles, we first claim that the fundamental group of $S_g$ with any fixed word metric is $(k,c)$-quasi-isometric to the universal covers $\tilde{(X,q)}$ of $(X,q) \in K$ via a maps $\phi_(X,q)(h)$ that send $h \in \pi_1(K)$ to $h(x,q)$, where $x \in X$ is chosen in a continuous section of the universal disk bundle over $QD^1(\mathcal{T}_g)$, i.e. the universal cover of the surface bundle over $QD(\mathcal{T}_{g,n})$. The uniformity of quasi-isometries then follows from finiteness of the set of elements $\{h \subset \pi_1(S_g): \exists \in K ~\mathrm{s.t.}~ d_{\tilde{q}}(x,hx) < R\}$.\\

\noindent Any saddle connection can be modified to become a simple closed curve by changing one of its endpoints, at the cost of changing its length by a uniformly bounded amount (the diameter of the surface.) A lower bound on the resulting simple closed curve then implies, by quasi-isometry with uniform constants, that the resulting simple closed curve came from a finite list of conjugacy classes in the fundamental group of $S$.\\

\noindent Second, we claim the following: if $q_n \in K$ all have the same Delaunay triangulation, and $q_n \to q$, then the common Delaunay triangulation of the $q_n$ is a possibly degenerate Delaunay triangulation of $q$. This follows because the conditions for a triangulation to be Delaunay are nonstrict inequalities in the absolute values of period coordinates, which persist under taking limits. It follows that there are only finitely many period coordinate systems that can be Delaunay triangulations in a neighborhood of each point in $K$, because there are only finitely many possibly degenerate Delaunay triangulations at each point, and only finitely many ways to undegenerate each (since only finitely many labelled ribbon graphs with an upper bound on number of edges can contract to a given ribbon graph by contracting edges). Thus, a simple compactness argument tells us that there are only finitely many Delaunay triangulations of surfaces in $K$. For each subset of $K$ sharing a common Delaunay triangulation we can fix a basis of period coordinates coming from edges of the triangulation, and we can write the period of every homotopy class of geodesic arc whose endpoints are singularities in terms of this basis.\\

\noindent By basically the same argument as above, only finitely many homotopy classes of geodesic arcs beginning and ending at saddle connections can have length $<D$. Finally, we claim that the saddle connections in a fixed homotopy class can only have finitely many representations in terms of any fixed basis for $H_1^{\mathrm{odd}}(\tilde{X},\tilde{\Sigma)}$ consisting of saddle connections that are edges of the $L^\infty$  triangulation, if we force $(X,q)$ to remain in $K$. Indeed, by uniform quasi-isometry of the fundamental group with $\tilde{X}$ (endowed with the pulled back $q$ metric), such a saddle connection can only pass through a fixed finite list of Delaunay triangles, and it cannot leave and then enter a triangle, since  triangles are geodesically convex in the $q$-metrics. Thus there is a finite list of sequences of triangles a saddle connection can pass through, and by developing this chain of triangles in the plane, we can write its period as one of finitely many sums of sides of periods of the triangles.$\Box$

\begin{corollary}\label{CoroFinCoordSys} For each compact $K \subset \mathcal{T}_{g,n}$ there is a finite set of triangulations $\{T_i\}_i \in I$ such that all quadratic differentials lying over $K$ belong to a $T_i$-convex set for some $i \in I$.\end{corollary}

\noindent Proof: the $q$-metric on $(X,q)$ is obtained by pulling back the Euclidean metric on $\rr^2$ from a collection of charts, but one could just as easily have pulled back the $L^{\infty}$-metric $d(a+bi,c+di) = \max(|a-c|,|b-d|)$ to form an $L^\infty$ $q$-metric on $X$. We may consider Delaunay triangulations with respect to the $L^\infty$ flat metric, following \cite{FG}. Every triangle in the $L^\infty$ Delaunay triangulation has its vertices on the boundary of an open square that is maximal with respect to the conditions that its lifts to the universal cover are embedded and that its edges be vertical and horizontal. If a triangulation is the $L^\infty$ triangulation for a set of surfaces in $QD(\mathcal{T}_{g,n})$ that has nonzero Masur-Veech measure, then no edge of the triangulation is always horizontal or always vertical. So we may further restrict each saddle connection in the triangulation to have non-positive slope or non-negative slope, and include that in the combinatorial data of the triangulation. We may also disregard any triangulations that arise only on sets of measure zero.\\

\noindent We may now claim the conditions that a triangulation be $L^\infty$ Delaunay are given by (nonstrict) linear inequalitites in the period coordinates. Moreover, no triangle can have three edges of positive slope or three edges of negative slope. It follows that if we orient all edges so that the imaginary parts of periods are non-negative, then for each $L^\infty$  triangulation $T$, the space of surfaces for which $T$ is an $L^\infty$  triangulation is $T$-convex.\\

\noindent Moreover, for each triangle, there is a point $p$ that is equidistant from the three vertices in the $L^\infty$-metric, and those three vertices are the nearest three vertices to $p$ in the $L^\infty$ metric (possibly along with some others, if the $L^\infty$  triangulation is not unique). If this distance is $d$ then they do not all have positive slope or all negative slope, since they touch at least 3 sides of a $2d \times 2d$ square with vertical and horizontal sides centered at $p$. This imposes a uniform bound on the lengths of saddle connections in the $L^\infty$ Delaunay triangulation of $2\sqrt{2}$ times the diameter of the surface. Since it is obviously enough to consider unit area surfaces lying over $K$, which have bounded diameter, this imposes a bound on the number of possible $L^\infty$ Delaunay triangulations. $\Box$\\

\noindent Proof of Proposition \ref{PropLocallyFinite}: There are locally finitely many $L^\infty$-Delaunay triangulations; for each such triangulation we can locally fix names of the singularities. The collection of all singularities in the universal cover of $(X,q)$ is quasi-isometric to the fundamental group of the compact surface $\bar{X}$, uniformly for all $(X,q) \in K$. In particular, there are finitely many ways to choose the names of the endpoints of a saddle connection of length at most $L$ in the metric universal cover of $(\bar{X},q)$, up to the action of $\pi_1(\bar{X})$, and therefore only finitely many systems of saddle connections of length at most $L$ for each $L^\infty$ Delaunay triangulation that occurs in $K$. $\Box$

\section{\texorpdfstring{$\delta$}{d}-clusters and the Euclidean Metric}

\noindent First, we would like to establish a basic fact about systems of period coordinates persisting under perturbation. This will be useful throughout:\\

\begin{proposition}\label{PropMultStablePrep}Let $\mathcal{Q}(\mu)$ be a stratum of cluster differential. Then there is a constant $\epsilon > 0$, depending only on $\mu$, such that whenever $\Gamma$ is a length-minimizing system of saddle connections on the vertex set of singularities of $q \in \mathcal{Q}(\mu)$, there is a holomorphic coordinate chart $U$ given by the logs of the periods $\log(P_1),...,\log(P_r)$ of saddle connections in $\Gamma$, then $U$ contains an $\epsilon$-ball about $q$, in which the coordinates are logs of periods of saddle connections.\end{proposition}

\noindent Proof: First, we reduce to the case in which the stratum has no poles, except for a higher order pole at $\infty$. In the case of a cluster differential there is a length-minimizing tree of saddle connections $\Gamma$, and we claim that the length-minimzing tree $\tilde{\Gamma}$ for the double cover $(\tilde{X},\tilde{q})$ of $(X,q)$ branched only over the pole contains the lifts of all the edges of the base graph $\Gamma$, and therefore a small perturbation of the logs of the periods of $\Gamma$ arises as the quotient of a perturbation of the logs of the periods of $\tilde{\Gamma}$ which descends. The proof is by following greedy algorithm for length-minimizing spanning trees: repeatedly pick the shortest edge that does not form a cycle with previous edges. Every saddle connection of $(\tilde{X},\tilde{q})$ is a lift of a saddle connection of $(X,q)$. Following the greedy algorithm to build $\tilde{\Gamma}$ will thus start with a lift of the shortest edge of $\Gamma$ and then the other lift; this does not form a cycle because the double cover of the graph $\Gamma$ branched over the pole is still a tree. Now, continuing inductively, assume that the first $2k$ edges of $\tilde{\Gamma}$ are the lifts of the first $k$ edges of $\Gamma$. The shortest saddle connection of $\Gamma$ that has not been picked has lifts which do not form cycles with the previously picked edges of $\tilde{\Gamma}$, because the double cover of $\Gamma$ is a tree. However, any shorter edge would project down to an edge $e$, which, when added to the first $k$ edges of $\Gamma$ creates a cycle. Adding the two lifts of $e$ to $\Gamma$ and $\Gamma^\prime$ would therefore produce a graph with first homology group of rank $2$, since it would be obtained by taking two disjoint copies of a graph with first integral homology group $\zz$ and identifying a pair of vertices. Since deleting one edge of a graph can only decrease the rank of the first homology group by $1$ it follows that adding just one lift of $e$ would create a cycle with the first $2k$ edges of $\Gamma$, so therefore the lift of the saddle connection $e$ is not available to pick as the next edge of $\tilde{\Gamma}$. It therefore follows that taking the two lifts of the $k+1^{\mathrm{st}}$ edge of $\Gamma$ can be chosen by the greedy algorithm, and it is possible to choose one and then the other as the next two edges of $\tilde{\Gamma}$ while following the greedy algorithm. By induction, our claim follows, as does the reduction to the case with no poles.\\

\noindent Now, suppose the proposition is false. We can find a sequence of counterexamples $\{q_m\}$ would converge in the sense of Definition \ref{DefCompactification}, and the graph $\Gamma$, as well as which collections vertices correspond to vanishing clusters, are the same for all $m$. There would be perturbation vectors $h_m$ of the logs of periods with $\|h_m\| \to 0$, such that $\{(\log P_1(q_m),...,\log P_r(q_m)) + t h_m: t \in [0,1)\}$ represent log-period coordinates of saddle connections with the graph $\Gamma$, but $(\log P_1(q_m),...,\log P_r(q_m)) + h_m$ does not. Now, if the sequence $\{q_m\}$ has no vanishing clusters, this follows easily from the fact that $\log P_i$ is a system of period coordinates. So now we induct on the number of nested clusters.\\

\noindent First, we note that the edges of $\Gamma$ contain a spanning tree for each vanishing cluster of $\{q_m\}$, i.e. there are $n-1$ edges that are internal to each vanishing cluster with $n$ singularities, since otherwise there would be a trivial improvement of $\Gamma$ by adding an edge that joined two singularities in the same vanishing cluster that were not connected by a path in $\Gamma$ that does not go outside the cluster, and deleting some other edge. In fact, the edge set of $\Gamma$ is naturally in bijection with the disjoint union of the edge sets of length-minimizing trees of saddle connections on the quadratic differentials corresponding the the various cluster differentials $(X_\ell,q_\ell)$ associated to the limit of the sequence $\{q_m\}$ in the compactification of \ref{DefCompactification}. To obtain the saddle connections corresponding to a length-minimizing tree in $(X_\ell,q_\ell)$, simply delete all singularities outside $S_\ell$ and contract all edges belonging to each $S_j \subsetneq S_\ell$; the remaining graph will have edges corresponding to $S_\ell$.\\

\noindent We will show that for large enough $m$, each singularity is bounded away from each edge of $\Gamma$ of which it is not a vertex along $\gamma_m([0,1)) := \{(\log P_1(q_m),...,\log P_r(q_m)) + t h_m: t \in [0,1)\}$, contradicting the fact that the coordinate system does not extend the point where $t = 1$.\\

\noindent STEP 1: We observe the following: for any two singularities $a$ and $b$, $d_{q_m}(a,b)$ and the distance between $a$ and $b$ along the metric graph $\Gamma$ is uniformly comparable as $m \to \infty$, i.e. $d_{q_m}(a,b) \dot{\asymp} \left. d_{q_m} \right|_\Gamma(a,b)$. The implied constants in $\dot{\asymp}$ depend on the sequence but not on $m$. Moreover, both of these distances are comparable to $(\dot{\asymp})$ the $q_m$-diameter of the smallest vanishing cluster $S_\ell(a,b)$ containing $a$ and $b$, or the entire singularity set if no such vanishing cluster exists, and comparable to the length of any saddle connection that is internal to $S_\ell(a,b)$ but not to a proper vanishing subcluster.\\

\noindent Proof: This is clear from the definitions and proposition \ref{PropVanishInEither}, since any saddle connection whose limit is less than every positive multiple of $\diam_q(S_\ell(a,b))$ would belong to a proper subcluster.\\

\noindent STEP 2: Let $d_{q_m + t h_m}$ denote the metric at time $t$ along $\gamma_m$. For any two singularities $a$ and $b$, $\frac{d}{dt}\log d_{q_m + t h_m}(a,b)$ is defined for almost every $t$, and converges to $0$ uniformly along points of definition as $m \to \infty$. Moreover, $\log d_{q_m + t h_m}(a,b)$ is locally absolutely continuous, so the integral of this derivative represents the change in length.\\

\noindent Proof: Since the angle formed by two saddle connections varies real-analytically, the set of times at which the configuration of saddle connections on the geodesic from $a$ to $b$ changes is discrete, consisting only of isolated times at which two of the saddle connections form an angle of $\pi$. This establishes the final claim, so we must now bound the logarithmic derivatives of distances. This is done by induction on the nesting of the cluster containing $a$ and $b$.\\

\noindent First, we note that the distance from $a$ to $b$ is the sum of the lengths of finitely many saddle connections, and a saddle connection from $c$ to $d$ has period equal to the sum of the periods of saddle connections of the path from $c$ to $d$ in $\Gamma$ with appropriate signs. There are therefore only finitely many ways to express the periods of saddle connections comprising the geodesic from $a$ to $b$ as linear combinations of periods from $\Gamma$, and since the geodesic does not revisit an edge, there is a bound on how many saddle connections are in the geodesic. So the derivatives (but not necessarily logarithmic derivatives) of distances are all uniformly bounded and tend to $0$. In particular, this means that for each $\delta > 0$, every maximal vanishing cluster remains a $\delta$-cluster along $\gamma_m$ for almost every $m$.\\

\noindent Now, by induction, assume that some vanishing cluster $S_\ell$ is, for each $\delta > 0$ a $\delta$-cluster in the singularity set of each point in $\gamma_m$ for all sufficiently large $m$, and moreover, assume that for all proper subclusters $S_k \supsetneq S_\ell$, the conclusion holds for any singularities that belong to $S_k$ but not to a common subcluster of $S_k$. Then along all but finitely many $\gamma_m$, the unique geodesic joining each pair of singularities in $S_k$ consists entirely of saddle connections that are internal to $S_k$. We can then rescale each $q_m$ to assume $S_k$ has $q_m$-diameter $1$, and ignore all singularities outside of $S_k$ and by the same arguments, the conclusion holds for all pairs of singularities in $S_k$ that are not in any proper vanishing subcluster. By induction the conclusion holds for all vertices $a,b$ of $\Gamma$.\\

\noindent STEP 3: By step 2 no singularities collide along $\gamma_m$ as $t \to 1$, $\gamma_m$ may be extended to all $[0,1]$, but with the possibility that some edges of $\Gamma$ may degenerate to concatenations of saddle connections rather than saddle connections at $\gamma(1)$. Moreover, for each pair of singularities $a$ and $b$, the direction of the geodesic from $a$ to $b$ varies a.e. differentiably, locally absolutely continuously, and with derivative uniformly tending to $0$ on points of definition as $m \to 0$. In particular, for each saddle connection that exists along $\gamma$, the angle changes absolutely continuously, and differentiably with uniformly bounded derivative.\\

\noindent Proof: It is sufficient to show that the imaginary part of the logarithmic derivative of the period of every saddle connection goes to $0$ as $\gamma_m \to \infty$. If $a$ and $b$ are joined by a saddle connection at some point along $\gamma_m$, and $m$ is sufficiently large, let $S_\ell(a,b)$ be the minimal vanishing cluster containing $a$ and $b$. For all points $q_m + t h_m$ on $\gamma$ we have $d_{q_m + t h_m}(a,b) \dot{\asymp} \diam_{q_m}(S_\ell(a,b))$. The period of the saddle connection joining $a$ and $b$ is a bounded linear combination of periods of saddle connections of $\Gamma$ internal to $S_\ell(a,b)$, all of which are at most a bounded multiple of $d_{q_m + t h_m}(a,b)$, and each of these periods has logarithmic derivative going to $0$ as $m \to \infty$. It follows that the logarithmic derivative of the saddle connection joining $a$ to $b$ goes to $0$.\\

\noindent STEP 4: By STEP 2 no singularities collide along $\gamma_m$ as $t \to 1$ so it must be the case that some the distance from some singularity $a$ to an edge of $\Gamma$ joining two vertices $b,c \neq a$ goes to zero along $\gamma_m$ as $t \to 1$. The only way this can happen is if $\frac{d_{q_m}(a,b)}{d_{q_m}(b,c)} \to 0$ or $\frac{d_{q_m}(a,c)}{d_{q_m}(b,c)} \to 0$.\\

\noindent Proof: The length-minimizing property implies each saddle connection of $\Gamma$ is a side of two equilateral singularity-free triangles in $q_m$. We must have $d(a,b) + d(a,c) \to d(b,c)$ along $\gamma_m$ as $t \to 1$, but since distances are preserved up to a small multiplicative error, it must be the case that $\lim\limits_{m \to \infty} \frac{d_{q_m}(a,b) + d_{q_m}(a,c)}{d_{q_m}(b,c)} \to 1$. Now, if $b,c$ belong to a vanishing cluster not containing $a$ this is impossible, and if the only vanishing clusters containing $a$ and one of $b$ and $c$ contains all three, then $a$ is not on the geodesic joining $b$ and $c$ in the cluster differential associated to this vanishing cluster.\\

\noindent STEP 5: Assume that as $t \to 1$ along infinitely many $\gamma_m$, the singularities $a,b$ belong to a vanishing cluster $S_\ell$, and $b$ is joined a vertex $c \notin S_\ell$ by an edge of $\Gamma$. Then the angle formed at $b$ by the geodesic from $b$ to $c$ and the geodesic from $b$ to $a$ remains bounded away from $0$ along $\gamma_m$.\\

\noindent Proof: This is obvious from step $3$, since the geodesic joining $a$ to $b$ is a concatenation of saddle connections. This completes the proof, since as $t \to 1$, the distance between $a$ and any saddle connection in $\Gamma$ not containing $a$ does not tend to $0$ along $\gamma_m$ as $t \to 1$. $\Box$

\begin{proposition}\label{PropMultStable} Let $K$ be a compact subset of a $\mathcal{M}_{g,n}$ or a moduli space of cluster differentials. Then there is some $\epsilon(K) > 0$ such that for all half-translation surfaces $(X,q,\Sigma)$ length-minimizing systems of saddle connections whose pairs of lifts added with opposite signs form a basis for $H_{\mathrm{odd}}^1(\tilde{X},\tilde{\Sigma})$ persist when the logs of their periods all change by at most $\epsilon(K)$.\end{proposition}

\noindent Proof: The first thing we observe is that given a sequence that converges in the sense of Definition \ref{DefCompactification}, the collection of saddle connections that belong to saddle connections internal to vanishing clusters is either linearly independent or has exactly one dependence relation, which we will describe. Cutting $\tilde{X}$ along this system of saddle connections produces either one connected component or two, and if produces two components, the relation is that the sum of the lifts of the edges that divide length-minimizing tree into two components, each of which contains an odd number of cone points of cone angle an odd multiple of $\pi$, taken with appropriate signs, is the shared boundary of two surfaces. Therefore all but the longest such edge will be included in any length-minimizing system of saddle connections. By \ref{PropMultStable} all of the saddle connections in the length-minimizing system internal to vanishing clusters will persist, since a perturbation that changes the log-periods of all $k$ of the remaining saddle connections by at most $\epsilon$ changes the log-period of the longest one by $O(k\epsilon)$. The proof for that saddle connections that are not internal to vanishing clusters persist is essentially the same as in the cluster differential case, since $(\tilde{X},\tilde{q})$ admits a cover branched only over the poles, and the universal cover of this is complete and non-positively-curved, so there is a unique geodesic joining each pair of points. $\Box$\\

\noindent Proof of Lemma \ref{LemAsympClust}: This is immediate from Lemma \ref{LemOutsiders}, Proposition \ref{PropScale}, and Proposition \ref{PropMultStable}.

\begin{definition}A length minimizing period-coordinate system for $(X,q)$ is \emph{persistent} for $(X^\prime,q^\prime)$ if for each coordinate, the change in the log-period is less than the quantity $\epsilon(K)$ for some compact set $K$ containing $(X,q)$ in its interior.\end{definition}

\noindent Proof of Proposition \ref{PropPersistentCluster}: First, we note that there is an upper bound on the systole of $(X,q)$ and a lower bound on the diameter of $(X,q)$ coming from $K$. We can fix a compact set $K^\prime \subset QD^1(\mathcal{M}_{g,n})$ which contains the projection of $K$ to the moduli space $\mathcal{M}_{g,n}$, defined to be all surfaces $(X,q)$ with some fixed lower bound on the $d_q$-systole and some fixed upper bound on the $d_q$-diameter. For all period coordinate systems on $K^\prime$ with upper bound $L$ on the length of the saddle connections used, as $(X,q)$ varies, the $q$-distances between pairs of singularities, the systole, and the diameter of $(X,q)$ are uniformly Lipschitz with respect to the path metric $d_\mathrm{Euclidean}$, since there are locally finitely many period coordinate systems by \ref{PropLocallyFinite}. Since there is a neighborhood of $K$ which projects to $K^\prime$, and the path $\gamma$ can be assumed to lie entirely in $K^\prime$, the proposition follows easily from the Lipschitz property. $\Box$\\

\noindent Proof of Proposition \ref{PropNoCancellation}: Clearly, the only way a sequence of counterexamples is possible is if $$\frac{d_\mathrm{Euclidean}(X_m,Y_m)}{\epsilon_m} \to 0.$$ Obviously a sequence of counterexamples contains a subsequence in which $\{X_m\}$ converges in the sense of Notation \ref{NotaCompactConverge} by passing to a subsequence. Moreover, we can assume that for every vanishing cluster $D$ of $\{X_m\}$ the ratio $\frac{\diam_{q_m}(D)}{\epsilon_m}$ converges in $[0,\infty].$\\

\noindent Moreover, we can pick a system of period coordinates that consists of a maximal set of saddle connections internal to vanishing clusters whose periods are linearly independent in $H^1(\tilde{X},\tilde{\Sigma};\cc)$, and a length-minimizing set of complementary geodesic arcs in $X_m$ whose endpoints are singularities, that complete the system of period coordinates. Note that the saddle connections defining coordinate system are internal to vanishing clusters if and only if their lengths go to $0$. Then, in these systems of period coordinates, we make the following claim:

\begin{claim}
The complementary geodesic arcs can be chosen (Hausdorff) continuously on the straight line-segment path from $X_m$ to $Y_m$, and the change in their periods is $\dot{\asymp}\epsilon_m$ as $m \to \infty$.
\end{claim}

\noindent Proof of claim: The only periods that change are those whose lengths are bounded away from zero, and the perturbations are close to zero. Thus the length-minimizing coordinate system is persistent. The claim is thus immediate from Proposition \ref{PropMultStable} and the fact that distance in each of the finite permissible coordinate systems is uniformly comparable to distance in the length-minimizing system.\\

\noindent If a sequence of counterexamples exists then we can of course choose them to lie in the principal stratum, since the principal stratum is dense on each coordinate chart.\\

\noindent Now, we assume that we have a counterexample sequence $\{(X_m,Y_m)\}_{m = 1}^\infty$ with $$d_\mathrm{Euclidean}(X_m,Y_m) < \epsilon_m^\prime = o(\epsilon_m),$$ i.e. $\lim\limits_{m \to \infty} \frac{\epsilon_m^\prime}{\epsilon_m} = 0$. Then we will show that pairs of corresponding complementary saddle connections of $X_m$ and $Y_m$ in the persistent coordinate chart have periods differing by $o(\epsilon_m)$.\\

\noindent Along a path $\gamma_m(t)$ from $X_m$ to $Y_m$ that is rectifiable, parametrized by arc length, and of length $\epsilon_m^\prime$ with respect to $d_\mathrm{Euclidean}$. For $1 \leq j \leq 4g - 4 + 2n$, make Hausdorff continuous choices of singularities $a_m^j(t)$ with the property that for each $t$, the collection of cone points $a_m^j(t)$, counted with multiplicity, are the singularities of $\gamma_m(t)$, where $a_m^j(t), 1 \leq j \leq n$ correspond to the marked points and $a_m^j(t), n + 1 \leq j \leq 4g - 4 + n$ correspond to the zeros with multiplicity (or multiplicity - 1 if they collide with marked points). These choices are not unique, but they are unique up to a permutation of the singularities that preserves maximal vanishing clusters, since all non-uniqueness is caused by collisions of singularities, and this can only happen to two singularities in the same vanishing cluster. For a proof that it is possible to continuously choose a root along the path $\gamma$, see e.g. Lemma 1 of \cite{PickRoot}.\\

\noindent Then $t$ runs from $0$ to $\epsilon_m^\prime$, and let $\ell_m(t)$ be a continuous choice of simple (non-self-intersecting) geodesic joining $a_m^j(t)$ and $a_m^k(t)$ such that $\ell_m(0)$ is one of the complementary saddle connections that is part of the length-minimizing system for $X_m$. We note that this simple geodesic is not quite determined by its endpoints and homotopy class because there may be poles, but there are at most two choices for each pair of endpoints, and they can be chosen continuously by passing to an appropriate finite branched cover branched only over marked points and branched over all marked points, if such exist. (In particular, this operation is defined over the base surface, so this finite cover varies in a Lipschitz manner over $QD(\mathcal{T}_{g,n}$.) In fact, we may assume that we picked the functions $a_m^r(t)$ to be projections of continuous choices of singularities in this branched double cover; each such choice determines a unique choice of $\ell_m(t)$. For the remainder of the proof we may assume $\ell_m(t)$ is chosen continuously in the universal cover of some branched cover.\\

\noindent In particular, the length of $\ell_m(t)$ is bounded above and bounded away from $0$. Then, by the definition of the Euclidean path metric, the following are uniformly Lipschitz in $t$, i.e. Lipschitz with a constant that holds for all sufficiently large $m$:

\begin{enumerate}
\item The length of $\ell_m(t)$
\item The distance between any two choices of $\alpha_m^r(t)$, and $\beta_m^r(t)$, if these are two possible ways of choosing $a_m^r(t)$, for each $r$
\item The diameter of the maximal vanishing cluster containing $a_m^r$ along $\gamma_m(t)$, for each $r$
\item The distance from any singularity $a_m^r(t)$ to $\ell_m(t)$
\item The slope of the longest segment of $\ell_m(t)$, measured as an angle in $\rr/\pi\zz.$
\end{enumerate}

\noindent Let $b_m^j$ be the singularity that corresponds to $a_m^j(0)$ in the persistent coordinate system for $X_m$ that extends to $Y_m$. In order for us to get a contradiction, it must be the case that for some choice of $a_j,a_k$ the singularity $a_m^j(\epsilon_m)$ in $Y_m$ is not equal to $b_m^j$, since if the length and angle of $\ell_m(\epsilon_m)$, and hence also the period of the corresponding saddle connection would each differ by $O(\epsilon_m)$ from $\ell_m(0)$. In fact, it must be the case that for some $j$, $\epsilon_m^\prime/d_Y(b_m,a_m^j(\epsilon_m^\prime)) \to 0$. Otherwise, we could move the endpoints of the geodesic joining $a_m^j(t)$ and $a_m^k(t)$ to $b_m^j$ and $b_m^k$ and the length of the geodesic, and the slope of the longest segment, would vary in a manner that is Lipschitz in the change of the endpoints. The resulting geodesic joining $b_m^j$ and $b_m^k$ would be a single saddle connection, whose length and angle are $o(\epsilon_m)$ away from the period of the saddle connection joining $a_m(0)$ and $b_m(0)$.\\

\noindent Finally, we show that $a_m^j(\epsilon_m^\prime)$ and $a_m^k(\epsilon_m^\prime)$ are $O(\epsilon_m^\prime)$ away from the singularities $b_m^j$ and $b_m^k$ in $Y_m$, which makes this contradicion impossible. If we vary the endpoints in $Y_m$ along rectifiable curves, the distance and angle of the longest segment again vary in a uniformly Lipschitz manner with respect to the endpoints. So we can therefore move the endpoints to the singularities corresponding to $a_m$ and $b_m$, and conclude that their lengths and angles have changed by $o(\epsilon_m)$. However, this is also impossible, since a singularity that is not an endpoint of a saddle connection of length $M$ in a length minimizing system and is distance $d$ away from both ends of the saddle connection is distance at least $\min((\sqrt{3}/2)d,M)$ away from the saddle connection. (This is due to the fact that the original saddle connection is a side of two singularity-free equilateral triangles). This gives a contradiction to item 4 above. $\Box$

\bibliographystyle{amsalpha}
\bibliography{lhc}

\end{document}